\documentclass[10pt,a4paper,reqno]{amsart}  
\usepackage{graphicx}
\usepackage{pdfpages}	
\usepackage{amsmath}	
\usepackage{amssymb}
\usepackage{amsthm}
\usepackage{mathtools}
\usepackage{csquotes}
\usepackage[T1]{fontenc}
\usepackage[latin1]{inputenc}
\usepackage[english]{babel}
\usepackage{listings}
\usepackage{graphicx}
\usepackage{times}
\usepackage{tikz}
\usepackage{enumitem}
\usepackage{multicol}
\usepackage{scrextend}
\usepackage{relsize}
\usepackage{tikz-qtree}
 
\usepackage{float}
\usepackage{hyperref}
\usepackage{cite}
\usepackage{tikz-cd}
\usepackage{dsfont}
\usepackage{xcolor}
\usepackage{url}
\usepackage[left=1.2in,right=1.2in,top=1.2in,bottom=1.2in]{geometry}
\numberwithin{equation}{section}

\makeatletter
\@namedef{subjclassname@2020}{%
	\textup{2020} Mathematics Subject Classification}
\makeatother

\newcommand{\NN}{\mathcal{N} \mathcal{N}}
\newcommand{\R}{\mathbb{R}}
\newcommand{\N}{\mathbb{N}}
\newcommand{\T}{\mathcal{T}}
\newcommand{\C}{\mathbb{C}}
\newcommand{\Z}{\mathbb{Z}}
\newcommand{\F}{\mathcal{F}}
\newcommand{\M}{\mathcal{M}}

\newcommand{\x}{v}
\newcommand{\TT}{\T_b}

\newcommand{\strbdepsph}{\epsilon_{\mathbb{P}h\mathrm{B}}^{\mathrm{s}}(\mathrm{p})}
\newcommand{\strbdepsp}{\epsilon_{\mathbb{P}\mathrm{B}}^{\mathrm{s}}(\mathrm{p})}
\newcommand{\nh}{\text{NH}}
\newcommand{\gprob}{\Gamma^{\mathrm{ran}}}
\newcommand{\pr}{\mathbb{P}}

\theoremstyle{plain}
\newtheorem{thm}{Theorem}[section]
\newtheorem{lem}[thm]{Lemma}
\newtheorem{prop}[thm]{Proposition}
\newtheorem{corr}[thm]{Corollary}

\theoremstyle{definition}
\newtheorem{defn}[thm]{Definition}
\newtheorem{exmp}[thm]{Example}
\newtheorem{rem}[thm]{Remark}

\newcommand*{\MG}[1]{#1}
\makeindex

\title[Generalised hardness of approximation and the SCI hierarchy]
{Generalised hardness of approximation and the SCI hierarchy  \\--\\ On determining the boundaries of training algorithms in AI}

\author{Luca Eva Gazdag} 
\address{Department of Mathematics, University of Oslo}
\email{lucaeg@student.matnat.uio.no}

\author{Alexander Bastounis} 
\address{Department of Mathematics, King's College London}
\email{alexander.bastounis@kcl.ac.uk}

\author{Anders C. Hansen} 
\address{Department of Applied Mathematics and Theoretical Physics, University of Cambridge}
\email{a.hansen@damtp.cam.ac.uk}
\newcommand{\catop}{\mathbin\Vert}
\begin{document}
	
	\keywords{Generalised hardness of approximation, phase transitions, Solvability Complexity Index hierarchy, boundaries of AI, foundations of computational mathematics}
	\subjclass[2020]{65Yxx, 03D55 (primary) and 90C26, 15A29, 68Q87, 68W20 (secondary)}

	\maketitle

	\begin{abstract}
		Generalised hardness of approximation (GHA) is the phenomenon that one can easily compute an $\epsilon$-approximation to a  solution of a computational problem for $\epsilon > \epsilon_1 > 0$, but for 
		$\epsilon < \epsilon_1$ (the approximation threshold) it suddenly becomes hard, for example, non-computable or intractable (non-polynomial time).  In this paper we demonstrate the phenomenon that GHA happens \MG{when using AI techniques for solving inverse problems, namely training neural networks (NNs) to optimally perform on the training data}. In particular, for any non-zero underdetermined linear inverse problem the following phase transition can occur: \MG{For a certain} family of training sets $\Omega$, one can prove the existence of optimal NNs for solving the inverse problem for each $\T \in \Omega$, 
		however, these optimal \MG{neural} networks can only be computed to a certain accuracy $\epsilon_1 > 0$. Below the approximation threshold $\epsilon_1$, not only does it become intractable to compute the NNs, it becomes impossible regardless of computing power, and no randomised algorithm can solve the problem with probability better than 1/2. Moreover, despite the existence of a stable optimal NN, any attempts of computing it below two times the approximation threshold $2\epsilon_1$ will yield an unstable NN. 
		Our results use and extend the current mathematical framework of the Solvability Complexity Index 
		(SCI) hierarchy and initiate a program for analysing the GHA phenomenon throughout computational mathematics and AI. GHA generalises the phenomenon of  hardness of approximation in discrete computations to arbitrary computational problems.
	\end{abstract}
	
\vspace{3mm}
	
\noindent Communicated by Shmuel Weinberger.	
	
	\section{Introduction} \label{section:introduction}
	
	In this paper we add to the foundations theory on \emph{generalised hardness of approximation} (GHA) and study the phenomenon in connection with AI methods for underdetermined inverse problems. In particular, we show that the GHA phenomenon occurs when training optimal neural networks for underdetermined inverse problems. GHA is the phenomenon that the feasibility of computing an $\epsilon$-approximation to a problem may change dramatically with the \MG{approximation accuracy} $\epsilon$. In particular, we have an
	\emph{approximate computational problem}: Given an $\epsilon > 0$, the $\epsilon$-approximate computational problem is the problem of computing an approximation that is no more than $\epsilon$ away from the true solution -- in some appropriate predefined metric. 
	Suppose that we have a computational problem and two classes of approximate computational problems $S_1$ and $S_2$ with $S_1 \cap S_2 = \emptyset$ and \MG{thresholds} $\epsilon_1 \geq \epsilon_2 > 0$. For example, we could have 
	\[
	S_1 = P \text{ (polynomial time solvable), } \quad S_2 = P^c \text{ (the complement of $P$)}.
	\]
	We say that the computational problem has an $(S_1,S_2)$-phase transition at $(\epsilon_1, \epsilon_2)$ if we have the following:
	\begin{equation}\label{eq:GenHardAppr}
		\begin{split}
			&\text{The approximate computational problem} \in S_1, \text{ for } \epsilon > \epsilon_1, \\
			&\text{The approximate computational problem} \in S_2, \text{ for } \epsilon < \epsilon_2.
		\end{split}
	\end{equation}
	If $\epsilon_1 = \epsilon_2$ in \eqref{eq:GenHardAppr} we say that the phase transition is sharp and call $\epsilon_1$ the approximation threshold. Schematically, the concept of generalised hardness of approximation with a sharp phase transition can be visualised as follows: 
	\begin{equation}\label{eq:GenHardApprDiag}
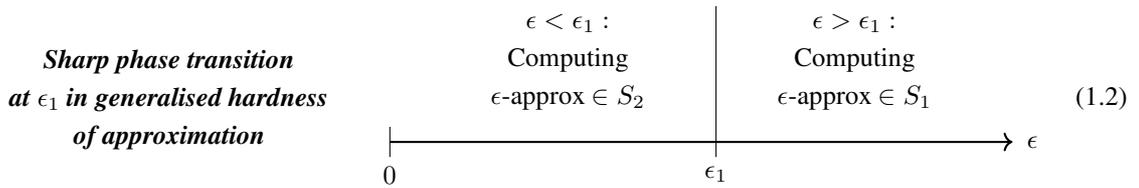
		\begin{minipage}{0.33\textwidth}
			\begin{tabular}{c}
				\bf \emph{Sharp phase transition} \\
				\bf \emph{at $\epsilon_1$ in generalised hardness} \\
				\bf \emph{of approximation} \\
				%\phantom{space}
			\end{tabular}
		\end{minipage}
		\begin{minipage}{0.60\textwidth} 
			\begin{tikzpicture}[xscale=1.3]
				\draw [thick,->] (3,0) -- (9.3,0);
				\draw (3,-.2) -- (3,.2);
				\draw (6.3,-.2) -- (6.3, 1.8);
				\node[align=center, below] at (3,-.2) {$0$};
				\node at (9.5,0) {$\epsilon$};
				\node[align=center, below] at (6.3,-.2) {$\epsilon_1$};
				\node[align=center, above] at (7.7,0.3) 
				{$\epsilon > \epsilon_1:$ \\ Computing \\ $\epsilon$-approx $\in S_1$ };
				\node[align=center, above] at (4.8,0.3)%
				{$\epsilon < \epsilon_1:$ \\ Computing \\ $\epsilon$-approx $\in S_2$ };
			\end{tikzpicture}
		\end{minipage}
	\end{equation}
	This definition can of course be generalised to any family of collections $S_1, \hdots, S_{k}$, $k > 1$ of computational problems with $S_j \cap S_i = \emptyset$ for $j \neq i$ and $\epsilon_1, \hdots, \epsilon_{2(k-1)} > 0$ with $\epsilon_1 \geq \epsilon_2 > \epsilon_3 \geq \epsilon_4 > \hdots > \epsilon_{2k-3}  \geq \epsilon_{2(k-1)}$ as follows. We say that we have an $(S_1, \hdots, S_k)$-phase transition at $(\epsilon_1, \hdots, \epsilon_{2(k-1)})$ if we have the following:
	\begin{equation}\label{eq:GenHardAppr2}
		\begin{split}
			&\text{The approximate computational problem} \in S_1, \text{ for } \epsilon > \epsilon_1,  \\
			&\text{The approximate computational problem} \in S_2, \text{ for } \epsilon_3 < \epsilon < \epsilon_2,\\
			& \qquad \qquad \qquad \qquad \qquad \vdots\\
			&\text{The approximate computational problem} \in S_{k-1}, \text{ for } \epsilon_{2k-3} < \epsilon < \epsilon_{2k-4},\\
			&\text{The approximate computational problem} \in S_k, \text{ for } \epsilon < \epsilon_{2(k-1)},
		\end{split}
	\end{equation}
	where we say that the phase transition in \eqref{eq:GenHardAppr2} is sharp at $\epsilon_{2j-1}$ if $\epsilon_{2j-1} = \epsilon_{2j}$ for some integer $j$ and call $\epsilon_{2j-1}$ an approximation threshold. \MG{Two examples of GHA follow below. See} \S \ref{sec:perspectives} for a more complete list of examples and discussions. 
	
	\begin{exmp}[\underline{Optimisation and Smale's 9th problem with extensions}]
		To the best of our knowledge, the GHA phenomenon in optimisation was first discovered in \cite{comp}, 
		in connection with Smale's 9th problem and its extensions, where the phenomenon was 
		documented in a large collection of convex optimisation problems including linear programmes (LPs). 
		In particular, for LPs, the problem is to compute an element
		\begin{equation}\label{problems}
			z \in \mathop{\mathrm{arg min}}_{x} \langle x , c \rangle \text{ subject to } Ax=y, \quad x \geq 0,
		\end{equation}
		given inexact -- yet arbitrarily fine precision -- representations of the input. 
		In short, \cite{comp} establishes that, for any $0 < \epsilon_2 \leq \epsilon_1$, there are collections of LPs such that for 
		\[
		S_1 = P , \quad S_2 = \Delta_1\setminus P, \quad S_3 = \Delta_1^c
		\]
		these LPs have a sharp $(S_1, S_2, S_3)$-phase transition at $(\epsilon_1, \epsilon_2)$. Here, we have
		\[
		\Delta_1 = \text{ Set of computable problems}, \quad P = \text{Set of polynomial time solvable problems}, 
		\]
		and $\Delta_1^c$ denotes the complement of $\Delta_1$, see \S \ref{section:preliminaries} for a review of the classes (such as $\Delta_1$) in the Solvability Complexity Index (SCI) hierarchy.
		This work was expanded in \cite{Matt2}, showing that the GHA phenomenon occurs when computing certain neural networks through optimisation problems, see also \cite{Choi2}. 
		Determining all the potential GHA classes $(S_1, \hdots, S_k)$ and the approximation thresholds 
		$(\epsilon_1, \hdots, \epsilon_{2(k-1)})$ for different classes of LPs is a delicate open problem (see Problem 5 (J. Lagarias) in \cite{AIM}).
	\end{exmp}
	\begin{exmp}[\underline{Hardness of approximation in computer science}]
		GHA generalises the concept of hardness of approximation (HA) \cite{Sudan, Arora_JACM_98_2, Lovasz_JACM_96, Arora2007, Sudan_Overview_2009, Hastad_Acta} in computer science to general settings in computational mathematics. The much celebrated PCP theorem \cite{Arora2007, Lovasz_JACM_96, Arora1_Goedel, Arora_JACM_98_2} implies that there are large collections of discrete combinatorial optimisation problems for which there is a sharp $(S_1, S_2)$-phase transition at $\epsilon_1 \geq 0$ (where $\epsilon_1$ depends on the problem), where 
		\[
		S_1 = P, \quad S_2 = P^c\MG{.}
		\]
		The issue of strict inequality $\epsilon_1 > 0$ typically depends on the P vs. NP question. 
	\end{exmp}
	
	In this paper we study the GHA phenomenon in connection with AI methods for 
	underdetermined inverse problems. More specifically, 
	we show that the GHA phenomenon occurs when training optimal neural networks for underdetermined 
	linear inverse problems.

	{\it Acknowledgements.} LG acknowledges support from the Niels Henrik Abel and C. M. Guldbergs  memorial fund. ACH acknowledges support from the Simons Foundation Award No. 663281 granted to the Institute of Mathematics of the Polish Academy of Sciences for the years 2021-2023, from a Royal Society University Research Fellowship, and from the Leverhulme Prize 2017.
	The authors want to thank the anonymous referee for their highly detailed and useful reviews with very helpful suggestions that have substantially improved the manuscript.
	
\section{Optimality of neural networks for underdetermined linear systems}
	
	We study the following underdetermined systems of equations. Let $A: \R^N \to \R^m$ be a linear mapping with non-trivial kernel and let $\mathcal{M}_1 \subset \mathbb{R}^N$ be some subset that we call initial domain. We now consider the following inverse problem:
	\begin{align} \label{problem}
		\text{Given measurements} \; y = Ax + e \;  \text{of} \; x \in \mathcal{M}_1 , \; \text{recover} \;  x \in \mathcal{M}_1,
	\end{align}
	where $e \in \mathbb{R}^m$ is a potential noise vector \MG{(for our purposes, we consider the case where $e =0$ since this is generally an easier problem to compute and thus our negative results become stronger)}. 
	These types of problems have been extensively studied in sparse recovery and compressed sensing when 
	$\mathcal{M}_1$ is, for example, a collection of sparse vectors or vectors with some structured sparsity 
	\cite{Juditsky_2012, Juditsky_2011, AdcockHansenBook, candes2006robust, donoho2006compressed, adcock2016generalized, Boyer_2016, BOYER_ACHA_2019, Breaking, Bast_2017, Bast_SIAM_NEWS, fannjiang2020numerics, candes2013phaselift}. However, recent developments have led to a plethora of  
	AI techniques \cite{Genzel,jin17, mccann2017convolutional, arridge2019solving, 
		Baraniuk2020, hammernik2018learning} for solving these types of problems. More specifically, deep learning approaches have become popular over the last years, where one trains a neural network 
	$\mathbf{N}: \mathbb{R}^m \rightarrow \mathbb{R}^N$ from some training set $\mathcal{T} \subset \R^N \times \R^m$ of the form 
	$\mathcal{T} = \{(x^j,Ax^j)\}_{j=1}^{\MG{\ell}}$, for some $\{x^j\}_{j=1}^{\MG{\ell}} \subset \mathcal{M}_1$. 
	However, it has been established that there is a fundamental stability-accuracy trade-off for such methods \cite{antun2020instabilities,gottschling2023troublesome, Matt2}. Indeed, great accuracy on certain inputs may cause unstable behaviour and hallucinations in the reconstruction in the form of false information in other reconstructed objects. Thus, it becomes important to establish the optimal choice of neural network in order to optimise performance. Such optimal maps have already been a focal point in approximation theory. 
	
	Indeed, approximation theory has a rich tradition in the theory of optimal approximations and optimal reconstruction maps. A particular example is the seminal work of A. Cohen, W. Dahmen and R. DeVore \cite{DeVore}, where they define the concept of optimal reconstruction maps for underdetermined inverse problems. 
	
	\begin{defn}[Optimality in the sense of Cohen, Dahmen \& DeVore \cite{DeVore}] \label{optmap:trainingset}
		Let $A : \R^N \rightarrow \R^m$ be linear, $\mathcal{M}_1 \subset \R^N$ and 
		\[
		\mathcal{M}_2 := A(\mathcal{M}_1).
		\]
		Define the optimality constant for the pair $(A,\mathcal{M}_1)$ as
		\begin{equation*}
			c_{\mathrm{opt}}(A,\mathcal{M}_1)
			= \inf_{\varphi \colon\! \mathcal{M}_2  \rightrightarrows \R^N} \sup_{x \in \mathcal{M}_1} d_1^{H}(\varphi(Ax), x),
		\end{equation*}
		where $d_1^H$ denotes the Hausdorff metric associated with a metric $d_1$ \MG{(in particular, for a non-empty set $M \subset \R^N$ we have $d_1^{H}(M,x) = \sup_{\tilde x \in M} \|\tilde x - x\|_{\MG{2}}$ and if $\varphi$ is single-valued then $d_1^{H}(\varphi(Ax),x) = d_1(\varphi(Ax),x)$)}. Here, the double arrow notation $\rightrightarrows$ denotes that the mapping can be multivalued. 
		We define a family of approximately optimal maps of $(A,\mathcal{M}_1)$ as follows. We say that 
		$\varphi_{\epsilon} : \mathcal{M}_2 \rightrightarrows \R^N$ is a family of approximate optimal maps of $(A,\mathcal{M}_1)$ if for all $\epsilon \in (0,1]$, 
		\begin{equation}\label{eq:approx_optimal}
			\sup_{x \in \M_1}  \, d^H_1(\varphi_{\epsilon}(Ax),x) \leq c_{\mathrm{opt}}(A,\mathcal{M}_1) + \epsilon,
		\end{equation}
		that $\varphi_{\epsilon}$ is $\epsilon$-optimal, and that $\varphi_0$ is an optimal map for $(A,\mathcal{M}_1)$ if $\varphi_0$ satisfies \eqref{eq:approx_optimal} with 
		$\epsilon = 0$.
	\end{defn}

	The first key question is whether there exist neural networks that are optimal maps -- or $\epsilon$-optimal -- for different underdetermined inverse problems. 
	Moreover, one can ask if 
	such an optimal neural network can be trained from training data. In particular, given $A$ as above and 
	$\ell \in \N$, the class $T_{\ell}(A)$ of all training 
	sets with $\ell$ elements is defined as
	\begin{align}\label{eq:T_ell}
		T_{\ell}(A) := \{ \T \in\MG{ (\R^N \times \R^m)^\ell} \, \vert \, & \T = ((x^k,y^k))_{k = 1}^{\ell} \text{ with }
		\|x^k\|_{\MG{2}}, \|y^k\|_{\MG{2}} \leq 1 \: \text{and } \: \MG{y^k = Ax^k},\\& (x^k,y^k) \neq (x^j,y^k) \text{ for every } j,k \in \{1,2,\dotsc,\ell\} \text{ with } j \neq k\}.\notag
	\end{align}

	We \MG{impose} the norm restriction to the elements in $T_{\ell}(A)$ above, because we are interested in 
	studying the errors of algorithms relative to the magnitude of the elements in the training sets.
	
	\MG{We slightly abuse notation by using standard notation for sets when the meaning is clear from context. Some examples include the following: for $\T  = ((x^k,y^k))_{k = 1}^{\ell} \in \T_{\ell}(A)$ we write $(x,y) \in \T$ whenever there exists a $j \in \{1,2,\dotsc,\ell\}$ with $(x,y) =  (x^j,y^j)$. Similarly, we write} 
	\MG{\[[\T \setminus {(x^j,y^j)}] \cup {(x',y')}:= ((x^1,y^1),(x^2,y^2),\dotsc (x^{j-1},y^{j-1}),(x',y'),(x^{j+1},y^{j+1}),\dotsc,(x^\ell,y^\ell))\]}
	\MG{and, for functions $f$ with inputs in $\R^N \times \R^m$, we write \[
	f(\T):= (f(x^1,y^1),f(x^2,y^2),\dotsc,f(x^\ell,y^\ell)) \text{ when }\T = ((x^1,y^1),(x^2,y^2),\dotsc,(x^\ell,y^\ell)).\]}
	\MG{Particular functions of tuples that we will frequently use are the projections $\pi_1(\T)$ and $\pi_2(\T)$. For $\T$ as above we set \begin{align}
			\pi_1(\T) &:= (x^1,x^2,x^3,\dotsc,x^\ell) \label{eq:pi1def}\\
			\pi_2(\T) &:= (y^1,y^2,y^3,\dotsc,y^\ell).\label{eq:pi2def}\end{align}
	}
	\MG{We use the standard notation $\catop$ to mean the \emph{concatenation} of two tuples so that, in particular, 
	\[\T \catop ((w,z)) := ((x^1,y^1),(x^2,y^2)\dotsc,(x^\ell,y^\ell),(w,z)), \quad (w,z) \in \R^N \times \R^m,
	\]
	and similarly for arbitrary tuples.}
	\MG{Finally, we define $|\T|$ as the number of distinct elements in a tuple $\T$, analogously to set notation. }
	
	Furthermore, let $\mathcal{F}$ be a collection of underdetermined inverse problems $(A,\mathcal{M}_1)$ and let $\NN_{m,N}^{\F}$ be the set of NNs (of arbitrary size) mapping $\R^m \to \R^N$ that are bounded on 
	$\bigcup_{(A, \M_1) \in \mathcal{F}} A(\M_1),$ where we will throughout the paper use the more compact notation
	\begin{equation}\label{eq:Cups}
	\bigcup_{\mathcal{F}} \M_2 := \bigcup_{(A, \M_1) \in \mathcal{F}} A(\M_1), \qquad \bigcap_{\mathcal{F}} \M_2 := \bigcap_{(A, \M_1) \in \mathcal{F}} A(\M_1).
	\end{equation}
	The boundedness assumption on $\bigcup_{\mathcal{F}} \M_2$ \MG{for $\NN_{m,N}^{\F}$}  \MG{is added to ensure that the metric we define later (see \eqref{eq:Xi_norm}, \eqref{metric} and \eqref{breakdown_metric})  on $\NN_{m,N}^{\F}$ is well-defined.}
	We also require that the NNs in $\NN_{m, N}^{\F}$ have a finite fixed set of nonlinearities that are computable (see Definition \ref{def:NNs_with_rhos}).
	\MG{\begin{rem}
			Equation \eqref{eq:T_ell} ensures that we work on bounded training data; furthermore, single-valued computable functions are necessarily continuous (and hence neural networks of interest in this paper are continuous). Thus, the boundedness assumption on  $\bigcup_{\mathcal{F}} \M_2$ for $\NN_{m,N}^{\F}$ is a very natural one since continuous functions on compact domains are bounded.
	\end{rem}}
	 We then have the following questions: Let $A \in \R^{m \times N}$ be given, let $\Omega \subset \T_{\ell}(A)$
	be a collection
	of training sets, and assume that there is some prescribed way of associating to each training sample 
	$\T \in \Omega$ an initial domain $\M_1(\T) \subset \R^N$. 
	\vspace{1mm}
	\begin{displayquote}
		{\it 
			\begin{itemize}
				\item[(i)] Does there exist, for each $\T \in \Omega$ and each $\epsilon \: \geq\:  0$, a neural network that is $\epsilon$-optimal for the inverse problem $(A, \M_1(\T))$?
				\item[(ii)] Does there exist an algorithm $\Gamma$ that for each 
				$\T \in \Omega$ and $\epsilon > 0$ produces a neural network 
				$\Phi_{\T,\epsilon} = \Gamma(\T, \epsilon)$ that is $\epsilon$-close to an optimal neural network for the inverse problem $(A, \M_1(\T))$?
			\end{itemize}
		}
	\end{displayquote}
	\vspace{1mm}
	Note that Question (ii) is only sensical if the answer to Question (i) is `yes' with $\epsilon = 0$. However, this paper is mostly focusing on answering Question (ii), given an affirmative answer to Question (i) with $\epsilon = 0$. In this paper we throughout consider the case where 
	\begin{align} \label{M1_dep}
		\M_1(\T) := \{x \, \vert \, (x,y) \in \T\}.
	\end{align}
	That is, we study the problem of computing neural networks that have optimal \MG{performance} on the training set $\T$ for each $\T \in \Omega$ \MG{(note the similarity to $\pi_1(\T)$)}. The motivation for this is to get the strongest possible negative results. In particular, if the answer to Question (ii) is negative for the choice of \MG{$\M_1(\T)$ as in \eqref{M1_dep}}, any hope of computing an optimal NN on more general choices of $\M_1(\T)$ will be gone \MG{(under the weak assumption that each $(x,y) \in \T$ has $x \in \M_1(\T)$ -- in some sense, the training set is relevant to the inverse problem)}. We will also use the notation
	\[
	\mathcal{M}_2(\mathcal{T}) \coloneqq A(\M_1(\T)).
	\]

	\begin{rem}[The basic mechanism: Why the answer to Question (ii) is "no"]
		The proofs of our GHA results rely on some delicate positive results (`there exists an algorithm') and negative results (`no algorithm exists'). However, all the negative results have a \MG{core idea} (combined with Proposition \ref{computability_prop}) that we sketch here. 
		For simplicity, \MG{consider a non-zero $A$ with non-trivial nullspace and training sets with two elements. Take} non-trivial vectors $w \in (\text{ker } A)^\bot$ and $v \in \text{ker } A$. Consider the training sets
		\begin{align*}
			\T_n = \left\{ \left(0,0\right),\left(v+w/n,A(v+w/n)\right) \right\}, \quad 
			\quad  \T^* = \{(0,0),(v,0)\}.
		\end{align*}
		Then it is not hard to show that any optimal network $\Phi^*$ for $\T^*$ (i.e., any network minimizing the training error on the set $\T^*$) needs to satisfy $\Phi^*(0) = v/2$, whereas any optimal network $\Phi_n$ for $\T_n$ needs to satisfy $\Phi_n(0) = 0$.
		However, since $v + \frac{w}{n} \xrightarrow[n \to \infty]{} v$ and 
		$A(v + \frac{w}{n}) = \frac{1}{n}A(w) \xrightarrow[n \to \infty]{} 0$, it follows that the training set $\T_n$ comes arbitrarily close to the training set $\T^*$ as $n \to \infty$, so that an algorithm based on 
		\emph{inexact inputs} cannot distinguish between $\T_n$ and the set $\T^*$ \MG{for large $n$}. This simple phenomenon gives rise to the many delicate GHA issues described in Theorem \ref{thm_intro:second} and Theorem \ref{thm_intro:first}, and the negative results are universal in the sense that they hold for all randomised general algorithms (see Definitions \ref{def:General_alg} and \ref{definition:ProbablisticAlgorithm}) that are much more powerful than any randomised Blum-Shub-Smale (BSS) machine \cite{Complexity_and_real_comp} or any randomised Turing machine \cite{Turing_Machine}. 
		We give an intuitive explanation for why we need to consider inexact inputs and what 
		we mean by inexact inputs in the next \MG{subsection}.
	\end{rem}
	
	\subsection{The model of computation}\label{sec:model_comp}
	
	Fix a linear map $A: \R^N \to \R^m$, and consider a collection $\Omega$ of training sets $\mathcal{T}$,
	\[
	\mathcal{F} = \bigcup_{\T \in \Omega} \MG{\{}(A,\mathcal{M}_1(\T))\MG{\}},
	\]
	which is the 
	corresponding collection of inverse problems, and a mapping 
	\begin{equation}\label{eq:Xi}
		\Xi: \Omega \rightrightarrows  \NN_{m,N}^{\F}, \text{  such that  } \,\,  
		\Xi(\mathcal{T})  = \left\{ \mathbf{N}_{opt}^{\M_1(\T)} \: \vert \: \mathbf{N}_{opt}^{\M_1(\T)} 
		\: \text{ is optimal for } (A,\mathcal{M}_1(\T)) \right\}.
	\end{equation}
	Here, the double arrow notation $\rightrightarrows$ denotes that a map can be multivalued. \MG{This choice is made because the optimal neural network may not be unique and thus the map $\Xi$ can be multivalued.}  We will in some cases be interested in bounding the Lipschitz-constant of the optimal mappings to ensure stability. Thus, we define 
	$\Xi_{D}: \Omega \rightrightarrows \NN_{m,N}^{\F}$, with $D>0$, 
	by
	\begin{equation}\label{eq:Xi_D} 
		\Xi_{D}(\mathcal{T}) \coloneqq \left\{ \mathbf{N}_{opt}^{\M_1(\T)} \in  
		\Xi(\mathcal{T}) \; \vert \; L(\mathbf{N}_{opt}^{\M_1(\T)}) \leq D \right\},
	\end{equation}
	where, for a function $N: \R^m \to \R^N$, the Lipschitz constant $L(N)$ is given by
	\begin{align} \label{Lipdef}
		L(N) \coloneqq \sup\limits_{y \neq \hat y}
		\frac{\| N(y) - N(\hat y)\|_{\MG{2}}}{\|y -\hat y\|_{\MG{2}}}.
	\end{align}
	Note that $\Xi_{D}(\mathcal{T})$ could potentially be empty for small values of $D$. 
	
	Our goal is to compute 
	$\mathbf{N}_{opt}^{\M_1(\T)}$ -- or more precisely, one of the optimal neural networks.  
	We take the word `compute' literally, and thus we wish to work with models which process approximate inputs. Indeed, often a training set $\mathcal{T} = \MG{(}(x^j,y^j)\MG{)}_{j=1}^{\MG{\ell}}$ will not be \MG{exactly representable} on a computer. This is because $A$ could have rows from the discrete Fourier transform, for example -- as in accelerated (subsampled) Magnetic Resonance Imaging (MRI) -- and thus $A$ contains irrational numbers.  Another issue is that an overwhelming amount of modern software used is based on floating-point arithmetic, and hence even if the input is rational, there will be inexactness due to the floating-point representation. For example, $1/3$ can only be approximated in finite base-2 arithmetic, 
	giving rise to round-off approximation.
	
	\MG{We model this inexactness in the following way. Let $\tilde \T$ be $\ell$ ordered sequences of $2$-tuples, that is \[\tilde \T = (\{(\tilde x^1(n),\tilde y^1(n))\}_{n=1}^{\infty},\{(\tilde x^2(n),\tilde y^2(n))\}_{n=1}^{\infty},\{(\tilde x^3(n),\tilde y^3(n))\}_{n=1}^{\infty},\dotsc, \{(\tilde x^\ell(n),\tilde y^\ell(n))\}_{n=1}^{\infty}).\]
	\MG{We say that $\tilde \T$ \emph{corresponds} to $\T \in T_{\ell}(A)$ if for each $j \in \{1,2,\dotsc,\ell\}$ and $n \in \mathbb{N}$ we have}
	\begin{equation}\label{eq:approx}
		\|\tilde x^{j}(n) - x^j\|_{\MG{2}} \leq 2^{-n},\qquad \|\tilde y^{j}(n) - y^j\|_{\MG{2}} \leq 2^{-n}.
	\end{equation}
	Informally, the algorithm is not allowed to work with the exact input $\T$ but instead can only  access $(\tilde x^j(n),\tilde y^j(n))$ for finitely many $n$ (that must be chosen in a (potentially) adaptive but recursive way -- see \S\ref{inexact} and \S\ref{comp_problem} for details) and every $j \in \{1,2,\dotsc,\ell\}$ via e.g. an oracle. We require that a successful algorithm should work \MG{correctly} on any $\tilde T$ which corresponds to $\T$.
	}
	
	Given an input set $\Omega \subseteq T_{\ell}(A)$ consisting of training sets, we denote by 
	$\tilde \Omega$, the corresponding input set consisting of approximated training sets, more specifically
\MG{\begin{align}
		\tilde \Omega \coloneqq \{ \tilde \T \, \vert \, \exists \, \T \in \Omega, \; \tilde \T \text{ corresponds to } \T\}.
	\end{align}
In fact, we will use the notation $\tilde \T$ throughout to refer to an inexact training set above that corresponds to a training set $\T$.}
	Note that this extended computational model of having inexact input is standard and can be found in many areas of the mathematical literature, and we mention only a small subset here including the work in \cite{bishop1967foundations, Cucker_Smale97, Fefferman_Klartag, Fefferman_Klartag2, Ko1991ComplexityTO, lovasz1987algorithmic}.

	\begin{rem}[Oracles]\label{rem:norms_for_Xi}
		The above model means that for each training set $\mathcal{T} = \MG{(}(x^j,y^j)\MG{)}_{j=1}^{\MG{\ell}} \in \Omega$, there \MG{are infinitely many possible $\tilde \T$ which correspond to $\T$.}
		A sequence of approximations is provided to the algorithm through an `oracle'. For example, in the case of a Turing machine \cite{Turing_Machine}, this would be through an oracle input tape (see \cite{Ko1991ComplexityTO} for the standard setup), or in the case of a Blum-Shub-Smale (BSS) machine 
		\cite{Complexity_and_real_comp}, this would be through an oracle node. The algorithm can thus ask for an approximation to any given accuracy as in 
		\eqref{eq:approx}, and use as many queries as desired. 
	\end{rem} 
	
	\begin{rem}[Computability of the `oracle'] \label{oracle1}
		We want to emphasise that our results continue to hold when instead of requiring an algorithm to ``work correctly'' on every $\tilde \T$ corresponding to $\T \in \Omega$, we consider inexact inputs $\tilde \T$ which are computable, i.e. when 
		$\tilde \T = \{(\tilde x^{j}\MG{(n)},\tilde y^{j}\MG{(n)}\}_{n\in \mathbb{N}}$ is a computable sequence for $j = 1,\hdots, \MG{\ell}$. See also Remark \ref{rem:TurVsMarkov} and Remark \ref{remark:CompDelta1Markov}. 
	\end{rem}
	
	\begin{rem}[Markov model -- When the `oracle' is the input as a finite string] 
		\label{oracle2}
		A Markov algorithm \cite{MarkovModel} cannot -- as opposed to a Turing machine -- take an infinite string as an oracle input. However, if an infinite string represents a computable number, where the elements in the string are approximations to the computable number, the code or algorithm producing the string can be used as a finite input. For example, in the Markov model, suppose that $\mathcal{T} = \MG{(}(x^j,y^j)\MG{)}_{j=1}^{\MG{\ell}}$ contains the number $\pi$ in some of the vector components. The input to the Turing machine (or Markov algorithm) would be algorithms -- represented by finite strings -- that would provide arbitrarily fine precision approximations to $\pi$. See Remark \ref{rem:Markov} \MG{regarding our results in the Markov model}.
\end{rem}

	\begin{defn}[Computing the neural network $\mathbf{N}_{opt}^{\M_1(\T)}$ to $\epsilon$-accuracy]
		We say that the mapping $\Xi$ in \eqref{eq:Xi} can be computed to $\epsilon$-accuracy if there exists an algorithm $\Gamma$ such that for any 
		$\mathcal{T} \in \Omega$,
		\begin{equation}\label{eq:Xi_norm}
			\inf_{\mathbf{N}_{opt}^{\M_1(\T)} \in \Xi(\mathcal{T}) } \sup_{y \in \bigcup_{\F} \M_2} \| \Gamma(\tilde{\T}, \epsilon)(y) - \mathbf{N}_{opt}^{\M_1(\T)}(y) \|_{\MG{2}} \leq \epsilon, \quad \forall \,  \mathcal{\tilde T} \MG{\in \tilde \Omega \text{ so that } \tilde \T \text{ corresponds to } \T}. 
		\end{equation} 
	\end{defn}
	
	\begin{rem}[Slight abuse of notation for $\Gamma(\mathcal{T}, \epsilon)$]
		As the algorithm must work \MG{correctly} on any such representation $\tilde{\T}$ of $\mathcal{T}$ we will frequently abuse notation and write $\Gamma(\mathcal{T}, \epsilon)$ 
		instead of $\Gamma(\tilde{\mathcal{T}}, \epsilon)$. Also, if a result holds for only one specific $\epsilon \in \R_+$, we will frequently drop writing $\epsilon$ as an input to $\Gamma$, writing $\Gamma(\T)$ in place of $\Gamma(\T, \epsilon)$.
	\end{rem}

	%%%%%%%%%%%%%%%

	\section{Main results -- Generalized hardness of approximation (GHA)} \label{main_results_section}
	
	Our main results are gathered in two main theorems, both presented formally in \S \ref{section:proofs}. 
	The first main theorem asserts that there exist phase transitions for a large class of inverse 
	problems for the computational problem 
	described above, in the sense that for any $\epsilon_1 \in (0, 3/16]$, there exists classes of 
	training sets $\Omega$ such that
	\begin{itemize}
		\item[(i)] One cannot compute an optimal neural network for all $\T \in \Omega$ for any accuracy 
		\MG{better than}\footnote{Throughout, we use the wording `better than' and `worse than' when discussing accuracies to refer to accuracy smaller than a certain value (better than) or accuracy larger than a certain value (worse than).} $\epsilon_1$.
		\item[(ii)] The Lipschitz constants of the computed networks will not be well behaved \MG{for any accuracy that lies} between 
		$\epsilon_1$ and $2\epsilon_1$.	
		\item[(iii)] We can compute an optimal neural network for all $\T \in \Omega$, with uniformly bounded Lipschitz constants, for any accuracy \MG{worse than} $2\epsilon_1$.
	\end{itemize}
	
	The second main theorem
	asserts fundamental computational barriers when attempting to train neural networks to solve 
	\eqref{problem} in the standard computational model with inexact input.
	In addition, the second theorem illustrates how small changes in 
	the training set can lead to the collapse of the accuracy of a working algorithm. In particular, the theorem demonstrates how phase transitions  in 
	\eqref{eq:GenHardAppr} and \eqref{eq:GenHardApprDiag} can suddenly change with the training data.
	
	For all impossibility results we use a generalized model for computation -- that is also used in 
	\cite{comp,Matt1,Matt2,Hansen_JAMS, Ben_Artzi2022, ben2021computing} -- in order to obtain universal lower bounds regardless of the computational model. More precisely, when we refer to an algorithm we mean a so called \emph{general algorithm} (the formalism of this is defined in \S \ref{algorithm}). However, for all positive results -- of the form "there exists an algorithm" -- we use the Turing model to achieve the strongest upper bounds possible.

	\begin{thm}[Generalized hardness of approximation -- Phase transitions for computing optimal NNs]\label{thm_intro:second}
		For any integers $N > m$, and any $\beta > 0$, 
		consider any fixed non-zero linear map 
		$A: \R^N \to \R^m$ such that the spectrum 
		$\mathrm{Sp}(AA^*) \subset [\beta^2,\infty)$. Then, for any rational 
		$\epsilon_1 \in (0, 3/16]$ and any integer \MG{$\ell \: > \; m +2$}, there exists a domain 
		$\Omega \subset T_{\ell}(A)$ (as described in \eqref{omega}) of training sets and a set of corresponding initial domains 
		$\{\M_1(\T) : \T \in \Omega\}$ (as described in \eqref{M1_dep}), such that the following occurs.
		For the mapping
		$\Xi:\Omega \rightrightarrows \NN_{m,N}^{\F}$ (as described in \eqref{eq:Xi}), we have that
		$\Xi(\mathcal{T}) \neq \emptyset$, for each $\T \in \Omega$.
		Moreover, each of the following happen simultaneously:
		\begin{itemize}
			\item[(i)]
			No algorithm, not even randomised, that takes any \MG{$\mathcal{\tilde T} \in \tilde \Omega$ (corresponding to $\mathcal{T} \in \Omega$)} as inputs, can approximate an optimal neural network 
			$\mathbf{N}_{opt}^{\M_1(\T)} \in \Xi(\mathcal{T})$ for all 
			$\T \in \Omega$ to accuracy $\epsilon_1$ (with probability greater than $p > 1/2$ in the randomised case -- this is even the case if the algorithm has a non-zero probability of not halting).
			\item[(ii)] There exists an algorithm $\Gamma$ such that for all 
			$\tilde \T \in \tilde \Omega$ and dyadic
			$\epsilon > 2 \epsilon_1$ we have that $\mathbf{N}_{\tilde \T,\epsilon} = 
			\Gamma(\tilde \T,\epsilon)$  is a NN that approximates an optimal neural network 
			$\mathbf{N}_{opt}^{\M_1(\T)} \in \Xi(\mathcal{T})$ \MG{(where $\tilde \T$ corresponds to $\T \in \Omega$)} to accuracy 
			$\epsilon$ \MG{(in the sense of \eqref{eq:Xi_norm})} and is such that the Lipschitz constant $L(\mathbf{N}_{\tilde \T,\epsilon})$
			is uniformly bounded by $2/\beta$.
			\item[(iii)]  However, there exists a $\T_1 \in \Omega$ and a $D>0$ such that $\Xi_{D}(\T_1) \neq \emptyset$ (recall \eqref{eq:Xi_D}), that is, there exists an optimal neural network $N_{opt}^{\M_1(\T_1)}$ 
			for $\T_1$ such that
			the Lipschitz constant $L(N_{opt}^{\M_1(\T_1)}) \leq D$. Yet, there is a $\tilde \T_1 \in \tilde \Omega$ \MG{corresponding to $\T_1$}, such that for any $K >0$, $\delta \in (0, \epsilon_1)$ and any algorithm $\hat \Gamma$ such that 
			$\mathbf{N}_{\tilde \T,\epsilon} = \hat \Gamma(\tilde \T,\epsilon)$  is a NN that approximates an optimal NN 
			$\mathbf{N}_{opt}^{\M_1(\T)} \in \Xi(\mathcal{T})$ to accuracy 
			$\epsilon \in (\epsilon_1, 2\epsilon_1 - \delta]$ \MG{(in the sense of \eqref{eq:Xi_norm})} \MG{for all $\T \in \Omega$ and $\tilde \T$ corresponding to $\T$}, we have the following: 
			The vector $0 \in \pi_2(\T_1)$ but for every $\eta > 0$ there exists a $y \neq 0$ with $\|y\|_{\MG{2}} \leq \eta$ and
			\begin{align}\label{eq:Kom}
				\frac{\|\mathbf{N}_{\mathcal{\tilde T}_1,\epsilon}(0) - 
					\mathbf{N}_{\mathcal{\tilde T}_1,\epsilon}(y)
					\|_{\MG{2}}}{\|y\|_{\MG{2}}} > K.
			\end{align}
		\end{itemize}
	\end{thm}
	
	See Theorem \ref{thm:second_formal} for the technical statement of Theorem \ref{thm_intro:second}. 	
	
	\begin{rem}[Stronger result\MG{s} in \MG{(i) and} (ii) of Theorem \ref{thm_intro:second}] \MG{Our result in part (i) is stronger than stated in that the failure of any algorithm is in a sharper sense than \eqref{eq:Xi_norm}; see also Remark \ref{rem:whereFail}.} Note \MG{also} that we prove that the algorithm in (ii) of Theorem \ref{thm_intro:second} is actually independent of the ordering of the training sets $\T$. That is, it will work \MG{correctly} on any reordering of the training sets, and the algorithm is agnostic to the reordering. 
	\end{rem}

	\begin{rem}[The stable NN exists but cannot be computed] Part (iii) of Theorem \ref{thm_intro:second} says that \MG{a stable optimal NN} exists, but it cannot be computed to accuracy \MG{better than} $2\epsilon_1$. Yet, according to (ii), \MG{a stable optimal NN} can be computed to accuracy \MG{worse than} $2\epsilon_1$. Moreover, the instability of the attempted computed NN in (iii) will occur on the training set. In particular, the blow up of the Lipschitz constant happens `at a point' in the training data, as described in \eqref{eq:Kom}. 
	\end{rem}

	\begin{rem}[No $p=2/3$ randomised algorithm with non-zero probability of not halting]\label{rem:2/3}
		Note that the phase transition\MG{s} described above are different to the phase transitions described in 
		\cite[Theorem 3.3]{comp}. Indeed, the phase transitions \MG{discussed} in part (ii) of Theorem 3.3 in 
		\cite{comp} show constructively the existence (in that setting) of randomised algorithms with a non-zero probability of not halting which succeed with probability $p=2/3$. This is not the case for the generalised hardness of approximation phenomenon described in Theorem \ref{thm_intro:second}. Hence, the phenomenon described in Theorem 
		\ref{thm_intro:second} is `harder' and different to the generalised hardness of approximation phenomenon first discovered in \cite{comp}. Furthermore, \cite{comp} does not address the neural network setting of this paper. 
	\end{rem}

	\begin{rem}[Consequences of Theorem \ref{thm_intro:second} -- GHA implies an accuracy-stability trade-off] 
		The accuracy-stability trade-off in AI methods for inverse problems is well documented empirically \cite{Genzel, antun2020instabilities, gottschling2023troublesome, Matt2} and to some extent theoretically \cite{gottschling2023troublesome, Matt2}. Theorem \ref{thm_intro:second} provides a new theoretical understanding of this phenomenon, and demonstrates that any attempt \MG{to compute NNs too accurately} will necessarily yield arbitrarily unstable NNs. However, as long as \MG{one ensures that} the accuracy of the computed NN is \MG{worse than} a certain threshold, stability can be achieved. Thus, there is an accuracy-stability trade-off and overperformance necessarily yields instabilities.
	\end{rem}
	\begin{thm}[Phase transitions -- Rapid changes and the Halting problem] \label{thm_intro:first}
		Given any $N,m \in \mathbb{N}$ \MG{(with $N \geq 3$)}, a fixed non-zero linear map $A: \R^N \to \R^m$ with non-trivial kernel and an integer $\ell \geq 2$, 
		there exists a domain 
		$\Omega \subset T_{\ell}(A)$ (as described in \eqref{omega}) of training sets $\T$, a set of corresponding initial domains 
		$\{\M_1(\T) \, \vert \, \T \in \Omega\}$ (as described in \eqref{M1_dep}) such that the following occur.
		For the mapping $\Xi:\Omega \rightrightarrows \NN_{m,N}^{\F}$ (as described in \eqref{eq:Xi}), we have $\Xi(\mathcal{T}) \neq \emptyset$ for each $\T \in \Omega$. However, the following happens simultaneously:
		\begin{itemize}[leftmargin=7mm]
			\item[$(i)$] The mapping $\Xi$ cannot be computed. In particular, no (randomised) algorithm that 
			takes any \MG{$\mathcal{\tilde T} \in \tilde \Omega$ (corresponding to $\mathcal{T} \in \Omega$)} as inputs, can produce a neural network that approximates 
			any $\mathrm{N}_{opt}^{\M_1(\T)} \in \Xi(\T)$ to an accuracy of 0.1 (see \eqref{eq:1}) for all $\mathcal{T} \in \Omega$ 
			(with probability $p > 1/2$, see \eqref{eq:2}).

			\item[(ii)] 
			There exists an infinite sequence of \MG{pairwise distinct} training sets 
			$\{\T^1_i\}_{i \in \N} \subset \Omega$ and a \MG{(potentially non-halting)} algorithm $\Gamma^1$ that computes\footnote{\MG{Here `computes' is understood in the sense of computing to an arbitrary prescribed accuracy, i.e. $\Gamma^1(\cdot,\epsilon)$ produces a NN that approximates an optimal NN with accuracy (in the $2$-norm, cf. \eqref{metric}) better than $\epsilon$, where $\epsilon > 0$ is an arbitrary rational number.}} an optimal NN in $\Xi(\T^1_i)$ for each \MG{$i \in \mathbb{N}$ and for each inexact input $\tilde \T^1_i$ (corresponding to $\T^1_i$)}. However, for any infinite sequence 
			$\{ \T^2_{k} \}_{k \in \N} \subset \Omega$ of \MG{pairwise distinct} training sets and any algorithm $\Gamma^2$ such that: 
			\begin{enumerate} 
				\item $\Gamma^2$ takes inputs in $\tilde \Omega$;
				\item for any $k \in \N$, $\Gamma^2$ produces  -- on input $\tilde \T^2_{k}$ (corresponding to $\T^2_{k}$) -- an approximation to an optimal NN in 
			$\Xi(\T^2_{k})$ \MG{to accuracy (in the $2$-norm, cf. \eqref{metric}) better than 0.1};
			\end{enumerate} we have the following.
			For each $j \in \mathbb{N}$, there exists an integer $i$ and $(x,y) \in \T^2_i$, as well as an element 
			$(x', y') \in \R^N \times \R^m$ with 
			$\| x'\|_{\MG{2}}, 
			\| y' \|_{\MG{2}} \leq 1$ and 
			$\| (x,y) - (x',y')\|_{\MG{2}} \leq \sqrt{2}/4^{j}$, 
			such that if we replace $(x,y)$ with $( x', y')$ then we obtain a new training set 
			$\T^{\prime} = [\T^2_i \setminus \{(x,y)\}] \cup \{( x',  y')\} \in \Omega$ 
			such that 
			\begin{align}\label{eq:failure222}
				\sup_{y \in \M_2(\T^{\prime})}\| \Gamma^2(\tilde \T^{\prime})(y) - \mathbf{N}_{opt}^{\mathcal{M}_1(\T^{\prime})}(y)  \|_{\MG{2}}  > 10^{-1}, 
			\end{align} 
			for some $\tilde \T^{\prime} \in \tilde \Omega$ (corresponding to $\T^{\prime}$),
			where $\mathbf{N}_{opt}^{\mathcal{M}_1(\T^{\prime})}$ is any optimal neural network for the inverse problem $(A,\M_1(\T^{\prime}))$.  In particular, the failure in \eqref{eq:failure222} applies to $\Gamma^1$ and $\{\T^1_i\}_{i \in \N}$. 
			\item[(iii)] Consider the arbitrary sequence $S = \{ \T^2_{k} \}_{k \in \N}$ and algorithm $\Gamma^2$ from (ii). Then, there is a $\T^2  = \{(x^{\prime,l},y^{\prime,l})\}_{l=1}^{\ell} \in \Omega$ such that for any $j$, we have 
			$
			\max_{1 \leq l \leq \ell}\| (x^l,y^l) - (x^{\prime,l},y^{\prime,l})\|_{\MG{2}} \leq \sqrt{2}/4^{j}, 
			$
			where $\{(x^l,y^l)\}_{l=1}^{\ell} = \T^2_n$,
			for some $n = n(j) \in \mathbb{N}$, and so that \MG{even if $\Gamma^2$ could compute optimal NNs (to arbitrary accuracy)} for training data in $S$, it cannot \MG{do so} for training data in $S \cup \{\T^2\}$. In \MG{fact,  }
			\[
			\sup_{y \in \M_2(\T^2)}\| \Gamma^2(\tilde \T^2)(y) - \mathbf{N}_{opt}^{\mathcal{M}_1(\T^2)}(y)  \|_{\MG{2}}  > 10^{-1}, \, \text{ for some } \tilde \T^2 \in \tilde \Omega \MG{\text{ corresponding to } \T^2},
			\] 
			where $\mathbf{N}_{opt}^{\mathcal{M}_1(\T^2)}$ is any optimal neural network for the inverse problem $(A,\M_1(\T^2))$.	
			
			\item[(iv)] If \MG{every entry of} $A$ is computable, the problem of computing approximations to $\Xi:\Omega \rightrightarrows \NN_{m,N}^{\F}$, to an accuracy of 0.1, is comparable to the Halting problem and is at least as hard.
		\end{itemize}
	\end{thm} 
	
	See Theorem \ref{thm:first:technical} for the technical statement of Theorem \ref{thm_intro:first}.

	\begin{rem}
		Note that by (i) in Theorem \ref{thm_intro:first} it follows that the algorithm $\Gamma^1$ in (ii) must fail on some inputs in $\tilde \Omega$. The failure of the algorithm $\Gamma^1$ we construct is that \MG{it may not halt on some inputs}. 
	\end{rem}

	\begin{rem}
		Note that the algorithms constructed in Theorem \ref{thm_intro:second} and
		Theorem \ref{thm_intro:first} (the positive results) are implementable on a Turing machine that is \emph{given an oracle for the inexact inputs} (that is, an oracle for the approximations of the elements in $\Omega$). This oracle does not need to be computable (implementable on a Turing machine). 
		However, the negative results hold when the oracle is computable, and even when the algorithm can access the code for the oracle (specifically in the Markov model, see also Remarks \ref{rem:norms_for_Xi}, \ref{oracle1} and \ref{oracle2}).
	\end{rem}
	
	\subsection{Consequences of Theorem \ref{thm_intro:first} and the SCI hierarchy}  For a quick review of the SCI hierarchy, see \S \ref{section:preliminaries}, and \cite{Hansen_JAMS, SCI, CRAS, Matt1, colbrook_spec_meas} for detailed discussions.

	\begin{rem}[Phase transitions can change rapidly with training data]
	Part (ii) of Theorem \ref{thm_intro:first} shows that replacing elements in the training set can change the phase transition and the approximation threshold rapidly. Moreover, the training data that causes the change in the phase transition can be arbitrarily close to elements that are already in the original training set. Hence, the phase transitions can be highly unstable and thus the training process itself can be unstable.
	\end{rem}
	
	\begin{rem}[Our results and the Halting problem]
		The non-computability of the Halting problem asserts that there is no Turing machine that can take any algorithm (described as a finite string) with an input (also described as a finite string) such that the Turing machine will halt and determine 'yes' if the algorithm would halt on the input, and 'no' otherwise. Likewise, our GHA results say that if \MG{the approximation threshold $\epsilon$ satisfies $\epsilon < \epsilon_1$, then} there does not exist a Turing machine (or Markov algorithm), that takes as input any algorithm (described as a finite string) that computes a computable input $\T$ to the problem function $\Xi$, and produces an 
		$\epsilon$-approximation to $\Xi(\T)$. Hence, our results and the Halting problem are indeed very similar, however, 
		it is not immediate that they are comparable, see Remark \ref{halting_comparable} for further details.
	\end{rem}

	\begin{rem}[Training NNs beyond the approximation threshold is as hard as the Halting problem] 
		\label{halting_comparable}
		The classical Halting problem is often used to show non-computability results -- in the cases where a comparison is possible. It is somewhat of a tradition to determine if a given problem is comparable to the Halting problem. A recent example is the problem of deciding the spectral gap in quantum mechanics \cite{Cubitt1, Cubitt2}, which is comparable to the Halting problem.  Note that it is not always the case that a comparison can be made. In fact, Turing reduction only provides a partial ordering on decision problems, hence \MG{sometimes} two problems may not be comparable.  Thus, at a first glance, it may be somewhat surprising that the problem of computing NNs for inverse problems is comparable to the Halting problem -- a problem that has little to do with computational problems in analysis and AI. Indeed, the Halting problem is in $\Sigma_1\setminus \Delta_1$ in the SCI hierarchy. However, as the metric used to analyse the problem of computing optimal NNs for inverse problems does not allow for any concept of convergence 'from above and below' -- which is needed to define the $\Sigma_j$ and $\Pi_j$ classes -- only the $\Delta_j$ classes make sense for computing NNs. Hence,   the possible comparison of the two problems is not obvious. 
		The fact that our problem of interest is comparable to the Halting problem raises two issues:
		\begin{itemize}
			\item[(I)] Showing (i) in Theorem \ref{thm_intro:first}, in the Turing/Markov model, can be done with more classical tools, without any reference to the key proposition from the SCI framework -- since we establish a comparison to the Halting problem.  However, (i) in Theorem \ref{thm_intro:first} is actually much stronger, as it holds in any model of computation, including the more general BSS model. The more general the computational model, the harder it is to show non-computability. Also, the probabilistic result cannot be proven through a comparison with the Halting problem. For general computational models and the probabilistic results, Proposition \ref{computability_prop} -- which is a main driver for GHA results -- is needed. 
			\item[(II)] Since the problem of computing the NN in Theorem \ref{thm_intro:first} is as hard as the Halting problem (which is in $\Sigma_1\setminus \Delta_1$), the question of where the problem is in the SCI hierarchy is immediate. Moreover, one can ask where arbitrary problems of computing optimal NNs for different inverse problems lie in the SCI hierarchy:
			\begin{displayquote}
				\vspace{1mm}
				{\it Where in the SCI hierarchy are the problems of computing optimal neural networks for inverse problems? }
				\vspace{1mm}
			\end{displayquote}
			Note that the SCI classification theory of spectral problems is very rich \cite{Ben_Artzi2022,ben2021computing,colbrook_spec_meas,Matt1,Hansen2016ComplexityII,CRAS,SCI,Hansen_JAMS}. It remains to \MG{be seen} whether computing decoders for inverse problems may have an equally delicate classification theory. 
		\end{itemize}
	\end{rem}

	\section{GHA in mathematics and connections to our results}\label{sec:perspectives}
	
	The mathematics behind GHA in the various areas of the mathematical sciences is in its infancy, as this phenomenon is a rather new area in computational mathematics and scientific computing. However, we provide a summary below helping to put our results in perspective. We have already mentioned optimisation and HA in \S \ref{section:introduction}, and below follows a complementing list of \MG{areas of mathematics}. 
	\begin{itemize}[leftmargin=8pt]
		\item[] \underline{\emph{Compressed sensing and statistics}}:  A natural question to ask in scientific computing is the following: Given a subclass $\Omega$ of computational problems, does the GHA phenomenon 
		occur? And if so, what are the classification classes $(S_1, \hdots, S_k)$ and the approximation thresholds $(\epsilon_1, \hdots, \epsilon_{2(k-1)})$. As an example, consider the classical compressed sensing problem \cite{AdcockHansenBook} of computing an approximation to an element 
		\begin{equation}\label{problems3}
			z \in \Xi(y,A) = \mathop{\mathrm{arg min}}_{x} \|x\|_1 \text{ subject to } \|Ax - y\|_2 \leq \delta, \qquad \delta \in [ 0,1],
		\end{equation}
		where $\Omega$ consists of \MG{all} matrices $A \in \mathbb{R}^{m \times N}$ that satisfy the \emph{$\ell_2$-robust nullspace property of order $s$} with parameters $\rho \in (0,1)$ and $\tau>0$, $y = Ax$ \MG{for a vector $x$ which is $s$-sparse, where} $m$ depends on $N$ and $s \in \mathbb{N}$ (see \cite{AdcockHansenBook, comp} for details). 
		\begin{thm}[Theorem 7.1 in \cite{comp}]
			\label{thrm:BP}
			(Paraphrased) Consider the setup above with the given $\Omega$ and with 
			$\delta \in [0,(1-\rho)/16\tau]$ and further assume that $\rho \in (1/3,1)$ and $\tau > 10$. Then the basis pursuit problem \eqref{problems3} has a
			$(S_1,S_2)$-phase transition at $(\epsilon_1, \epsilon_2)$ with $\epsilon_2 = \delta/2$ and $\epsilon_1 = 16\tau \delta(1-\rho)^{-1}$, where $S_1 = P$ and $S_2 = \Delta_1^c$.		
		\end{thm}

		\item[] The GHA phenomenon also \MG{occurs} in statistics, in particular in the LASSO problem for feature selection. In \cite{Alex1, Alex2}, this was investigated in relation to a newly defined condition number for the LASSO problem. In this special case, the problem of computing the condition number becomes a $\Sigma_1$ problem in the SCI hierarchy. As a consequence, despite the GHA phenomenon, one can compute the feature selection problem in statistics, when the condition number is finite.  \\

		\item[] \underline{\emph{Computer assisted proofs and blow up of 3D Euler}}:  GHA becomes crucial in computer assisted proofs, see \cite{AIM}. Indeed, although it may seem counter-intuitive at a first glance, non-computable problems can be used in computer assisted proofs, where famous examples include the Dirac-Schwinger conjecture \cite{fefferman1994_2, fefferman1996interval} and the Kepler conjecture \cite{hales2005proof, hales2017formal} (see \cite{SCI} for details).  A third example is the newly announced computer assisted proof, by Chen and Hou \cite{Hou2,Hou3}, of the century-long open question on blow-up of the 3D Euler equation with smooth initial data. In connection with this, T. Hou poses, in Problem 2 in \cite{AIM}, the question: ``{\it In problems where mathematical analysis precedes a computer-assisted step: How can we ensure that the formulation of the problem is correctly posed so that computability and non-computability of a problem can be determined?}" This problem turns out to be a problem in GHA. Consider the 3D Euler equation 
		\begin{equation}\label{eq:Euler}
			\boldsymbol{u}_t + (\boldsymbol{u} \cdot \nabla) \boldsymbol{u} = -\nabla p, \quad \nabla \cdot \boldsymbol{u} = 0,
		\end{equation}
		with smooth initial data. The delicate issue is that a blow-up of the solution to \eqref{eq:Euler} implies instability in terms of unboundedness \cite{Vasseur2019BlowUpST} of the forward operator taking the initial data to the solution at a given time. One might think, as suggested in \cite{Vasseur2019BlowUpST}, that this means that 3D Euler blow-up cannot be determined through computational means, as this complication could hinder the prospects of a computer-assisted proof, where the validity of the computational step needs to be verifiable. Indeed, unboundedness of the solution operator typically yields non-computability of the PDE solution \cite{SCI}. 
		However, the key theorem (Theorem 4 in \cite{Hou2}) -- that implies the blow-up of the 3D Euler equation -- becomes a theorem in GHA, as it establishes an upper bound of the approximation threshold $\epsilon_1 \leq 10^{-3}$ for the problem of computing the solution to a rescaled 3D Euler equation (as given in (6.13) of \cite{Hou2}) for all time and for a specific set of initial values. It is the fact that $\epsilon_1 \leq 10^{-3}$, which allows for the deduction of the blow-up of the original 3D Euler equation (not rescaled), finally resolving the problem.

		Note that the fact that the proof establishes the upper bound needed on the approximation threshold ends the speculations in \cite{Vasseur2019BlowUpST} that non-computability would hinder a computer assisted proof of 3D Euler. Indeed, GHA helps the understanding of how potential non-computability can be dealt with in computer assisted proofs.
	\end{itemize}
	
	\begin{rem}[The mathematical phenomena behind GHA] The mathematical mechanisms behind GHA are very diverse. For example, the results in \cite{comp, Matt2, Choi2} are due to multivaluedness of the function to be approximated, and these problems are easier than the Halting problem \cite{CRP}. However, in the case of differential equations, as in the 3D Euler case \cite{Hou2,Hou3, Vasseur2019BlowUpST}, the mechanism is often unboundedness/discontinuity. This is also the case in the results in our paper, where the function to be approximated is multivalued, yet there is also a discontinuity phenomenon, and the problems can be compared to the Halting problem. A different phenomenon occurs in feature selection in statistics \cite{Alex1, Alex2}. Here, multivaluedness and discontinuities both occur, but in a different way compared to the results in our paper. In classical hardness of approximation, the mechanism is different and depends on the P vs NP question. The different mechanisms cause some highly delicate intricacies regarding the use of randomised algorithms. 
	\end{rem}
	
	In the language of GHA -- just as in the case of Theorem \ref{thrm:BP} from compressed sensing -- parts of Theorem \ref{thm_intro:second} can be stated (informally) as follows.
	\begin{corr}[Parts of Theorem \ref{thm_intro:second} (informally) in the GHA language]
		Consider the setup in Theorem \ref{thm_intro:second}. Then, for the problem of computing the mapping 
		$\Xi:\Omega \to \NN_{m,N}^{\F}$, there is a $(S_1,S_2)$-phase transition at $(2\epsilon_1, \epsilon_1)$, where $S_1 = \Delta_1$ and $S_2 = \Delta_1^c$.  
	\end{corr}
	To describe the phase transition at $2\epsilon_1$, when considering the problem of computing $\Xi_D$, one needs to include getting a bounded Lipschitz constant as a part of the computational problem. We avoid the details here, as the purpose of this section is to provide an informal introduction to GHA.

	\section{Historical background, the mathematics of the SCI \MG{hierarchy} and related work}\label{sec:history}
	
	The results in this paper can be viewed as a continuation of Smale's program 
	\cite{Complexity_and_real_comp, smale_question, Smale2} 
	on the foundations of computational mathematics. Smale \MG{posed} several 
	fundamental questions on the foundations of computations, among them his 18th problem: what are the limits of artificial intelligence? -- which bares similarities with the famous Turing paper from 1950 
	\cite{Turing}. Our work can be viewed as a step towards answering this question. There are several results \MG{in the literature} that are very much related to \MG{the present} paper.
	
	\begin{itemize}[leftmargin=12pt]
		\item[] {\bf  \emph{The mathematics behind the SCI hierarchy:}} GHA is part of the greater program on the mathematics behind the SCI hierarchy, and this foundations program provides the framework for our results and proofs. The SCI framework was introduced in \cite{Hansen_JAMS} and continued in the work by J. Ben-Artzi et al. \cite{Ben_Artzi2022, ben2021computing}, in the work by M. Colbrook et al.  
		\cite{colbrook_spec_meas, Matt1} as well as in the work by O. Nevanlinna 
		\cite{Hansen2016ComplexityII, CRAS, SCI} and co-authors. See also the work by S. Olver and M. Webb \cite{webb2021spectra}, and \cite{colbrook2021computingSIREV}. The SCI hierarchy is directly related to S. Smale's \cite{Smale2, Smale_Acta_Numerica} program on the foundations of computational mathematics and scientific computing that initiated 
		the early work by C. McMullen \cite{McMullen1, McMullen2, Smale_McMullen} and P. Doyle \& C. McMullen 
		\cite{Doyle_McMullen} on polynomial root-finding. These are pioneering classification results in the SCI hierarchy that were obtained prior to the introduction of the SCI hierarchy. See also classification results in the SCI hierarchy by S. Weinberger \cite{Weinberger}. Note that some of the mathematics behind the SCI hierarchy can be traced back to K. G\"odel \cite{godel1931formal} and A. Turing \cite{Turing_Machine}, however, the new techniques developed allow for any model of computation.
		
		\vspace{1mm}
		
		\item[] {\bf  \emph{Instability in AI:}} Our results are intimately linked to the instability phenomenon in AI methods -- which is widespread \cite{heaven2019deep, finlayson2019adversarial, Hallucination_NatureM2, Choi, SzZ-14, DezFaFr-17} -- and our results add theoretical understandings to this vast research program. There are also particular links to the work by B. Adcock and N. Dexter \cite{adcock2020gap}, P. Grohs and F. Voigtlaender \cite{Felix}, and V. Antun et al.\cite{antun2020instabilities}. See also recent developments by D. Higham, I. Tyukin et al. 
		\cite{THGW21, tyukin2020adversarial}. 
		
		\vspace{1mm}
		
		\item[] {\bf  \emph{Existence vs computability of NNs:}} There is a substantial literature on existence results of NNs \cite{yarotsky2018optimal, boelcskei2019optimal, petersen2018optimal}, see for example the review papers by A. Pinkus \cite{pinkus} and the work by R. DeVore, B. Hanin, and G. Petrova \cite{devore2020neural} and the references therein.  However, as established in \cite{Matt2} by M. Colbrook, V. Antun et al.\@, only a small subset of the NNs that can be proven to exist can be computed by algorithms. However, following the framework of A. Chambolle and T. Pock \cite{Chambolle_Alg, Chambolle_Alg2}, the results in \cite{Matt2} demonstrate how -- under specific assumptions -- stable and accurate NNs can be computed. See also the work by P. Niyogi, S. Smale and S. Weinberger \cite{Smale_Weinberger} on existence results of algorithms for learning.
		
		\vspace{1mm}
		
		\item[] {\bf  \emph{The PCP theorem:}} The 2001 G\"odel Prize was awarded to S. Arora, U. Feige, S. Goldwasser, C. Lund, L. Lov\'{a}sz, R. Motwani, S. Safra, M. Sudan, and M. Szegedy for their work on the much celebrated PCP theorem \cite{Arora2007, Lovasz_JACM_96, Arora1_Goedel, Arora_JACM_98_2} and its connection to hardness of approximation. The PCP theorem implies that -- subject to P$\neq$NP -- there are large collections of combinatorial optimisation problems for which there is a sharp phase transition at some $\epsilon_1 > 0$ (where $\epsilon_1$ depends on the problem). 
	\end{itemize}

	\section{Tools for the proofs -- Mathematical preliminaries from the SCI hierarchy} \label{section:preliminaries}
	The SCI hierarchy and the mathematical framework that \MG{supports} it have been very useful in order to \MG{chart} the boundaries of computational mathematics in the theory of spectral problems, inverse problems, optimisation, AI etc. The SCI hierarchy is based on the concept of a computational problem that we formally define below. This is described by a (potentially multi-valued) function 
	\[
	\Xi:\Omega \rightrightarrows \mathcal{M}
	\] that we want to compute, where $\Omega$ is some domain, and $(\mathcal{M},d)$ is a metric space. 
	In this section we recall the basic concepts of the SCI-hierarchy and we follow the lines of \cite{comp} in our presentation.

	\begin{defn} (\cite[Definition 8.2]{comp})\label{def:GenCompProblem}.
		Let $\Omega$ be a set, which we call the domain. Let $\Lambda$ be a set of complex valued functions 
		$f:\Omega \to \C$ such that for $\iota_1,\iota_2 \in \Omega$, $\iota_1 = \iota_2$ if and only if 
		$f(\iota_1) = f(\iota_2)$ for all $f \in \Lambda$, called an evaluation set. Let $(\M,d)$ be a metric 
		space, and finally let $\Xi: \Omega \rightrightarrows \M$ be a (potentially multi-valued) function which we call the problem function. 
		We call the collection $\{ \Xi, \Omega, \M, \Lambda \}$ a computational problem.
	\end{defn}
	\subsection{Algorithms} \label{algorithm}
	
	Whenever we aim to use neural networks to solve a computational problem, there is a process of constructing an, in some sense, optimal neural network. This is 
	generally called the training process. We formalize this process in terms of an algorithm. 
	
	\begin{defn}(General algorithm, \cite[Definition 8.3]{comp}). \label{def:General_alg}
		Given a computational problem $\{ \Xi, \Omega, \M, \Lambda \}$, a general algorithm is a mapping 
		$\Gamma: \Omega \to \M \cup \{\text{NH}\}$, with NH $\notin \M$, such that, for each $\iota \in \Omega$, the following conditions hold:
		\begin{itemize}
			\item[$(i)$] There exists a subset of evaluations $\Lambda_{\Gamma}(\iota) \subseteq \Lambda$,
			and, whenever $\Gamma(\iota) \neq  \text{NH}$, we have $|\Lambda_\Gamma(\iota)|<\infty$,
			\item[$(ii)$] The action of $\Gamma$ on $\iota$ is uniquely determined by 
			$\{f(\iota) \}_{f \in \Lambda_{\Gamma}(\iota)}$,
			\item[$(iii)$] for every $\iota' \in \Omega$ such that $f(\iota') = f(\iota)$ for every 
			$f \in \Lambda_{\Gamma}(\iota)$, it holds that $\Lambda_{\Gamma}(\iota') = \Lambda_{\Gamma}(\iota)$.
		\end{itemize}
	\end{defn}
	
	\begin{rem}(The purpose of \MG{the notion of} a general algorithm: universal impossibility results, \cite[Remark 8.4]{comp}).
		The purpose of a general algorithm is to have a definition that will encompass any model of computation, and that will allow impossibility results to become universal. Given that there are several non-equivalent models of computation, \MG{our impossibility results will be shown with respect to this general definition of an algorithm.}
	\end{rem}
	Owing to the presence of the special non-halting ''output'' NH, we have to extend the metric 
	$d_{\mathcal{M}}$ on $\mathcal{M}\times \mathcal{M}$ to 
	$d_{\mathcal{M}}:(\mathcal{M} \cup \{\text{NH}\}) \times 
	(\mathcal{M} \cup \{\text{NH}\}) \to [0,\infty]$ in the following way:
	\begin{equation}\label{eq:extended-metric}
		d_{\mathcal{M}}(x,y) = \begin{cases} d_{\mathcal{M}}(x,y) & \text{ if } x,y \in \mathcal{M} \\
			0 & \text{ if } x = y = \text{NH}\\
			\infty & \text{ otherwise.} \end{cases}
	\end{equation}
	In the case of $\eqref{eq:extended-metric}$ $d_{\mathcal{M}}$ is, strictly speaking, not a metric 
	any more, since it can attain the value $\infty$, however, for simplicity we will abuse terminology slightly and keep referring to the extended function $d_{\mathcal{M}}$ as a metric or an extended metric.
	
	There are a myriad of different types of machines that can be used to model an algorithm: the Turing machine \cite{Turing_Machine} (and its cousins including the Markov model \cite{MarkovModel}), the BSS machine \cite{BSS_machine}, the von Neumann architecture  \cite{von_Neumann}, the real RAM 
	\cite{realRAM}, etc. as well as their randomised versions. Indeed, since randomised methods, such as for example randomised gradient descent, are often used when training neural networks, we need to consider randomised algorithms in order to achieve full generality. However, these models are not equivalent when it comes to computability. Thus, to create universal impossibility results we use general algorithms from 
	Definition \ref{def:General_alg} and randomised general algorithms from Definition \ref{definition:ProbablisticAlgorithm} that encompass any reasonable definition of a computational model in the way that they are more powerful than any standard machine, therefore making the impossibility results stronger. Formally we define a randomised algorithm as follows:

	\begin{defn}(Randomised General Algorithm, \cite[Definition 8.25]{comp}).
		\label{definition:ProbablisticAlgorithm} 	Given a computational problem $\{\Xi,\Omega,\mathcal{M},\Lambda\}$, where $\Lambda = \{f_k \, \vert \, k \in \mathbb{N}, \, k\leq |\Lambda| \}$, a \emph{randomised general algorithm} (RGA) is a collection $X$ of general algorithms $\Gamma:\Omega\to\mathcal{M}\cup \{\nh\}$, a sigma-algebra $\mathcal{F}$ on $X$, and a family of probability measures $\{\mathbb{P}_{\iota}\}_{\iota \in \Omega}$ 
		over the $\sigma$-algebra $\mathcal{F}$ so that the following conditions hold:
		\begin{enumerate}[label=(\roman*)]
			\item For each $\iota \in \Omega$, the mapping $\gprob_{\iota}:(X,\mathcal{F}) \to (\mathcal{M}\cup \{\nh\}, \mathcal{B})$ defined by $\gprob_{\iota}(\Gamma) = \Gamma(\iota)$ \MG{is measurable}, 
			where $\mathcal{B}$ is the sigma-algebra on $\mathcal{M}\cup \{\nh\}$ given by
			
			\begin{align*}
				\mathcal{B} \coloneqq \{ T \: : \: T \subset \M \: \text{Borel}\} \cup 
				\{T \cup \{\text{NH}\} \:  : \: T \subset \M \: \text{Borel}\}. 
			\end{align*} \label{property:PAlgorithmMeasurable}
			\vspace{-5mm}
			\item For each $n \in \mathbb{N}$ and $\iota \in \Omega$, we have $\lbrace \Gamma \in X \, \vert \, T_{\Gamma}(\iota) \leq n \rbrace \in \mathcal{F}$, where 
			\begin{align*}
				T_{\Gamma}(\iota) := \sup\{ m \in \N \: | \: f_m \in \Lambda_{\Gamma}(\iota) \}
			\end{align*}
			is the minimum amount of input information.
			\label{property:PAlgorithmRTMeasurable}
			
			\item For all $\iota_1,\iota_2 \in \Omega$ and $E \in \mathcal{F}$ so that, for every $\Gamma \in E$ and every $f \in \Lambda_{\Gamma}(\iota_1)$, we have $f(\iota_1) = f(\iota_2)$, it holds that $\mathbb{P}_{\iota_1}(E) = \mathbb{P}_{\iota_2}(E)$.
			\label{property:PAlgorithmConsistent}
		\end{enumerate}
		It is not immediately clear whether condition \ref{property:PAlgorithmRTMeasurable} for a given RGA $(X,\mathcal{F},\{\pr_\iota\}_{\iota\in\Omega})$ holds independently of the choice of the enumeration of $\Lambda$. This is indeed the case and was established in \cite[Lemma 9.2]{comp}.
	\end{defn}
	
	\MG{\begin{rem}
			Strictly speaking, a general algorithm is a pair $\Gamma$ and a suitable map $\Lambda_{\Gamma}$ taking elements of $\Omega$ to subsets of $\Lambda$. We will abuse notation somewhat and refer to a general algorithm as just the map $\Gamma$, but it is worth considering this subtlety when considering Definition \ref{definition:ProbablisticAlgorithm}. Similarly, we will refer to the mapping $\gprob_{\iota}$ as an RGA even though the full definition of an RGA requires $X$, $\mathcal{F}$ and $\{\mathbb{P}_{\iota}\}_{\iota \in \Omega}$ to be specified. 
	\end{rem}}
	
	We should justify this model before we move on: (i) and (ii) are measure theoretic \MG{conditions} that ensure that natural sets that one might construct for randomized algorithms (such as the minimum 
	amount of input information) are measurable sets\MG{.} \MG{T}hese assumptions are satisfied by any classical probabilistic 
	model such as randomized Turing machines or randomized BSS machines. The third point ensures the consistency of the model, in the sense that if the model 
	reads the same information for two inputs then the probability distribution of the outputs should be the same.
	
	\subsection{Inexact input and breakdown epsilons} \label{inexact}
	We now introduce the notion of inexact input for general computational problems. \MG{For the purpose of doing this, we will use the notation $\mathbb{D}$ to represent the dyadic numbers and similarly $\mathbb{D}_n:= \{ k 2^{-n} \: | \: k \in \Z\}$.}

	\begin{defn}($\Delta_1$-information, \cite[Definition 8.12]{comp}).\label{def:delta1_infromation}
		Suppose we are given a computational problem $\{\Xi, \Omega, \M,\Lambda\}$ and that $\Lambda = \{ f_i\}_{i \in I}$ where $I$ is some index set that can be 
		finite or infinite. As previously mentioned, in many cases we are forced to deal with inexact inputs. In this case we can not access $f_i(\iota)$, but rather an approximation $f_{i,n}(\iota)$ where 
		$f_{i,n}(\iota) \to f_i(\iota)$ as $n \to \infty$. Throughout this 
		paper we will assume that this can be done with error control. More precisely, we assume that for each 
		$n \in \N$ and each $i \in I$ there exists an
		$f_{i,n}: \Omega \to \mathbb{D}_n + i\mathbb{D}_n$, such that
		\begin{align} \label{inexact_input}
			\| \{ f_{i,n}(\iota) \}_{i \in I} - \{ f_i(\iota) \}_{i \in I} \|_{\infty} \leq 2^{-n} \quad \forall \; \iota \in \Omega.
		\end{align}
		If \MG{for each $n \in \N$ and $i \in I$ we have a function $f_{i,n}: \Omega \to \mathbb{D}_n + i\mathbb{D}_n$ such that \eqref{inexact_input} holds for $\{f_{i,n}\}_{i \in I}$, then}
		we say that the \MG{family} $\hat{\Lambda} = \{ f_{i,n} \: | \: i \in I, \: n \in \N \}$ provides \emph{$\Delta_1$-information} for $\{\Xi, \Omega, \M,\Lambda\}$.
		Moreover, we denote the family of all such $\hat{\Lambda}$ by $\mathcal{L}^1(\Lambda)$. 
	\end{defn}
	
	We want our algorithms to be able to deal with inexact input, in other words we want to have algorithms that can handle the computational problems 
	$\{\Xi, \Omega, \M,\hat{\Lambda}\}$ for all possible choices of $\hat{\Lambda} \in \mathcal{L}^1(\Lambda)$. In order to formalize this we introduce computational problems with $\Delta_1$-information:
	
	\begin{defn} (\cite[Definition 8.15]{comp}). \label{def:delta1}
		Given $\{\Xi, \Omega, \M,\Lambda\}$ with $\Lambda = \{ f_i \}_{i \in I}$ the corresponding computational problem with $\Delta_1$-information is defined as 
		\begin{align*}
			\{ \Xi, \Omega, \M, \Lambda\}^{\Delta_1} = \{\tilde{\Xi}, \tilde{\Omega}, \M, \tilde{\Lambda} \},
		\end{align*}
		where 
		\begin{equation}\label{eq:tildeOmega}
			\tilde{\Omega} = \{  \tilde{\iota} = \{ f_{i,n}(\iota)\}_{i \in I, n \in \N} \:  | \: \iota \in \Omega, \: 
			\{f_{i,n}\}_{i \in I} \; \text{satisfying \eqref{inexact_input} for all } n \in \N \},
		\end{equation}
		$\tilde{\Xi}(\tilde{\iota}) = \Xi(\iota)$
		and 
		$\tilde{\Lambda} = \{\tilde{f}_{i,n} \}_{i \in I, n \in \N}$ where 
		$\tilde{f}_{i,n}(\tilde{\iota}) = \tilde{\iota}_{i,n}$. Given an 
		$\tilde{\iota} \in \tilde{\Omega}$, there is a unique $\iota \in \Omega$ for which 
		$\tilde{\iota} = \{f_{i,n}(\iota) \}_{i \in I, n \in \N}$. We say that this 
		$\iota \in \Omega$ corresponds to $\tilde{\iota} \in \tilde{\Omega}$.
	\end{defn}
	
	We interpret the computational problem $\{ \Xi, \Omega, \M, \Lambda\}^{\Delta_1}$ as follows: 
	The domain $\tilde{\Omega}$ is the collection of all the 
	sequences approximating the inputs in $\Omega$, and we say that an algorithm $\Gamma$ works on inexact input if the algorithm returns an approximation to the target value, with an arbitrary desired accuracy, for all $\tilde{\iota} \in \tilde{\Omega}$, that is, for any sequence approximating $\iota$. 
	For later use, we need to specify how the input $\tilde{\iota} \in \tilde{\Omega}$ is passed to a Turing machine $T$ as an input. For this we need to assume that the index set $I$ for $\Lambda$ is countable\MG{. In }the context of this paper this is \MG{immediate}, since in our case the index set $I$ is  always finite. With this observation in place $\tilde{\iota}$ is 
	represented by an oracle tape that $T$ can access, 
	where on input $(i,n) \in I \times \N$ the oracle returns the unique finite binary string representing 
	$\tilde{f}_{i,n}(\tilde{\iota})$.

	\begin{rem}(\MG{Computable input}, \cite[Remark 8.13]{comp}). \label{rem:TurVsMarkov}
		It is possible to consider a restriction of $\Delta_1$ information wherein \MG{for each $\iota \in \Omega$ and} $j \in I$, $\{f_{j,n}\MG{(\iota)}\}_{n\in\mathbb{N}}$ forms a computable sequence. This restriction strengthens the negative results and weakens the positive results. See Remark \ref{remark:CompDelta1Markov} for further details.
	\end{rem}
	
	We wish to investigate the constructibility of optimal neural networks in the theory of underdetermined systems. That is, investigate the existence of training algorithms that \MG{produce} a neural network that 
	solves the computational problem stated in \eqref{problem}. 
	To do this we need the notion of breakdown epsilons.

	\begin{defn}\label{def:StrongBDEPS}(Strong and probabilistic strong breakdown epsilons 
		\cite[Definition 8.18 and Definition 8.28]{comp}).
		Given a computational problem
		$\{\Xi,\Omega,\mathcal{M},\Lambda\}$, the strong breakdown-epsilon 
		$\epsilon_{\mathrm{B}}^{\mathrm{s}} \: \in [0,\infty]$ is given by
		\begin{align*}
			\epsilon_{\mathrm{B}}^{\mathrm{s}} = \sup \{ \epsilon > 0 \: | \: \forall \: \MG{\text{general }}\text{algorithms} \: 
			\Gamma, \exists \: \iota \in 
			\Omega \: \text{such that} \: \inf_{\xi \in \Xi(\iota)} d_{\M}(\Gamma(\iota),\xi)>\epsilon \}\MG{.}
		\end{align*}
		The strong probabilistic breakdown-epsilon $\epsilon_{\mathbb{P}\mathrm{B}}^{\mathrm{s}}:[0,1) \to 
		\: [0,\infty]$ is given by
		\begin{align*}
			\epsilon_{\mathbb{P}\mathrm{B}}^{\mathrm{s}}(p) = 
			\sup \{ \epsilon &> 0 \: | \: \forall \: \MG{(X,\mathcal{F},(\mathbb{P}_{\iota})_{\iota })} \in \text{RGA} \: \exists \: 
			\iota \in \Omega \: \text{such that} \: \mathbb{P}_{\iota}(\inf_{\xi \in \Xi(\iota)} 
			d_{\M}(\gprob(\iota),\xi) > 
			\epsilon) > p \}.
		\end{align*}
		In both definitions we use the convention $\sup \emptyset = 0$.
	\end{defn}
	
	As established in \cite{comp}, impossibility results for randomised algorithms can differ if one considers only those algorithms that halt on every input, leading to the following two definitions.
	
	\begin{defn}(Halting randomised general algorithms, \cite[Definition 8.29]{comp}). 
		\label{def:halting_randomised}
		A  randomised general algorithm $\gprob$ for a computational problem 
		$\{\Xi,\Omega,\mathcal{M},\Lambda\}$ is called a \emph{halting randomised general algorithm} (hRGA) if $\mathbb{P}_{\iota}(\gprob_{\iota} = \nh) = 0$, for all $\iota\in\Omega$. We denote the class of all halting randomised general algorithms by $\mathrm{hRGA}$.
	\end{defn}

	\begin{defn}(Probabilistic strong halting breakdown epsilon, \cite[Definition 8.30]{comp}). 
		\label{prob_strong_break_halt}
		Given the computational problem 
		$\{\Xi,\Omega,\mathcal{M},\Lambda\}$, where $\Lambda =\{f_k \, \vert \, k \in \mathbb{N}, \, k\leq |\Lambda| \}$, we define the \emph{halting probabilistic strong Breakdown-epsilon} $\epsilon_{\mathbb{P}h\mathrm{B}}^{\mathrm{s}}: [0,1) \to \: [0,\infty]$ according to
		\begin{align*}
			\epsilon_{\mathbb{P}h\mathrm{B}}^{\mathrm{s}}(\mathrm{p}) = \sup\{&\epsilon \geq 0, \, \vert \, \forall \, \gprob \in \mathrm{hRGA} \,\,\exists \, \iota \in \Omega \text{ such that } 
			\mathbb{P}_{\iota}(\inf_{\xi \in \Xi(\iota)} d_{\M}(\gprob(\iota),\xi) > \epsilon) > \mathrm{p}\},
		\end{align*}
		where $\gprob_{\iota}$ is defined in \ref{property:PAlgorithmMeasurable} in Definition \ref{definition:ProbablisticAlgorithm} and where, as in Definition \ref{def:StrongBDEPS}, we treat $\sup \emptyset = 0$. 
	\end{defn}

	Throughout this paper we wish to prove the strongest possible lower bounds of computability. Thus we allow our general algorithms to perform arbitrary general 
	operations and we prove all our lower bounds for such algorithms. On the other hand, we also wish to have the strongest possible positive results, thus when we construct algorithms, these algorithms will always be recursive (That is, possible to implement on a Turing machine).
	
	\subsection{Formal description of the problem} \label{comp_problem}
	
	In order to prove the \MG{non-computability results} asserted in \S \ref{main_results_section} we formulate our problem as a computational problem (in the sense of Definition \ref{def:GenCompProblem}).
	We start by formally defining what we mean by a neural network. Deep learning is a rapidly developing field, where new constructions and 
	architectures of neural networks are constantly proposed. Our aim is to capture as many of these constructions as possible with our results, 
	\MG{and} thus we propose the following general formal definition of a neural network:
	\begin{defn}\label{def:NNs_with_rhos}
		A neural network $\mathbf{N}: \MG{\R}^m \to \MG{\R}^N$ is a map of the form
		\begin{align*}
			\mathbf{N}(y^0) = V^L(y^0, \dots, y^{L-1}),
		\end{align*}
		where $y^1,y^2,\dotsc,y^{L-1}$ satisfy the following. For $j=1,2,\dotsc,L$, there is an affine map $V^j: \MG{\R}^{N_{0}} \times \dots \times \MG{\R}^{N_{j-1}} \to \MG{\R}^{N_j}$ of the form 
		$V^j(y^0, \dots, y^{j-1}) =   W^{j,1}y^0 + \dotsb + W^{j,j-1}y^{j-2} + W^{j,j}y^{j-1} + b^j$ for $j = 1, \dots, L$, where $N_0 = m$, $N_L = N$ and $N_1,N_2,\dotsc,N_{L-1}$ are integers, where $W^{j,i} \in \MG{\R}^{N_j \times N_{i-1}}$ for $i \in \{1, \dots, j\}$ and where $b^j \in \MG{\R}^{N_j}$. Then
		$y^i = \rho_i(V^i(y^0, \dots, y^{i-1}))$ for $i \in 1, \dots, L-1$ with $\rho_i$ an activation function.
		The activation functions $\rho_i: \R \to \R$ are potentially non-linear functions that act component-wise. 
	\end{defn}
	\MG{\begin{rem}[More general neural networks]
			The proof techniques will typically apply to more general classes of neural networks (e.g. those that make use of the softmax activation function) -- for the sake of keeping the presentation concise, we do not comment further on this issue. 
		\end{rem}}
	\MG{\begin{rem}[Class of permissible activation functions]
			We assume that the activation functions $\rho_i$ for $i\in \{1,2,\dotsc,L-1\}$ (which can vary depending on the layer of the neural network) are all contained in some set $\mathcal{A}$ of permitted activation functions. The only requirements we make on this set $\mathcal{A}$ are that every function in $\mathcal{A}$ is computable (to ensure that the concept of computing an NN with an algorithm is well-defined -- see Remark \ref{rem:computeNN}), that the map $x \mapsto x^2$ for $x \in \R$ is contained in $\mathcal{A}$, and that there is an $f \in \mathcal{A}$ such that the restriction of $f$ to non-negative $t$ satisfies $f(t) = 1/(1+t)$: these last two conditions ensure that the radial basis function $\phi(x) = 1/(1+x^2)$ can be expressed as a neural network. The same arguments made in the paper could also be made with a different choice of radial basis function -- if this were done, the Lipschitz constants discussed in Theorem \ref{thm_intro:second} would change, and we would have slightly different requirements on the functions contained in $\mathcal{A}$ -- but this is beyond the scope of this paper. Since our results (both positive and negative) hold for any such choice of $\mathcal{A}$, we shall not explicitly mention this class again.
	\end{rem}}
	\begin{rem}
		The affine dependence of $V^j$ on $(y^0, \dots, y^{j-1})$ allows for skip connections from the input of previous layers to any layer of the neural network. This allows us to obtain standard architectures such as residual networks \cite{res1,res2,Genzel} with an arbitrary number of 
		layers between the skip connections.
	\end{rem}
	
	The domain $\Omega$ is the set of objects that gives rise to our 
	computational problem. In our setting these objects are ordered training sets. 
	Since the error of an algorithm is mostly interesting relative to the size and 
	the bounds of the training set, we will consider bounded training sets of fixed size throughout 
	this paper. In addition, the elements of a training set are processed by a training algorithm in some order\MG{. To} encode this we therefore assume that every training set is represented as a list of elements. More formally, let $A: \R^N \to \R^m$ be a linear map, $\ell \in \mathbb{N}$, and let
	\begin{align}\label{omega}
		\Omega \subset  \{ \T = ((x^k,y^k))_{k = 1}^{\ell} \, \vert \, \T \in T_{\ell}(A)\}, 
	\end{align}
	be the domain, where we recall $T_{\ell}(A)$ from \eqref{eq:T_ell}.

	The set of measurements $\Lambda$ is the collection of functions that provide us with the information we are allowed to read as an input to an algorithm. We define the measurements as follows: Given a training set $\T = \MG{((x^k,y^k))_{k = 1}^{\ell} }\in \Omega$, let $f_{x,i}^k$ be given by 
	$f_{x,i}^k(\T) = x^k_i$, where the index $i$ denotes the $i$'th coordinate of the vector $x^k$. 
	We define $f_{y,j}^k$ in the same way to measure the $y$-coordinates, more precisely 
	$f^k_{y,j}(\T) = y^k_j$. In summary, we define 
	$\Lambda$ to be the collection
	\begin{equation}\label{measurements}
		\Lambda = \{f^k_{y,j},f^k_{x,i}: \Omega \to \R \: | \: i = 1, \dots, N \;, \; j = 1, \dots, m, \; 
		\text{and} \; k = 1, \dots, \ell\}.
	\end{equation}
	Next, we define a precise notation for an inexact representation of the elements in the domain 
	$\Omega$. 
	Let $\hat \Lambda = \{ f_{n} \: | \: f  \in \Lambda, \: n \in \N \}$ be a set that provides 
	$\Delta_1$-information for $\Omega$ as defined in Definition \ref{def:delta1_infromation}.
	Then, for an arbitrary $\T \in \Omega$ and $k \in \{1, \dots, \ell\}$, we define the corresponding inexact representation of $(x^k,y^k)$ to be the pair of sequences 
	$(\tilde x^k,\tilde y^k)$, where
	\begin{align} \label{def:inexact_rep}
		\tilde x^k =  \{\tilde x^{k}(n)\}_{n \in \N} = \{\{f_{x,i,n}^k(\T)\}_{i =1}^{N}\}_{n \in \N} \quad \text{and} \quad 
		\tilde y^k =  \{\tilde y^{k}(n)\}_{n \in \N} =  \{ \{ f_{y,j,n}^k(\T)\}_{j = 1}^{m}\}_{n \in \N},
	\end{align}
	\MG{so that equation \eqref{inexact_input} can be interpreted as the requirement
		\begin{equation}\label{eq:linfestimateDelta1}
			\|\tilde x^{k}(n) - x^k\|_{\infty} \leq 2^{-n}, \qquad \|\tilde y^{k}(n) - y^k\|_{\infty} \leq 2^{-n}.
	\end{equation}}
	Note that it is clear that one can recursively transform the estimate \eqref{eq:linfestimateDelta1} into an $\ell_2$ estimate. Thus, we may without loss of generality assume that (see \eqref{eq:approx})
	\begin{equation}\label{eq:l2estimateDelta1}
		\|\tilde x^{k}(n) - x^k\|_{\MG{2}} \leq 2^{-n}, \qquad \|\tilde y^{k}(n) - y^k\|_{\MG{2}} \leq 2^{-n}.
	\end{equation}
\MG{	where a consequence of the recursive transformation is that $\tilde x^k(n)$ and $\tilde y^k(n)$ are now assumed to have each component in $\mathbb{D}$, as opposed to $\mathbb{D}_n$.}
	\MG{In particular, as in \S\ref{sec:model_comp}, an inexact representation $\tilde \T$ which corresponds to $\T \in \Omega$ can be listed as \[\tilde \T = (\{(\tilde x^1(n),\tilde y^1(n))\}_{n=1}^{\infty},\{(\tilde x^2(n),\tilde y^2(n))\}_{n=1}^{\infty},\{(\tilde x^3(n),\tilde y^3(n))\}_{n=1}^{\infty},\dotsc, \{(\tilde x^\ell(n),\tilde y^\ell(n))\}_{n=1}^{\infty})\]
	with $\tilde x^l(n)$ and $\tilde y^l(n)$ as in \eqref{def:inexact_rep} and satisfying \eqref{eq:linfestimateDelta1}, or equivalently, satisfying \eqref{eq:l2estimateDelta1} -- we work with whichever concept is most convenient for a particular argument.
	}
	
	At last, given a matrix $A \in \R^{m \times N}$ and a collection 
	$\{ (\T, \M_1(\T)) \: | \: \T \in \Omega \}$ \MG{ (where we recall $\M_1$ from \eqref{M1_dep})}
	we define the problem function $\Xi$ in the following way
	\begin{equation*}
		\Xi: \Omega \rightrightarrows  \NN_{m,N}^{\F}, \text{  such that  } \,\,  
		\Xi(\mathcal{T})  = \left\{ \mathbf{N}_{opt}^{\M_1(\T)} \: : \: \mathbf{N}_{opt}^{\M_1(\T)} 
		\: \text{ is optimal for } (A,\mathcal{M}_1(\T)) \right\},
	\end{equation*}
	where $\NN_{m,N}^{\F}$ is the set of 
	neural networks of real input dimension $m$ and real output dimension $N$ that are bounded on 
	$\bigcup_{\mathcal{F}} \M_2 
	= \bigcup_{(A, \M_1(\T)) \in \mathcal{F}} A(\M_1(\T))$. We consider $\NN_{m,N}^{\F}$ as a metric space $\mathcal{M}$ equipped with the following metric, for 
	$\mathbf{N_1}$,$\mathbf{N_2} \in \NN_{m,N}^{\F}$,
	\begin{equation} \label{metric}
		d(\mathbf{N_1},\mathbf{N_2}) = \sup_{z \in \bigcup_{\mathcal{F}} \M_2} \| \mathbf{N_1}(z) - \mathbf{N_2}(z)\|_{\MG{2}}, 
	\end{equation}
	where we identify -- \MG{for the purpose of defining the metric} -- the networks that are equal on $\bigcup_{\mathcal{F}} \M_2$.

	\begin{rem}[Computing a NN with an algorithm]\label{rem:computeNN}
		We need to specify how to interpret an algorithm $\Gamma$ producing a NN. In particular, \MG{we will discuss here} how the mapping
		\[
		\tilde \Omega \ni \mathcal{\tilde T} \mapsto \Gamma(\mathcal{\tilde T}) \in \NN_{m,N}^{\F}
		\]
		should be interpreted. In the Turing case, $\Gamma(\mathcal{\tilde T})$ is a string $\mathcal{S}_{\tilde \T}$ that identifies $\mathbf{N}_{\tilde \T} \in \NN_{m,N}^{\F}$ as follows. The string $\mathcal{S}_{\tilde \T}$ provides the affine mappings $V^j$ in Definition \ref{def:NNs_with_rhos} through the rational matrix and vector entries (or codes for these if they are computable non-rationals). In addition, the string provides codes (Turing machines) for the computable functions $\rho_j$ in Definition \ref{def:NNs_with_rhos}. This uniquely defines a neural network $\mathbf{N}_{\tilde \T}$, which is clearly Turing computable as a mapping from $\mathbb{R}^m$ to $\mathbb{R}^N$, since the affine mappings are clearly computable, and so is the composition of finitely many computable functions. Indeed, it is a standard argument (\cite{pour1989computability}, p. 29) that one can recursively combine Turing machines representing functions into one Turing machine computing the composite function. Hence, there is a Turing machine $\Gamma^1$ that produces a Turing machine $\Gamma^2$ that computes $\mathbf{N}_{\tilde \T}$ as follows: $\Gamma^2_{\tilde \T} = \Gamma\MG{^1}(\mathcal{S}_{\tilde \T})$ such that 
		$
		\|\Gamma^2_{\tilde \T}(\tilde y, 2^{-n}) - \mathbf{N}_{\tilde \T}(y)\|_{\MG{2}} \leq 2^{-n}, 
		$
		for any $\tilde y$ that is a $\Delta_1$ approximation to $y \in \mathbb{R}^m$. 
		Thus, we can throughout the paper use the slight abuse of notation 
		\[
		\Gamma(\mathcal{\tilde T})(y) \coloneqq \mathbf{N}_{\tilde \T}(y), \quad y \in \mathbb{R}^m  
		\]
		\MG{(or $\Gamma(\mathcal{\tilde T},\epsilon)(y)$ in the case where an additional parameter to the algorithm, e.g. $\epsilon$, is needed)}.
		In particular, $\Gamma(\mathcal{\tilde T})$ is used both for the actual NN and for the string $\mathcal{S}_{\tilde \T}$, however the meaning will be clear from the context.
	\end{rem}

	\begin{rem}[Algorithms that compute]\label{rem:computeNN2}
		All \MG{of} our lower bounds (negative results) hold for general algorithms \MG{which yields} universal lower bounds regardless of the computational model. However, all the upper bounds (positive results) hold for Turing machines \MG{which yields} the strongest possible results. A Turing machine $\Gamma$ that computes neural networks will take two inputs: $\tilde \T \in \tilde \Omega$ and a rational $\epsilon>0$ such that $\Gamma(\tilde \T, \epsilon)$ is within $\epsilon$ accuracy of the desired NN. However, in many cases we will -- to simplify notation -- suppress the second variable. Moreover, the statement: `An algorithm $\Gamma$, taking inputs in $\Omega$, that computes an optimal neural network' means that $\Gamma(\tilde \T, \epsilon)$ is within $\epsilon$ accuracy of \MG{ an optimal NN for all rational $\epsilon$, every $\T \in \Omega$ and every $\tilde \T \in \tilde \Omega$ corresponding to $\T$.}
	\end{rem}
	
	\subsection{An important preliminary result}
	
	We use the following important proposition (simplified from Proposition 9.5 in \cite{comp} to prove the non-existence of algorithms that construct optimal neural networks.
	
	\begin{prop}[Proposition 9.5 in \cite{comp}] \label{computability_prop}
		Let $\{\Xi, \Omega, \mathcal{M}, \Lambda\}$ be a computational problem with $\Lambda=\{f_k\,\vert\,k\in\mathbb{N}, k\leq|\Lambda| \}$ countable, let $\{\iota^1_n\}_{n=1}^{\infty}$ be a sequence in $\Omega$, and let $\iota^0 \in \Omega$. Assume that the following properties hold:
		\begin{enumerate}[leftmargin=8mm, label=(\alph*)]
			\item There are sets $S^1, S^2 \subset \mathcal{M}$ and $\kappa > 0$ such that
			$
			\inf_{x_1 \in S^1, x_2 \in S^2}d_{\mathcal{M}}(x_1,x_2) \geq \kappa
			$ and \MG{for all $n \in \mathbb{N}$}, $\Xi(\iota^1_n) \subset S^1$, as well as $\Xi(\iota^0) \subset S^2$. \label{property:MinimisersOfNonZero}
			\item For every $k\leq |\Lambda|$ we have that $|f_k(\iota^1_n) - f_k(\iota^0)|\leq 1/4^n$, for all $n \in\mathbb{N}$. \label{property:Del1Info}
		\end{enumerate}
		Then, there exists a 
		$\hat\Lambda\in \mathcal{L}^1(\Lambda)$
		such that, for the computational problem $\{\Xi,\Omega,\mathcal{M},\hat\Lambda\}$,  we have 
		$\epsilon^{\mathrm{s}}_{\mathrm{B}} \geq \strbdepsph \geq \strbdepsp \geq \kappa/2 $ for $\mathrm{p} \in [0,1/2)$.
	\end{prop}

	\begin{rem}
		Proposition \ref{computability_prop} is a result about how inexact input can \MG{preclude the existence of algorithms with good accuracy}, given that the domain and the problem function $\Xi$ have some unfortunate properties. 
		In a way the above statement is very intuitive and should be read as follows: If two \MG{inputs} that are arbitrarily close in the domain get mapped far apart by the problem function $\Xi$, then any algorithm that works with inexact input (that is, works with the computational problem 
		$\{ \Xi, \Omega, \M, \Lambda\}^{\Delta_1}$) will break down. In fact, the result even states something slightly stronger: There exists one specific 
		$\hat{\Lambda} \MG{= \{f_{i,n}(\iota)\}_{i \in I, n \in \N}}$ that gives $\Delta_1$-information for $\{ \Xi, \Omega, \M, \Lambda\}$ such that all algorithms $\Gamma$ break down for the corresponding computational problem $\{ \hat{\Xi}, \hat{\Omega}, \M, \hat{\Lambda}\}$, where 
		\begin{align*}
			\hat{\Omega} = \{ \hat{\iota} = \{ f_{i,n}(\iota)\}_{i \in I, n \in \N} \:  | \: \iota \in \Omega \; \},
		\end{align*}
		and $\hat{\Xi}(\hat{\iota}) = \Xi(\iota)$. As before, we notice that there is a canonical bijection between $\Omega$ and $\hat{\Omega}$ and by 
		identifying $\Omega \simeq \hat{\Omega}$ the mapping $\hat{\Xi}$ becomes the same as $\Xi$. Because of this identification we will often write 
		$\{ \Xi, \Omega, \M, \hat{\Lambda}\}$ in place of $\{ \hat{\Xi}, \hat{\Omega}, \M, \hat{\Lambda}\}$.
	\end{rem}

	\begin{rem} [Computability of $\Delta_1$ information]\label{remark:CompDelta1Markov}
		Suppose that
		$\Omega = \{\iota_n^1\}_{n \in \N} \cup \{\iota^0\}$ satisfies condition (b) in 
		Proposition \ref{computability_prop}. 
		Then, by assuming that $f_k(\iota)$ are computable numbers, in the sense of Turing-computability, for each $k \MG{\leq }|\Lambda|$ \MG{(when $\Lambda$ is a finite set, otherwise for each $k < |\Lambda|$)} and for each $\iota \in \Omega$, we have the following:  
		$\hat \Lambda$ can be chosen so that $\hat \Lambda = \{f_{k,n} \, \vert \, k \leq |\Lambda| \text{ and } 
		n \in \mathbb{N}\}$ \MG{(respectively, when $|\Lambda|$ is infinite, $\hat \Lambda = \{f_{k,n} \, \vert \, k \in \mathbb{N} \text{ and } 
			n \in \mathbb{N}\}$) } and so that, for any $\iota \in \Omega$ and each $f_k \in \Lambda$ there exists a Turing machine, that given input $n \in \mathbb{N}$, outputs $ f_{k,n}(\iota)$, see 
		\cite[Remark 9.6]{comp} for details.
		
	\end{rem}
	
	\begin{rem}[`Oracle' as input and Markov algorithms]\label{rem:Markov}
		The Markov model (see Remark \ref{oracle2}) -- where the algorithm for the computable input is the actual input to the Turing machine -- dates back to Turing's legendary 1936 paper \cite{Turing_Machine}. However, this model is usually referred to as the Markov model \cite{MarkovModel}. \MG{The Markov model, and the model where the Turing machine takes a computable infinite string representing a computable number, are equivalent when considering computable functions on any $[a,b] \cap \mathbb{R}_c$ (an interval $[a,b]$ intersecting the computable numbers, and similarly for vector valued computable functions). The equivalence between the two models was established in \cite{Shoenfield, Markov}, however, see \cite{Hertling} (Theorem 6) for a modern exposition. This means that, in many cases, to establish impossibility results in the Markov model, it suffices to demonstrate impossibility results given computable $\Delta_1$-information as input. In view of Remark \ref{remark:CompDelta1Markov}, many of our results (when restricted to Turing machines/Markov algorithms and not general algorithms) also follow in the Markov model. Moreover, the proof of Theorem \ref{thm:first:technical} (iv) is an example of a direct proof in the Markov model. }  
\end{rem}

	\subsection{Recalling the basics from the SCI hierarchy}
	
	Here we provide an informal review of the basics of the SCI hierarchy for an easy reference. 
	The mainstay of the hierarchy are the $\Delta^{\alpha}_k$ classes. The $\alpha$ is related to the model of computation, for example, when $\alpha = A$ we are in the classical Turing setup. In particular, given a collection $\mathcal{C}$ of computational problems, then
	\begin{itemize}
		\item[(i)] $\Delta^{\alpha}_0$ is the set of problems \MG{in $\mathcal{C}$} that can be computed \MG{exactly} in finite time, the SCI $=0$.
		\item[(ii)] $\Delta^{\alpha}_1$ is the set of problems \MG{in $\mathcal{C}$} that can be computed using one limit 
		(the SCI $\leq\: 1$) with control of the error, i.e. there exists a sequence of algorithms $\{\Gamma_n\}$ such that 
		$\text{dist}_{\M}(\Gamma_n(\iota), \Xi(\iota)) := \inf_{\xi \in \Xi(\iota)} d_{\mathcal{M}}(\Gamma_n(\iota),\xi) \leq 2^{-n}, \, \forall \iota \in \Omega$.
		\item[(iii)] $\Delta^{\alpha}_2$ is the set of problems \MG{in $\mathcal{C}$} that can be computed using one limit 
		(the SCI $\leq\: 1$) without error control, i.e. there exists a sequence of algorithms $\{\Gamma_n\}$ such that $\lim_{n\rightarrow \infty}\Gamma_n(\iota) = \Xi(\iota), \, \forall \iota \in \Omega$.
		\item[(iv)] $\Delta^{\alpha}_{m+1}$, for $m \in \mathbb{N}$, is the set of problems \MG{in $\mathcal{C}$} that can be computed by using $m$ limits, (the SCI $\leq m$), i.e. there exists a family of algorithms $\{\Gamma_{n_m, \hdots, n_1}\}$ such that 
		\begin{equation}\label{eq:SCI_limits}
			\lim_{n_m \rightarrow\infty}\hdots \lim_{n_1\rightarrow\infty}
			\text{dist}_{\M}(\Gamma_{n_m,\hdots, n_1}(\iota),\Xi(\iota)) = 0, \, \forall \iota \in \Omega.
		\end{equation}
	\end{itemize}
	In general, this hierarchy cannot be refined unless there is some extra structure on the metric space $\mathcal{M}.$ The hierarchy typically does not collapse, and we have:
	\begin{equation}\label{SCI1}
		\Delta_0^{\alpha} \subsetneq \Delta_1^{\alpha} \subsetneq \Delta_2^{\alpha} \subsetneq \hdots \subsetneq \Delta^{\alpha}_{m} \subsetneq \hdots.
	\end{equation}
	However, depending on the collection $\mathcal{C}$ of computational problems, the hierarchy \eqref{SCI1} may terminate for a finite $m$, or it may continue for arbitrarily large $m$. We will focus on the lower parts of the hierarchy in \eqref{SCI1} in this paper, however, for the interested reader we point out that for certain metric spaces $\mathcal{M}$ one can extend \eqref{SCI1} to the full hierarchy and define the $\Pi_j^{\alpha}$ and $\Sigma_j^{\alpha}$ classes for $j \in \mathbb{N}$. We then get the following hierarchy:
	\begin{equation}\label{SCI_hierarchy}
		\begin{tikzpicture}[baseline=(current  bounding  box.center)]
			\matrix (m) [matrix of math nodes,row sep=1.2em,column sep=1.5em] {
				\Pi_0^{\alpha}   &                    & \Pi_1^{\alpha} &    &  \Pi_2^{\alpha}&  & {}\\
				\Delta_0^{\alpha}&  \Delta_1^{\alpha} & \Sigma_1^{\alpha}\cup\Pi_1^{\alpha} & \Delta_2^{\alpha}&      \Sigma_2^{\alpha}\cup\Pi_2^{\alpha} & \Delta_3^{\alpha}& \cdots\\
				\Sigma_0^{\alpha}&                    & \Sigma_1^{\alpha} & &  \Sigma_2^{\alpha}&  &{} \\
			};
			\path[-stealth, auto] (m-1-1) edge[draw=none]
			node [sloped, auto=false,
			allow upside down] {$=$} (m-2-1)
			(m-3-1) edge[draw=none]
			node [sloped, auto=false,
			allow upside down] {$=$} (m-2-1)
			
			(m-2-2) edge[draw=none]
			node [sloped, auto=false,
			allow upside down] {$\subsetneq$} (m-2-3)
			(m-2-3) edge[draw=none]
			node [sloped, auto=false,
			allow upside down] {$\subsetneq$} (m-2-4)
			(m-2-4) edge[draw=none]
			node [sloped, auto=false,
			allow upside down] {$\subsetneq$} (m-2-5)
			(m-2-5) edge[draw=none]
			node [sloped, auto=false,
			allow upside down] {$\subsetneq$} (m-2-6)
			(m-2-6) edge[draw=none]
			node [sloped, auto=false,
			allow upside down] {$\subsetneq$} (m-2-7)

			(m-2-1) edge[draw=none]
			node [sloped, auto=false,
			allow upside down] {$\subsetneq$} (m-2-2)
			(m-2-2) edge[draw=none]
			node [sloped, auto=false,
			allow upside down] {$\subsetneq$} (m-1-3)
			(m-2-2) edge[draw=none]
			node [sloped, auto=false,
			allow upside down] {$\subsetneq$} (m-3-3)
			(m-1-3) edge[draw=none]
			node [sloped, auto=false,
			allow upside down] {$\subsetneq$} (m-2-4)
			(m-3-3) edge[draw=none]
			node [sloped, auto=false,
			allow upside down] {$\subsetneq$} (m-2-4)
			(m-2-4) edge[draw=none]
			node [sloped, auto=false,
			allow upside down] {$\subsetneq$} (m-1-5)
			(m-2-4) edge[draw=none]
			node [sloped, auto=false,
			allow upside down] {$\subsetneq$} (m-3-5)
			(m-1-5) edge[draw=none]
			node [sloped, auto=false,
			allow upside down] {$\subsetneq$} (m-2-6)
			(m-3-5) edge[draw=none]
			node [sloped, auto=false,
			allow upside down] {$\subsetneq$} (m-2-6)
			(m-2-6) edge[draw=none]
			node [sloped, auto=false,
			allow upside down] {$\subsetneq$} (m-1-7)
			(m-2-6) edge[draw=none]
			node [sloped, auto=false,
			allow upside down] {$\subsetneq$} (m-3-7);
			
		\end{tikzpicture}
	\end{equation}
	For details about \eqref{SCI_hierarchy} see \cite{Hansen_JAMS, SCI, CRAS}. 
	
	%\begin{rem}[Multivalued functions] When dealing with multivalued problems one needs a framework that can handle multiple solutions. As the setup above does not allow $\Xi$ to be multi-valued we need some slight changes.  We allow $\Xi$ to be multivalued, even though towers of algorithms are not. Hence, the only difference to the standard SCI hierarchy is that the last limit in \eqref{eq:SCI_limits}
	%	is replaced by
	%	\[
	%	\disM(\Xi(\iota),\Gamma_{n_m}(\iota)) \longrightarrow 0, \qquad n_m \rightarrow \infty, 
	%	\]
	%	where 
	%	$
	%	\disM(\Xi(\iota),\Gamma_{n_m}(\iota)) := \inf_{x \in \Xi(\iota)}d_{\mathcal{M}}(x,\Gamma_{n_m}(\iota)).
	%	$ 
	%\end{rem}

	\section{Formal statements and proofs of the main results} \label{section:proofs}
	
	Before we embark on the proofs of the theorems we will introduce some basic notation and discuss some vocabulary that will be used in the proofs.
	First, we \MG{recall} some standard projections. Indeed, \MG{for integers $N,m$ and $\ell$, as well as $\T \in(\R^N \times \R^m)^\ell$, we define} $\pi_1(\T)$ and $\pi_2(\T)$ \MG{ as in  \eqref{eq:pi1def} and \eqref{eq:pi2def}.}
	
	In addition, we introduce a notation to describe the set theoretic properties of training sets.
	Since we require a training set $\T$ to have bounded elements, we \MG{recall} the set $T_{\ell}$ which captures this property: Let $A: \R^N \to \R^m$ be a linear map and 
	\MG{\begin{align*}T_{\ell}(A) := \{ \T \in(\R^N \times \R^m)^\ell \, \vert \, & \T = ((x^k,y^k))_{k = 1}^{\ell} \text{ with }
	\|x^k\|_{\MG{2}}, \|y^k\|_{\MG{2}} \leq 1 \: \text{and } \: y^k = Ax^k,\\& (x^k,y^k) \neq (x^j,y^k) \text{ for every } j,k \in \{1,2,\dotsc,\ell\} \text{ with } j \neq k\}.\notag
	\end{align*}}
	Secondly, recall Definition \ref{def:delta1_infromation} and Definition \ref{def:delta1} where we discuss the concepts of $\Delta_1$-information and computational 
	problems with $\Delta_1$-information. In particular, we have the original domain $\Omega$ and the corresponding domain $\tilde \Omega$ that contains all 
	sequences of approximations to the elements in $\Omega$. We will typically use the notation $\mathcal{T} \in \Omega$ and $\mathcal{\tilde T} \in \tilde \Omega$ \MG{where $\mathcal{\tilde T}$ is assumed to correspond to $\mathcal{T}$}. \MG{We also recall from \eqref{eq:extended-metric} that if an algorithm doesn't halt on an input then its error is infinity.}
	At last, we use $\| \cdot \|_{op}$ to denote the standard operator norm, more precisely for a linear map
	$M: \R^N \to \R^m$ we have that
	\begin{align} \label{op_norm}
		\|M\|_{op} = \sup \{ \|Mx\|_{\MG{2}} \: : \: \|x\|_{\MG{2}} \leq 1\}.
	\end{align}

	\subsection{Useful propositions and lemmas}
	We start by presenting a driving proposition, from which parts of our main results, 
	more specifically, part (i) of Theorem \ref{thm_intro:second} and the negative parts of (i)-(iii) of Theorem \ref{thm_intro:first}, will follow as corollaries. 
	
	\begin{prop} \label{prop:driving}
		Let $A \in \R^{m \times N}$ be a non-zero matrix with non-trivial kernel \MG{with $N \geq 3$}. Then, for any 
		$\kappa \in (0,3/8]$ and $\ell \geq 2$, 
		there exist infinitely many domains $\MG{\Omega \subset T_{\ell}(A)}$ 
		of training sets, which give rise to infinitely many computational problems 
		$\{ \Xi, \Omega, \M', \Lambda\}$, where $\Omega$, $\Lambda$ and $\Xi$ are as described in 
		\S \ref{comp_problem}, and where the metric space $\M = \NN_{m,N}^{\F}$ is replaced by $\M' = 
		\NN_{m,N}^{\F}$ 
		with the metric
		\begin{align} \label{breakdown_metric}
			d_{\M'}(\mathbf{N_1},\mathbf{N_2}) = \sup_{z \in \bigcap_{\mathcal{F}} \M_2} 
			\| \mathbf{N_1}(z) - \mathbf{N_2}(z)\|_{\MG{2}},
		\end{align}
		where we identify \MG{any two} neural networks that \MG{agree} on $\bigcap_{\mathcal{F}} \M_2$, such that we have the following.
		For each domain
		$\Omega$ the mapping in \eqref{eq:Xi} satisfies $\Xi(\iota) 
		\neq \emptyset$ for all $\iota \in \Omega$. Furthermore, there exists a sequence 
		$\{ \iota_n^1\}_{n \in \N},$ and an $\iota^0$ in $\Omega$ so that $\Omega = \{ \iota_n^1\}_{n \in \N} \cup \{\iota^0\}$
		and such that the conditions (a)-(b) in Proposition \ref{computability_prop} are satisfied for 
		$\{\iota_n^1\}_{n \in \N}$ and  $\iota^0$. Finally, if $A$ is computable then $\Omega$ can be chosen so that for any $f \in \Lambda$ there are algorithms \MG{(which can be implemented in the Turing Machine model)} $\Gamma^0: \N \to \mathbb{Q}$ and $\Gamma^1: \N \times \N \to \mathbb{Q}$ with 
		\begin{equation}\label{eq:ComputabilityOfTheIotas}
			|\Gamma^0(k) - f(\iota^0)|\leq 2^{-k}, \quad |\Gamma^1(k,n) - f(\iota^1_n)| \leq 2^{-k}.
		\end{equation}
		
	\end{prop}
	
	\begin{proof}
		The proof is divided into several steps.
		
		\textbf{Step I} (\emph{Defining $\Omega$, $\iota^0$ and $\{\iota^1_n\}_{n=1}^{\infty}$}).
		Since $A$ has a 
		non-trivial kernel, we can pick a $\x \in \ker(A)$ such that $\| \x \|_{\MG{2}} = 2\kappa$. 
		Since $A$ is non-zero, we can pick a unit vector 
		$w \in \ker(A)^{\bot}$. Let $\TT = ((x^k,y^k))_{k = 1}^{\ell-2} \in T_{\ell-2}$ 
		be a training set with distinct non-zero x-coordinates such that for all 
		$k = 1, \dots, \ell-2$ 
		we have that $(x^k,y^k) \in \TT$ is such that $y^k = Ax^k$ with $x^k \in \ker(A)^{\bot}$, 
		and such that $A(x^k) \neq \frac{\theta}{4^n}A(w)$ for any $n \in \N$, where 
		$\theta = \min\{\|A\|^{-1}_{op}, 1\}$. If $A$ is computable then $v$, $w$ and $\T_b$ can be made computable. This follows from the facts that an orthonormal basis for $\ker(A)$ (and similarly for $\ker(A)^{\bot}$) (see Theorem 11 in \cite{Ziegler}) and that $\|A\|_{op}$ (see Corollary 21 in \cite{Ziegler}) are computable when \MG{(the fixed matrix with fixed rank)} $A$ is computable.
		
		We now define training sets
		$\{\iota_n^1\}_{n \in \N}$ and $\iota^0$ as follows. Let
		\begin{equation}\label{iota} 
			\iota_n^1 = \T_b \MG{\catop} ((0,0),(\x \!+ \! \frac{\theta}{4^n}w,A(\x \!+ \! \frac{\theta}{4^n }w)))
			\quad \text{and} \quad  
			\iota^0 = \TT\MG{\catop}((0,0),(\x,0)), \quad n \in \N. 
		\end{equation}
		
		In other words, for $\iota_n^1$ we have the following: the first $\ell-2$ elements 
		are equal to the elements in 
		$\T_b$, $(x^{\ell-1},y^{\ell-1}) = (0,0)$, and $(x^{\ell},y^{\ell}) 
		= (\x \!+ \! \frac{\theta}{4^n}w,A(\x \!+ \! \frac{\theta}{4^n }w))$. 
		Similarly, for $\iota^0$ we have the following,
		the first $\ell-2$ elements are equal to the elements in 
		$\T_b$, $(x^{\ell-1},y^{\ell-1}) = (0,0)$, and $(x^{\ell},y^{\ell})=(\MG{\x},0)$.
		\MG{Note that $\|v\|_{\MG{2}} = 2\kappa \leq 6/8 \leq 1$, that $\| \x \!+ \! \frac{\theta}{4^n}w \|_{\MG{2}} \leq \|v\|_{\MG{2}} + \theta \|w\|_{\MG{2}}/4^n \leq 6/8 + \theta/4^n \leq 1$ and similarly that $\|A(\x \!+ \! \frac{\theta}{4^n }w)\|_{\MG{2}} = \theta \|Aw\|_{\MG{2}}/4^n \leq \|A\|_{op}^{-1}\|A\|_{op} \|w\|_2/4^n \leq 1$.}
		We \MG{thus} observe that $\iota^0 \in T_{\ell}$ and $\iota_n^1 \in T_{\ell}$ for all 
		$n \in \N$ and we define the domain $\Omega = \{\iota_n^1\}_{n \in \N} \cup \{\iota^0\}$. \MG{Because $N \geq 3$, at least one of $\ker(A)$ and $\ker(A)^{\bot}$ must have dimension at least $2$. Therefore, because the specific choice of $w$ and $v$ is not important,} it is clear that 
		there exist infinitely many choices for the domain $\Omega$.
		For $\iota \in \Omega$, we define the corresponding inverse problems $(A,\mathcal{M}_1(\iota))$ as follows: 
		$\mathcal{M}_1(\iota) = \pi_1(\iota)$, where $\pi_1(\iota)$ is defined in \MG{\eqref{eq:pi1def}}.
		
		\textbf{Step II} (\emph{Showing that $\Xi(\iota) \neq \emptyset$ for all $\iota \in \Omega$}).
		We show that $\Xi(\iota) \neq \emptyset$ for all $\iota \in \Omega$ as follows. By the fact that $\TT \in \MG{(\ker(A)^{\bot} \times A(\ker(A)^{\bot}))^{\ell-2}}$ with $y = Ax$ for all 
		$(x,y) \in \TT$, the fact that $A(x) \neq \frac{\theta}{4^n}A(w)$ for any $n \in \N$, and the fact 
		that $A$ is injective on $\ker(A)^{\bot}$, 
		it is clear that for all $n \in \mathbb{N}$ the elements in $\iota_n^1$ are distinct in both the 
		x and y-coordinates, or in other words that both $\pi_1(\iota_n^1)$ and $\pi_2(\iota_n^1)$
		each have $\ell$ distinct points. 
		Thus, we can achieve exact interpolation of every point by a smooth neural network according to 
		\cite[Theorem 5.1]{pinkus}. In particular, for each $n \in \mathbb{N},$ there is a neural network 
		$\mathbf{N}_n: \mathbb{R}^m \rightarrow \mathbb{R}^N$  such that for each pair $(\xi_n, \eta_n) \in \iota_n^1$ we have $\mathbf{N}_n(\eta_n) = \xi_n$, and it is clear from Definition 
		\ref{optmap:trainingset} that this is an optimal map for the inverse problem 
		$(A,\mathcal{M}_1(\iota_n^1))$. By the same argument we can find a neural network 
		$\mathbf{N}: \mathbb{R}^m \rightarrow \mathbb{R}^N$ such that for any 
		$(\xi, \eta) \in \TT$ we have $\mathbf{N}(\eta) = \xi$ and $\mathbf{N}(0) = v/2$.  
		We claim that $\mathbf{N}$ is optimal for $(A,\mathcal{M}_1(\iota^0))$. 
		Indeed, note that 
		\MG{\begin{equation}\label{eq:fix1}
			\inf_{\varphi \colon\! \mathcal{M}_2(\iota^0)  \rightrightarrows \R^N} 
			\sup_{x \in \mathcal{M}_1(\iota^0)} d_1^{H}(\varphi(Ax), x) \geq  
			\inf_{\varphi \colon\! \mathcal{M}_2(\iota^0)  \rightrightarrows \R^N} 
			\max_{x =0,x=v} d_1^{H}(\varphi(Ax), x) = \|v\|_{\MG{2}}/2. 
		\end{equation}}
		However, by the definition of $\mathbf{N}$ \MG{(and because $\mathbf{N}$ is single-valued)} we have that 
		\begin{equation}\label{eq:fix2}
		\MG{\sup_{x \in \mathcal{M}_1(\iota^0)} d_1^{H}(\mathbf{N}(Ax), x) }= 	\sup_{x \in \mathcal{M}_1(\iota^0)} \MG{\|\mathbf{N}(Ax) - x\|_{\MG{2}} }
			= \max_{x=0,x=v} \MG{\|\mathbf{N}(Ax) - x\|_{\MG{2}}}  = \|v\|_{\MG{2}}/2,
		\end{equation}
		proving our claim. This means that $\Xi(\iota) \neq \emptyset$ for $\iota \in \Omega$. 
		Thus the computational problem $\{ \Xi, \Omega, \M', \Lambda\}$ is now well defined.
		Next we show that the sequence $\{\iota_n^1\}_{n \in \N}$ and 
		$\iota^0$ satisfy points $a) -b)$ in Proposition \ref{computability_prop}:
		
		\textbf{Step III} (\emph{Proof of \ref{property:MinimisersOfNonZero} from Proposition \ref{computability_prop}}).
		We define 
		\begin{align*}
			S_1 = \bigcup_{n \in \N}\Xi(\iota_n^1) \quad \text{and} \quad 
			\MG{S_2 = \Xi(\iota^0)}.
		\end{align*}
		Using the metric defined in \eqref{breakdown_metric} we get that
		\begin{align*}
			\inf_{\mathbf{N_1} \in S_1,\mathbf{N_2} \in S_2} d_{\M'}(\mathbf{N_1},\mathbf{N_2}) &=
			\inf_{\mathbf{N_1} \in S_1,\mathbf{N_2} \in S_2} \sup_{z \in \bigcap_{\F} \M_2} 
			\| \mathbf{N_1}(z) - \mathbf{N_2}(z)\|_{\MG{2}} \\
			&\geq \inf_{\mathbf{N_1} \in S_1,\mathbf{N_2} \in S_2} \| \mathbf{N_1}(0) -\mathbf{N_2}(0)\|_{\MG{2}} \, \text{ (optimality of $\mathbf{N_j}$ and \eqref{eq:fix1}, \eqref{eq:fix2} give)}\\
			&=\|0 - \frac{1}{2}\x \|_{\MG{2}} = \frac{1}{2}\|\x\|_{\MG{2}} = \frac{1}{2} 2\kappa =  \kappa.
		\end{align*}
		Thus part \ref{property:MinimisersOfNonZero} in Proposition \ref{computability_prop} is satisfied.
		
		\textbf{Step IV} (\emph{Proof of \ref{property:Del1Info} from Proposition \ref{computability_prop}}).
		For each $n \in \N$, the elements in $\TT$ and $(0,0)$ have the same index in 
		$\iota_n^1$ and $\iota^0$. Further, we have that
		$(x^{\ell},y^{\ell}) = (\x \!+ \! \frac{\theta}{4^n}w,A(\x \!+ \! \frac{\theta}{4^n}w))$ in 
		$\iota_n^1$ and that $(x^{\ell},y^{\ell}) = (\x,0)$ in $\iota^0$, or in other words, that the 
		$\ell$'th element in $\iota_n^1$ is equal to 
		$(\x \!+ \! \frac{\theta}{4^n}w,A(\x \!+ \! \frac{\theta}{4^n}w))$ and that the $\ell$'th element 
		in $\iota^0$ is equal to \MG{$(\x,0)$}.
		Thus, we only need to show that the criteria in part (b) holds for each 
		$f_{x,i}^{\ell}$ and
		$f_{y,j}^{\ell}$ for $i=1,\dots,N$, $j = 1,\dots,m$.
		We start with $f_{y,j}^{\ell}$. For each $j$ we get
		$
		|f_{y,j}^{\ell}(\iota_n^1)-f_{y,j}^{\ell}(\iota^0)| \! = \! 
		|A(\frac{\theta}{4^n} w)_j \! - 0| 
		= \frac{\theta}{4^n}|A(w)_j| 
		\leq \frac{1}{4^n},
		$
		for each $n \in \N$. 
		Next, for each 
		$f_{x,i}^{\ell}$ we get
		$
		| f_{x,i}^{\ell}(\iota_n^1) - f_{x,i}^{\ell}(\iota^0)| = |\x_i + 
		\frac{\theta}{4^n}w_i - \x_i |
		\leq \frac{\theta}{4^n}\|w\|_{\MG{2}} 
		\leq \frac{\theta}{4^n} \leq \frac{1}{4^n}, 
		$
		for each 
		$n \in \N$. Thus we can conclude that part (b) in Proposition \ref{computability_prop} holds.
		
		\textbf{Step V} (\emph{The algorithms $\Gamma^0$ and $\Gamma^1$}).
		In the case where $A$ is computable, we have already discussed that each of $\x$, $w$ and $\T_b$ \MG{can be chosen to be} computable. In particular, since $A$ is computable, there exists an algorithm that takes $k$ and $n$ and returns an approximation to $(v+\theta w/4^n), A(v+\theta w/4^n)$ with $\ell_2$-norm error (or any other norm) at most $2^{-k}$. These observations yield the result.
	\end{proof}
	
	In order to prove the existence of algorithms in part (ii) of Theorem \ref{thm_intro:first} 
	we need to obtain a recursively constructable neural network 
	that is, a neural network computable by a Turing machine. 
	For this we will use the following technical lemma.
	
	\begin{lem} \label{lem:radial_network}
		Let $\ell \in \N$ and let $x^1, \dots, x^{\ell} \in \R^N$ and $y^1, \dots, y^{\ell} \in \R^m$ be such that $y^i \neq y^j$ for $j \neq i$. 
		Further, let $\phi: \R \to \R$ be given by the computable function
		$\phi(x) = \frac{1}{x^2+ 1}$, and let $\Phi: \R^m \to \R^{\ell}$ be the function given by 
		\begin{equation}\label{eq:Phi}
			\Phi(y) = 
			[\phi(\| y- y^1\|_{\MG{2}}),  
			\dots,
			\phi(\| y - y^{\ell}\|_{\MG{2}})]^T.
		\end{equation}
		At last, let
		the matrices $X \in \R^{N \times \ell}$ and $R \in \R^{\ell \times \ell}$ be given by
		\begin{equation}\label{eq:theX}
			X = (x^1,  \dots, x^{\ell})
		\end{equation}
		and 
		\begin{align}\label{eq:theR} 
			R = 
			\begin{pmatrix}
				\phi(\| y^{1} - y^1\|_{\MG{2}}) & \dots & \phi(\| y^{\ell} - y^1\|_{\MG{2}}) \\
				\vdots & & \vdots \\
				\phi(\|y^{1} - y^{\ell}\|_{\MG{2}}) & \dots & \phi(\| y^{\ell} - y^{\ell}\|_{\MG{2}})
			\end{pmatrix} 
			= 
			\begin{pmatrix}
				1 & \dots & \phi(\| y^{\ell} - y^1\|_{\MG{2}}) \\
				\vdots & & \vdots \\
				\phi(\|y^{1} - y^{\ell}\|_{\MG{2}}) & \dots & 1
			\end{pmatrix}.
		\end{align}
		\MG{Then $R$ is a symmetric, positive-definite, non-singular matrix.} Additionally, the function $s: \R^m \to \R^N$ given by $s(y) = X  R^{-1} \Phi(y)$ satisfies 
		$s(y^i) = X  R^{-1} \Phi(y^i) = x^i$ for all $i = 1, \dots, \ell$. Moreover, 
		$s \in \NN_{m,N}^{\F}$, that is, $s$ can be 
		represented as a neural network. Finally, \MG{for all $u,v \in \R^m$ we have \begin{equation}\label{eq:phiStability}|\phi(\|u\|_{\MG{2}}) - \phi(\|v\|_{\MG{2}})| \leq |\,\|u\|_{\MG{2}} - \|v\|_{\MG{2}}|
		\end{equation}}		
		 and thus any such $s$ is Lipschitz-continuous. 
	\end{lem}
	\MG{\begin{rem}\label{rem:ConstructedNNInNNmN}
		The claim that $s \in \NN_{m,N}^{\F}$ requires a little more information: the definition of $\NN_{m,N}^{\F}$ requires that the network is bounded on the set $\bigcup_{\mathcal{F}} \M_2$ (see the text surrounding \eqref{metric}). In fact, the network we produce in Lemma \ref{lem:radial_network} is bounded on any bounded set, so provided $\bigcup_{\mathcal{F}} \M_2$ is itself a bounded set (which it is for all cases of interest in this paper), the claim is valid.
	\end{rem}}
	\begin{proof}
		By \cite[Example 5.4]{meshfree} $R$ 
		is a symmetric positive definite nonsingular matrix when the $y^i$'s are unique, which they are by assumption. The fact that $s$ interpolates all the 
		pairs $(x^i,y^i)$ for $i = 1, \dots, \ell$ now follows by observing that 
		\begin{align*}
			s(y^i) = XR^{-1} \Phi(y^i) = X R^{-1}Re_i = Xe_i = x^i,
		\end{align*} 
		where $e_i \in \R^{\ell}$ is the $i$-th standard basis vector.
		It remains to argue that $s$ can be written as a neural network. 
		Indeed, we observe that
		\begin{equation}\label{eq:NN_opt}
			s(y) = \mathbf{N}_{\T}(y) := V_3\rho_2V_2\rho_1V_1(y),
		\end{equation}
		with the affine maps and non-linear functions defined as follows:
		\begin{equation}\label{eq:layers}
			\begin{split}
				V_1&: \,\, W_1 = [1,\hdots, 1]^T_{\ell}  \otimes  I_m \MG{\in \R^{\ell m \times m}}, \quad  b_1 = - \sum_{j=1}^\ell e^\ell_j \otimes \MG{y^j \in \R^{\ell m}},  \quad \rho_1(t) = t^2\\
				V_2 &: \,\, W_2 =  I_\ell \otimes [1,\hdots, 1]_{m} \MG{\in \R^{\ell \times \ell m}}, \quad b_2 = 0 \, \MG{\in \R^{\ell}}, \quad \rho_2(t) = 1/( t + 1), \quad t \in \mathbb{R_+}\\
				V_3 &: \,\, W_3 = X R^{-1} \MG{\in \R^{N \times \ell}}, \quad b_3 = 0 \MG{\in \R^N},
			\end{split}
		\end{equation}
		where $[1,\hdots, 1]_{m}$ is the row vector of length $m$ with ones, $e^\ell_j$ is the $j$-th coordinate column vector of dimension $\ell$, and $I_m$ is the $m$-dimensional identity matrix. The execution of the different layers in \eqref{eq:layers} for an input $y \in \R^m$ can be calculated as follows:
		\[
		y \mapsto z = V_2\rho_1V_1(y) = (\| y-y^1 \|_{\MG{2}}^2, \dots, \| y - y^{\ell} \|_{\MG{2}}^2 )^T \mapsto \rho_2(z)= \Phi(y) \mapsto V_3 \Phi(y) = XR^{-1}\Phi(y).
		\]
		
		The final claim \MG{\eqref{eq:phiStability} regarding $\phi$ follows from}
		\[
		|\phi(x) - \phi(w)|  = \left | \frac{1}{1+x^2} - \frac{1}{1+w^2} \right | = \frac{|w-x| |w+x|}{(1+x^2)(1+w^2)} \leq |w-x| \left(\frac{|w|}{1+|w|^2} + \frac{|x|}{1+|x|^2}\right) \leq |w-x|,
		\]
		where the final inequality follows from e.g. the Arithmetic-Geometric mean inequality applied to obtain $2|w| \leq 1 + |w|^2$. \MG{We thus see} that $s$ is Lipschitz-continuous \MG{since} $XR^{-1}$ is a linear map (and hence Lipschitz). 
	\end{proof}
	\subsection{Formal statement and proof of Theorem \ref{thm_intro:first}} \label{proof1}
	In this section we present the longest proof of the paper. The impossibility results of the theorem are direct consequences of the driving Proposition \ref{prop:driving}, while several of the other results require a substantial amount of work. As previously mentioned, we prove our second main result (Theorem \ref{thm_intro:first}) by proving the following more specific, but technical statement.

	\begin{thm} \label{thm:first:technical}
		Given any $N,m \in \mathbb{N}$ \MG{(with $N \geq 3$)}, a fixed non-zero linear map $A: \R^N \to \R^m$ with non-trivial kernel and an integer $\ell \geq 2$, 
		there exists a domain 
		$\Omega \subset T_{\ell}(A)$ (as described in \eqref{omega}) of training sets $\T$, a set of corresponding initial domains 
		$\{\M_1(\T) \, \vert \, \T \in \Omega\}$ (as described in \eqref{M1_dep}) such that the following occur.
		For the mapping $\Xi:\Omega \rightrightarrows \NN_{m,N}^{\F}$ (as described in \eqref{eq:Xi}), we have $\Xi(\mathcal{T}) \neq \emptyset$ for each $\T \in \Omega$. However, the following happens simultaneously:
		\begin{itemize}[leftmargin=7mm]
			\item[(i)]
			Given any algorithm 
			$\Gamma$ for the computational problem 
			$\{ \Xi, \Omega, \M', \Lambda\}^{\Delta_1}$, where $\M'$ is the metric space defined in \eqref{breakdown_metric}, and any $\delta > 0$, there exist 
			$\mathcal{\tilde T}_1 \in \tilde{\Omega}$,  $\mathcal{T}_1 \in \Omega$ such that \MG{$\tilde \T_1$ corresponds to $\T_1$ and}
			\begin{align}\label{eq:1}
				\MG{\inf_{\mathbf{N}_{opt}^{\M_1(\T_1)} \in \Xi(\T_1)}}\sup_{y \in\bigcap_{\F} \M_2} \| \Gamma(\mathcal{\tilde T}_1)(y) - 
				\mathbf{N}_{opt}^{\M_1(\T_1)}(y) \|_{\MG{2}} \geq 3/16 - \delta.
			\end{align} 
			Moreover, given any randomised algorithm $\Gamma^{\mathrm{ran}}$ for $\{ \Xi, \Omega, \M', \Lambda\}^{\Delta_1}$ \MG{and any $\delta > 0$}, there exist $\tilde{\T_2} \in \tilde{\Omega}$, $\mathcal{T}_2 \in \Omega$ \MG{such that $\tilde \T_2$ corresponds to $\T_2$ and}
			\begin{align}\label{eq:2}
				\mathbf{P}\left(\MG{\inf_{\mathbf{N}_{opt}^{\M_1(\T_2)} \in \Xi(\T_2)}}\sup_{y \in \bigcap_{\F} \M_2} \| \Gamma_{\tilde{\T_2}}^{\mathrm{ran}}(y) - 
				\mathbf{N}_{opt}^{\M_1(\T_2)}(y) \|_{\MG{2}} \geq 3/16 - \delta  
				\right) \geq \frac{1}{2} - \delta.
			\end{align}
			
			\item[(ii)] 
			There exists an infinite sequence of training sets 
			$\{\T^1_i\}_{i \in \N} \subset \Omega$ that are \MG{pairwise distinct} ($\T^1_{k} \neq \T^1_{j}$ when $k \neq j$) and a \MG{(potentially non-halting)} algorithm $\Gamma^1,$ taking inputs in $\tilde \Omega$, that computes\footnote{\MG{Here `computes' is understood in the sense of computing to an arbitrary prescribed accuracy, i.e. $\Gamma^1(\cdot,\epsilon)$ produces a NN that approximates an optimal NN with accuracy (in the $2$-norm, cf. \eqref{metric}) better than $\epsilon$, where $\epsilon > 0$ is an arbitrary rational number.}} an optimal neural network in $\Xi(\T^1_i)$ for each \MG{$i \in \mathbb{N}$ and input $\tilde \T^1_i$ (corresponding to $\T^1_i$)}. However, for any infinite sequence 
			$\{ \T^2_{k} \}_{k \in \N} \subset \Omega$ of \MG{pairwise distinct} training sets and any algorithm $\Gamma^2$ (taking inputs in $\MG{\tilde \Omega}$) that, for any $k \in \N$, produces  -- on input $\tilde \T^2_{k}$ (corresponding to $\T^2_{k}$) -- an approximation to an optimal NN in 
			$\Xi(\T^2_{k})$ \MG{to accuracy (in the $2$-norm, cf. \eqref{metric}) better than 0.1}, we have the following.
			For each $j \in \mathbb{N}$, there exists an integer $i$ and $(x,y) \in \T^2_i$, as well as an element 
			$(x', y') \in \R^N \times \R^m$ with 
			$\| x'\|_{\MG{2}}, 
			\| y' \|_{\MG{2}} \leq 1$ and 
			$\| (x,y) - (x',y')\|_{\MG{2}} \leq \sqrt{2}/4^{j}$, 
			such that if we replace $(x,y)$ with $( x', y')$ then we obtain a new training set 
			$\T^{\prime} = [\T^2_i \setminus \{(x,y)\}] \cup \{( x',  y')\} \in \Omega$ 
			such that 
			\begin{align}\label{eq:failure2}
				\sup_{y \in \M_2(\T^{\prime})}\| \Gamma^2(\tilde \T^{\prime})(y) - \mathbf{N}_{opt}^{\mathcal{M}_1(\T^{\prime})}(y)  \|_{\MG{2}}  > 10^{-1}, 
			\end{align} 
			for some $\tilde \T^{\prime} \in \tilde \Omega$ \MG{which corresponds to $\T^{\prime}$},
			where $\mathbf{N}_{opt}^{\mathcal{M}_1(\T^{\prime})}$ is any optimal neural network for the inverse problem $(A,\M_1(\T^{\prime}))$.  In particular, the failure in \eqref{eq:failure2} applies to $\Gamma^1$ and $\{\T^1_i\}_{i \in \N}$.

			\item[(iii)] Consider the arbitrary sequence $S = \{ \T^2_{k} \}_{k \in \N}$ and algorithm $\Gamma^2$ from (ii). Then, there is a $\T^2  = \{(x^{\prime,l},y^{\prime,l})\}_{l=1}^{\ell} \in \Omega$ such that for any $j$ \MG{there exists $n = n(j) \in \mathbb{N}$ } satisfying 
			\begin{equation}\label{eq:bound1}
				\max_{1 \leq l \leq \ell}\| (x^l,y^l) - (x^{\prime,l},y^{\prime,l})\|_{\MG{2}} \leq \sqrt{2}/4^{j}, \quad \{(x^l,y^l)\}_{l=1}^{\ell} = \T^2_n,
			\end{equation} 
			so that \MG{even if $\Gamma^2$ can compute optimal NNs (to arbitrary accuracy)} for training data in $S$ it cannot \MG{do so} for training data in $S \cup \{\T^2\}$. In \MG{fact},  
			\begin{align}\label{eq:bremnes}
				\sup_{y \in \M_2(\T^2)}\| \Gamma^2(\tilde \T^2)(y) - \mathbf{N}_{opt}^{\mathcal{M}_1(\T^2)}(y)  \|_{\MG{2}}  > 10^{-1}, \, \text{ for some } \tilde \T^2 \in \tilde \Omega \MG{\text{ corresponding to } \T^2},
			\end{align} 
			where $\mathbf{N}_{opt}^{\mathcal{M}_1(\T^2)}$ is any optimal neural network for the inverse problem $(A,\M_1(\T^2))$.

			\item[(iv)] If \MG{every entry of} $A$ is computable, the problem of computing approximations to $\Xi:\Omega \rightrightarrows \NN_{m,N}^{\F}$, to an accuracy of 0.1, is comparable to the Halting problem and is at least as hard.
		\end{itemize}
	\end{thm} 
	
	\begin{proof}[Proof of Theorem \ref{thm:first:technical}]
		Let $\kappa = 3/8$. We start by recalling and using Step I (verbatim) in the proof of Proposition \ref{prop:driving} to define $\Omega$, $\iota^0$ and $\{\iota^1_n\}_{n=1}^{\infty}$.
		The domain $\Omega$ gives rise to a computational problem $\{\Xi, \Omega, \M', \Lambda\}$  on the form 
		described in Proposition \ref{prop:driving}. 
		We proceed by proving each of the points (i) -- (iv), each in turn for $\{\Xi, \Omega, \M', \Lambda\}$ and slight variations of this computational problem.
		
		{\bf Proof of (i).}  
		From Proposition \ref{prop:driving}, it follows that the conditions \ref{property:MinimisersOfNonZero} - \ref{property:Del1Info} in Proposition \ref{computability_prop} are satisfied for 
		$\Omega = \{\iota_n^1\}_{n \in \N} \cup \{\iota^0\}$ with $\kappa = 3/8$. 
		Therefore, by Proposition \ref{computability_prop} we can conclude that 
		$
		\epsilon_B^s \geq \epsilon_{\mathbb{P}hB}^s(p) \geq \epsilon_{\mathbb{P}B}^s(p) \geq \frac{1}{2}\kappa,
		$
		for all $p \in [0,\frac{1}{2})$
		for the
		computational problem $\{\Xi, \Omega, \M', \Lambda\}^{\Delta_1}$.
		Hence, by the definition of the breakdown epsilons it follows that for any generalized algorithm 
		$\Gamma$ and $\delta > 0$ there exists a training set $\tilde \T_{\mathrm{fail}} \in \tilde \Omega$ (recall the definition of 
		$\tilde \Omega$ from \eqref{eq:tildeOmega} in Definition \ref{def:delta1}) \MG{corresponding to some $\T_{\mathrm{fail}} \in \Omega$} such that for every $\MG{\mathbf{N}_{opt}^{\M_1(\T_{\mathrm{fail}})} \in \Xi(\T_{\mathrm{fail}})}$
		\begin{equation}\label{eq:failure1}
			 \sup_{y \in \bigcap_{\F} \M_2} \| \Gamma(\tilde{\T}_{\mathrm{fail}})(y) - \mathbf{N}_{opt}^{\M_1(\T_{\mathrm{fail}})}(y) \|_{\MG{2}} \geq \kappa/2 - \delta \geq 3/16 - \delta,
		\end{equation}
		and for any randomised algorithm $\Gamma^{\mathrm{ran}}$ and any $\delta > 0$ there exists a training set 
		\MG{$\tilde \T_2 \in \tilde \Omega$ corresponding to some $\T_2 \in \Omega$} such that
		\begin{equation*}
			\mathbf{P}\left(\MG{\inf_{\mathbf{N}_{opt}^{\M_1(\T_2)} \in \Xi(\T_2)}}\sup_{y \in \bigcap_{\F} \M_2} 
			\| \Gamma^{\mathrm{ran}}(\tilde{\T_2})(y) - 
			\mathbf{N}_{opt}^{\M_1(\T_2)}(y) \|_{\MG{2}} 
			\geq \kappa/2 - \delta \geq 3/16 - \delta \right) \geq \frac{1}{2} - \delta.
		\end{equation*}

		{\bf Proof of (ii).} 
		We now proceed by proving part (ii) of the theorem in several steps.
		
		{\bf Step I.} \emph{Proving the computational breakdown in part (ii).} 
		Let $\{ \T^2_{k} \}_{k \in \N} \subset \Omega$ be an infinite sequence of different elements and $\Gamma^2$ any algorithm that takes inputs in $\Omega$ and \MG{attempts to} compute an optimal neural network in $\Xi(\T^2_{k})$ for all $k \in \N$. Since the elements in the sequence $S = \{ \T^2_{k} \}_{k \in \N}$ are all different, it follows by the choice of $\Omega$ that 
		$S$ contains elements of the form $\T_b \MG{\catop} ((0,0),(\x \!+ \! \frac{\theta}{4^n}w,A(\x \!+ \! \frac{\theta}{4^n }w)))$ for infinitely many different values of $n$. Hence, for any $j \in \mathbb{N}$ there is an $i \geq j$ such that $\T = \T_b \MG{\catop}((0,0),(\x \!+ \! \frac{\theta}{4^i}w,A(\x \!+ \! \frac{\theta}{4^i }w))) \in S$.
		Let 
		$(x',  y') = (\x,0)$ and $(x,y) = (\x + \frac{\theta}{4^{i}}w, A(\x + \frac{\theta}{4^{i}}w))$.
		It is not hard to see that 
		\begin{align}\label{eq:argue}
			\| ( x', y') - (x,y)\|_{\MG{2}} = \|(\frac{\theta}{4^{i}}w, A(\frac{\theta}{4^{i}}w)) \|_{\MG{2}} 
			\: \leq \: 
			\sqrt{\frac{1}{4^{2i}} + \frac{1}{4^{2i}}} = \frac{\sqrt{2}}{4^{i}} \leq \frac{\sqrt{2}}{4^{j}},
		\end{align}
		and that $\|x' \|_{\MG{2}}, \|y' \|_{\MG{2}} \leq 1$, which yields the bounds stated in (ii). Now let 
		$\T^{\prime} = [\T \setminus \{(x,y)\}] \cup \{(x',y')\}$, and observe that this yields that $\T^{\prime} = \TT \MG{\catop}((0,0),(\x,0)) = \MG{\iota^0}$. Let $\Omega^{\sharp} =\MG{ \{\T^{\prime}\} \cup (S\setminus \{\T\})}$, and observe that $\Omega^{\sharp}$ now contains $\TT \MG{\catop}((0,0),(\x,0))$ and some subsequence of $\{\iota_n^1\}_{n \in \mathbb{N}}$ from \eqref{iota}. By arguing as in the proof of (i), when considering the computational problem $\{\Xi, \Omega^{\sharp}, \M', \Lambda\}$, it follows by Proposition \ref{computability_prop} that $\{\Xi, \Omega^{\sharp}, \M', \Lambda\}^{\Delta_1}$ has the following lower bound on the strong breakdown epsilon: $\epsilon_B^s \geq \frac{1}{2}\kappa$. 
		This implies that \MG{there exists $\T^{\flat} \in \Omega^{\sharp}$ and $\tilde \T^{\flat} \in \tilde \Omega^{\sharp}$ (corresponding to $\T^{\flat}$) so that}
		\begin{equation}\label{eq:sharp}
			\sup_{y \in \M_2(\T^{\flat})}\| \Gamma^2(\tilde \T^{\flat})(y) - 
			\mathbf{N}_{opt}^{\mathcal{M}_1(\T^{\flat})}(y) 
			\|_{\MG{2}} \geq 
			\sup_{y \in \bigcap_{\F} \M_2}\| \Gamma^2(\tilde \T^{\flat})(y) 
			- \mathbf{N}_{opt}^{\mathcal{M}_1(\T^{\flat})}(y) 
			\|_{\MG{2}}  > 10^{-1},
		\end{equation}
		\MG{for every $\mathbf{N}_{opt}^{\mathcal{M}_1(\T^{\flat})} \in \Xi(\T^{\flat})$}. However, $\Gamma^2$ computes an optimal NN \MG{(to accuracy better than $10^{-1}$)} for each $\T \in S$. Thus, the failure input\MG{s} $\T^{\flat} \in \Omega^{\sharp}$ and $\tilde \T^{\flat} \in \tilde \Omega^{\sharp}$ in \eqref{eq:sharp} must be
		$\T^{\flat} = \T^{\prime}$ and $\tilde \T^{\flat} = \tilde \T^{\prime},$ and this finally establishes \eqref{eq:failure2}. 
		
		{\bf Step II.} \emph{Constructing $\Gamma^1$.} 
		We notice that since 
		$\M_1(\iota) = \pi_1(\iota)$ for all $\iota \in \Omega$, we get that any neural network $\mathbf{N} \in \NN_{m,N}^{\F}$ such that 
		$\mathbf{N}(y) = x$ for all $(x,y) \in \iota_n^1$ will be an optimal neural network for 
		$\M_1(\iota_n^1) = \pi_1(\iota_n^1)$ for all $n \in \N$. Thus, we next
		demonstrate how we build a (potentially non-halting) algorithm acting on $\tilde \Omega$, such that if $\tilde \iota \in \tilde \Omega$ corresponds to $\iota^1_n$ for some $n \in \N$ then the algorithm recursively approximates a neural network that interpolates all the points in 
		$\iota_n^1$ for each $n \in \N$. 
		To achieve this we use the radial basis function approach and recursively 
		approximate the interpolating neural network $s: \R^m \to \R^N$ obtained in Lemma \ref{lem:radial_network}. 
		
		We begin by recalling some notation and ideas that will be relevant to this proof. By Definition \ref{def:delta1_infromation} recall that for each $\tilde \iota \in \tilde \Omega$ there is an $\iota = (x^k,y^k)_{k = 1, \dots, \ell} \in \Omega$ as defined in \eqref{def:inexact_rep} such that if we write $\tilde \iota = (\tilde x^k,\tilde y^k)_{k = 1, \dots, \ell} $ then each $(\tilde x^k,\tilde y^k)$ corresponds to $(x^k,y^k)$.
		We will use the notation $\tilde x^{k}(j), \tilde y^{k}(j)$ so that $\tilde x^k = \{ \tilde x^{k}(j)\}_{j=1}^{\infty}$ (similarly $\tilde y^k = \{ \tilde y^{k}(j)\}_{j=1}^{\infty}$) and
		\begin{align} \label{xy_aprox}
			\|\tilde x^{k}(j) - x^k\|_{\MG{2}} \leq 2^{-j} \quad \text{and} \quad  
			\|\tilde y^{k}(j) - y^k\|_{\MG{2}} \leq 2^{-j},
		\end{align}
		for each $k = 1,\dots, \ell$.
		
		Recall the matrices $X,R$ and \MG{the map} $\Phi$ defined in \eqref{eq:theX}, \eqref{eq:theR}, and \eqref{eq:Phi} respectively.
		For a given $\tilde \iota \in \tilde \Omega$ corresponding to $\{(x^i, y^i)\}_{i=1}^{\ell} = \iota \in \Omega$, we consider the matrices \begin{equation}\label{eq:XRj}
			X_j = X(\tilde x^{1}(j),  \dots, \tilde x^{\ell}(j)), \quad R_j = R(\tilde y^{1}(j),  \dots, 
			\tilde y^{\ell}(j))
		\end{equation}
		which approximate the matrices $X = X(x^1,  \dots, x^{\ell})$ and $R = R(y^1,  \dots, y^{\ell})$. Similarly, we define
		\begin{equation}
			\Phi_j(\cdot) = \Phi(\cdot, \tilde y^{1}(j),  \dots, \tilde y^{\ell}(j)). 
		\end{equation}
		and
		$\Phi(\cdot) = \Phi(\cdot, y^1,  \dots, y^{\ell})$ from \eqref{eq:Phi}. With this notation, we can now define the algorithm $\Gamma^1$ taking inputs in $\tilde \Omega$. Broadly speaking, $\Gamma^1$ acts as follows:
		
		\underline{Inputs:} A training set $\tilde \T = \tilde \iota =  ((\tilde x^1,\tilde y^1),(\tilde x^2,\tilde y^2),\dotsc,\MG{(\tilde x^\ell,\tilde y^\ell)})$ corresponding to \\\indent $\iota = \T  = ((x^1,y^1),(x^2,y^2),\dotsc,\MG{(x^\ell,y^\ell)})$ and a dyadic $\epsilon = 2^{-\tilde j}$.

		\underline{Outputs:} A neural network  $N \in \NN_{m,N}^{\F}$.
		\begin{enumerate}[label = \arabic*.]
			\item Attempt to find $k_R \in \mathbb{N}$ so that $R$ has all of its singular values larger than $2^{-k_R}$.
			\item \MG{Set $r \in \mathbb{N}$ so that $2^{-r} \leq \frac{1}{16}\frac{1}{4}\frac{1}{\ell^{2}}2^{-2k_R}2^{-\tilde{j}}$.}
			\item Use the input oracles to read all $\tilde{y} \in \pi_2(\tilde{\iota})$ to \MG{precision $2^{-r}$} and construct the matrix 
			$R_{r}$. Similarly, use the input oracles to read all $\tilde{x} \in \pi_1(\tilde{\iota})$ \MG{to precision $2^{-r}$} to construct the matrix $X_r$.
			\item Return the neural network $X_rR_r^{-1} \Phi_r$\MG{.}
		\end{enumerate}

		We now describe each of these steps of the algorithm in detail. As we proceed, we show that this algorithm $\Gamma^1$ has the property that 
		\MG{for each }dyadic $\epsilon > 0$ and $n \in \N$, we have that 
		$
		\sup_{y \in \bigcup_{\F} \M_2} \| \Gamma^1(\tilde{\iota}_n^1, \epsilon)(y) - \mathbf{N}_{\iota_n^1}(y) \|_2 \leq \epsilon
		$, for all $n \in \mathbb{N}$, where $\mathbf{N}_{\iota_n^1}(y)$ is the interpolating neural network defined in \eqref{eq:NN_opt}. Note that this algorithm is defined on $\tilde \iota$ that correspond to $\iota^0$, but as we will see it may not terminate and we do not analyse its behaviour on these inputs.
		
		Define $k : \mathbb{N} \rightarrow \mathbb{N}$ by
		\begin{equation} \label{k_bound}
			\begin{split}
				k(j) &:= \min\big\{\mu \in \mathbb{N} \, \vert \, \forall r \geq \mu, \, \, R_r \, \text{is invertible and}\\
				& \qquad \qquad \qquad \: \sup_{y \in \bigcup_{\F} \M_2}
				\|  X_r R_r^{-1} \Phi_r(y) -X R^{-1} \Phi(y) \|_{\MG{2}} \leq 2^{-j} \big\}.
			\end{split}
		\end{equation}
		We note that $k$ is well defined whenever $\tilde \iota$ corresponds to $\iota^1_n$ \MG{for some $n \in \mathbb{N}$}: Firstly, since \MG{all the points} in $\pi_2(\iota_n^1)$ are distinct, there must exist a finite $\mu \in \N$ such that 
		$\tilde y^{i}(\mu) \neq \tilde y^{j}(\mu)$ for all $y^i,y^j \in \pi_2(\iota_n^1)$ with $i \neq j$, making 
		$R_r$ invertible for all $r \geq \mu$ by 
		Lemma \ref{lem:radial_network}.
		Secondly, from the definitions in \eqref{eq:theX} and \eqref{eq:theR}, 
		$X_j \rightarrow X$, $R_j \rightarrow R$. \MG{Similarly, using Lemma \ref{lem:radial_network} and equations \eqref{eq:Phi} and \eqref{eq:phiStability}, we have that for $i \in\{1,2,\dotsc,\ell\}$,
		\begin{align*}|(\Phi(y))_i - (\Phi_j(y))_i| &\leq  |\phi(\|y - y^i\|_{\MG{2}}) - \phi(\|y - \tilde{y}^i(j)\|_{\MG{2}})| \\&\leq |\,\|y - y^i\|_{\MG{2}} - \|y - \tilde{y}^i(j)\|_{\MG{2}}| \leq \|y^i-\tilde{y}^i(j)\|_{\MG{2}} \leq 2^{-j}
		\end{align*}
		so that
		$\Phi_j \rightarrow \Phi$ uniformly} as $j \rightarrow \infty$. \MG{We conclude that $k$ is well-defined.}
		
		We will \MG{work towards estimating} an upper bound for $k$, defined in \eqref{k_bound}, that can be recursively 
		(computable by a Turing machine) computed from $\{(\tilde x^i, \tilde y^i)\}_{i=1}^{\ell} $. In order to do that we need to start by recursively finding a lower bound for the smallest singular value 
		$\beta_R$ of $R$ and an upper bound for the largest singular value $\beta_X$ of $X$. The upper bound for $X$ can be derived in a simple  way. Indeed, by using the fact that 
		$\| X \|_{op} = \max_{\beta \in \mathrm{Sp}(X^*X)} \sqrt{\beta}$ \MG{and that the operator norm of a matrix is bounded from above by its Frobenius norm} we can observe that
		\begin{align} \label{X_bound}
			\sqrt{\beta} \leq \|X\|_{op} \leq 
			\sqrt{\sum_{i = 1}^{\ell} \| x^i \|_{\MG{2}}^2} 
			\leq \sqrt{\sum_{i = 1}^{\ell} 1} = \sqrt{\ell} 
			< \ell \quad \text{for all} \; 
			\beta \in \mathrm{Sp}(X^*X),  
		\end{align}
		where $\mathrm{Sp}(X^*X)$ denotes the spectrum of the operator $X^*X$.
		For the matrix $R$ the task is slightly more complicated, and we need to use sufficient approximations of $R$ in order to achieve the task of finding 
		a lower bound for the smallest singular value.  
		More specifically, we need to find a lower bound for $2^{-h(j)}$ where 
		$h: \N \to \N$ is given by
		\begin{align}\label{eq:defining_h}
			h(j) := \min\{\mu \in \mathbb{N} \, \vert \, \forall r \geq \mu \, , 
			\|R_r - R \|_{op} \leq 2^{-j}  \}.
		\end{align}
		
		\textbf{Claim:} 
		We claim that for $j,r \in \N$ we have that
		\MG{\begin{equation}\label{r_bound}
			2^{-r} \leq \frac{1}{2}\frac{1}{\ell}2^{-{j}} \, \Rightarrow \, h(j)\leq r,
		\end{equation}}
		where $h$ is defined in \eqref{eq:defining_h}. Indeed, we verify this by the following series of calculations. We first observe that
			$\|R_{r} - R \|_{\MG{2}}  = \| [\phi( \|  \tilde y^{j}(r) -  \tilde y^{i}(r) \|_{\MG{2}}) 
			- \phi( \|y^j - y^i \|_{\MG{2}})]_{i,j = 1}^{i,j = \ell} \|_{\MG{2}}$
		where 
		$ \tilde y^{j}(r)$ and 
		$\tilde y^{i}(r)$ are the approximations of $y^j$ and $y^i$
		to accuracy $2^{-r}$ as defined in \eqref{xy_aprox} and where the $\ell_2$-norm of a matrix should be interpreted as identical to the Frobenius norm.
		\MG{By Lemma \ref{lem:radial_network}, for any $i,j \in \{1,2,\dotsc,\ell\}$}, \MG{\begin{align*}\|\phi( \|  \tilde y^{j}(r) -  \tilde y^{i}(r) \|_{\MG{2}}) 
		- \phi( \|y^j - y^i \|_{\MG{2}})\| & \leq |\|  \tilde y^{j}(r) -  \tilde y^{i}(r) \|_{\MG{2}} - \|y^j - y^i \|_{\MG{2}}| \\
		&\leq \|\tilde y^{j}(r) - y^j\|_{\MG{2}} + \|\tilde y^{i}(r) - y^i\|_{\MG{2}} \leq 2\cdot 2^{-r},
		\end{align*}}
		and thus
		\MG{\begin{align*}
			\| R_{r} - R \|_{op} \leq \| R_{r} - R \|_{\MG{2}} = \sqrt{\sum_{i,j} [R_{r} - R]^2_{i,j}} \leq 
			\sqrt{\sum_{i,j} (\| \tilde y^{j}(r) - y^j \|_{\MG{2}} + 
				\| \tilde y^{i}(r) - y^i \|_{\MG{2}})^2}.
		\end{align*}}
		At last, since \MG{$2^{-r} \leq \frac{1}{2}\frac{1}{\ell}2^{-j}$}, we can conclude that
	\MG{	$
			\| R_{r} - R \|_{op} \leq  \sqrt{\ell^2(\frac{1}{2}\frac{1}{\ell}2^{-j}  \! 
				+ \! \frac{1}{2}\frac{1}{\ell}2^{-j})^2} 
			= \! \sqrt{\ell^2 \frac{1}{\ell^2} 2^{-2j}} =  2^{-j},
		$}
		which proves our claim. 
		
		Using this, we can now present a recursive algorithm that finds a lower bound for the smallest singular value of $R$. Note that $R$ is self-adjoint and positive definite (c.f. Lemma \ref{lem:radial_network}), thus the smallest singular value coincides with the smallest element in the spectrum of $R$:
		
		\vspace{1mm}
		
		\underline{Inputs:} Oracles for all $\tilde{y} \in \pi_2(\tilde{\iota})$.
		
		\underline{Outputs:} A lower bound for the \MG{spectral} values of $R$, assuming that $\tilde \iota$ corresponds to $\iota^1_n$ for some $n \in \N$. More precisely, some
		$k_R \in \N$ such that $2^{-k_R} \leq \beta_R$ for all $\beta_R \in Sp(R)$.
		\begin{itemize}
			\item[1.] Put $j = 1$.
			\item[2.] \MG{Choose $r \in \N$ such that $2^{-r} \leq \frac{1}{2}\frac{1}{\ell}2^{-(j+1)}$.}
			\item[3.] Use the input oracles to read all $\tilde{y} \in \pi_2(\tilde{\iota}_n^1)$ to \MG{precision $2^{-r}$} and construct the matrix 
			$R_{r}$. As we have seen, the criterion in step 2 as well as equations \eqref{eq:defining_h} and \eqref{r_bound} implies that 
			$\|R_{r} - R \|_{op}  \leq 2^{-(j+1)}$.
			\item[4.] Check -- by trying to compute the Cholesky decomposition -- whether the matrix $R_{r} - 2^{-j}I$ is positive definite \MG{(note that this is possible because $R_r - 2^{-j}I$ contains rational entries).}
			\begin{itemize}
				\item[a.] If it is, then we have found a lower bound for 
				the smallest spectral value of $R$, since then $\beta > 2^{-j}$ for all 
				$\beta \in \text{Sp}(R_r)$, and thus (using that $d_H(\text{Sp}(R_r),\text{Sp}(R)) \leq \|R - R_r\|_{op} \leq 2^{-(j+1)}$ \cite{Haagerup} (p. 723)\MG{)}, where $d_H$ is the Hausdorff metric,
				we get $\gamma > 2^{-(j+1)}$ for all $\gamma \in \text{Sp}(R)$. 
				Here $\text{Sp}(R)$ denotes the spectrum of $R$, 
				which equals the set of singular values of $R$ (since $R$ is symmetric and positive definite). Set $k_R = j+1$ and return $2^{-k_R}$. 
				\item[b.] If it is not, then increase $j$ by one and repeat steps 2, 3 and 4.\\
			\end{itemize} 
		\end{itemize}
		We notice that the above process will always terminate after a finite number of steps assuming that $\tilde \iota$ corresponds to $\iota^1_n$ for some $n$, since 
		$R$ is finite dimensional and positive definite
		by Lemma \ref{lem:radial_network}. \MG{Now that we have an algorithm that produces $k_R$ with $2^{-k_R} \leq \beta_R$ for all $\beta_R \in Sp(R)$, (in particular, $\|R^{-1}\|_{op} \leq 2^{k_R}$)}, we are ready to present a lower bound for $2^{-k(j)}$, where 
		$k$ is as defined in \eqref{k_bound}.
		
		\textbf{Claim:} We claim that for any fixed $j \in \N$, if $r \in \mathbb{N}$ is such that
		\MG{\begin{align} \label{strong_bound}
			2^{-r} \leq \frac{1}{16}\frac{1}{4}\frac{1}{\ell^{2}}2^{-2k_R}2^{-j}
		\end{align}}
		with $k_R$ set to be the output from the algorithm above, then $2^{-r}$ is a lower bound for $2^{-k(j)}$. In particular, note for such $r$ we must have 
		\begin{equation} \label{weak_bound}
			2^{-r} \leq \frac{1}{7}\frac{1}{\ell\sqrt{\ell}}2^{-k_R}2^{-j}\MG{.}
		\end{equation} 
		
		We now proceed with verifying this claim. \MG{We start by noting that our choice of $r$ in \eqref{strong_bound} implies that
		\begin{align*}
			2^{-r} \leq \frac{1}{16}\frac{1}{4}\frac{1}{\ell^{2}}2^{-2k_R}2^{-j}
			= \frac{1}{2}\frac{1}{\ell}2^{-(2k_R +j+5 + \log_2(\ell))} 
			\leq \frac{1}{2}\frac{1}{\ell}2^{-(2k_R +j+5 + \lfloor \log_2(\ell) \rfloor )}. 
		\end{align*}}
		\MG{By \eqref{eq:defining_h} and \eqref{r_bound}, and because 
		$5 + \lfloor \log_2(\ell) \rfloor \geq 4 + \log_2(\ell)$, this implies that
		\begin{equation}\label{eq:Rr-ROpBound}
			\|R_r - R\|_{op} \leq 2^{-(2k_R + j + 5 + \lfloor\log_2(\ell) \rfloor)}
			\leq 2^{-(2k_R + j + 4 + \log_2(\ell))} \leq 2^{-2k_R}2^{-j}\frac{1}{16}\frac{1}{\ell}.
		\end{equation}}
		\MG{Immediately, we see that this implies that $R_r$ is invertible.} Indeed, since $R$ is invertible and \MG{since}
		$\MG{\| R - R_r \|_{op}} \: < \: 2^{-k_R} \leq \MG{1/\| R^{-1}\|_{op}} $ we must have that 
		$R_r$ is also invertible, and we get that
		\begin{align*}
			\sup_{y \in \bigcup_{\F} \M_2} &\| X_{r}R_{r}^{-1}\Phi_{r}(y)\! -\! XR^{-1}\Phi(y)\|_{\MG{2}} \\
			&= \sup_{y\in \bigcup_{\F} \M_2}\| X_{r}R_{r}^{-1}\Phi_{r}(y) - XR^{-1}\Phi_{r}(y) + XR^{-1}\Phi_{r}(y) - XR^{-1}\Phi(y)\|_{\MG{2}} \\
			&\leq  \sup_{y \in \bigcup_{\F} \M_2}\| X_{r}R_{r}^{-1} \! - \! XR^{-1} \|_{op} \| \Phi_{r}(y)\|_{\MG{2}} \! + \! \|R^{-1}\|_{op}\|X\|_{op} 
			\|\Phi_{r}(y) \!- \Phi(y)\|_{\MG{2}}.
		\end{align*}
		Next, we observe that with the approximations
		$\tilde{y}^{r}$ as defined in \eqref{xy_aprox} 
		with $\|\tilde{y}^{i}(r) - y^i\|_{\MG{2}} \leq 2^{-r}$ for $i\in \{1,2,\dotsc,\MG{\ell}\}$, we get that \MG{(using \eqref{eq:phiStability})}:
		\begin{align}
			\sup_{y \in \bigcup_{\F} \M_2}\| \Phi_{r}(y) - \Phi(y) \|_{\MG{2}}\leq \sqrt{\ell}
			\sup_{i\in \{1,2,\dotsc,\MG{\ell}\}} \|\tilde y^{i}(r) - y^i \|_{\MG{2}} \leq  \sqrt{\ell} 2^{-r} 
			\leq \frac{1}{7}\frac{1}{\ell}2^{-k_{R}}2^{-j}\label{eq:Phir-PhiOpBound}.
		\end{align}
		Thus, using \eqref{X_bound},
		\begin{align*}\sup_{y \in \bigcup_{\F} \M_2}\|R^{-1} \|_{op} \| X \|_{op} \| \Phi_{r}(y) - \Phi(y) \|_{\MG{2}} \leq 2^{k_R}\ell 
		\sup_{y \in \bigcup_{\F} \M_2}\| \Phi_{r}(y) - \Phi(y) \|_{\MG{2}} &\leq 2^{k_R} \ell\frac{1}{\ell}2^{-k_R}\frac{1}{7}2^{-j} \leq \frac{1}{7}2^{-j},
		\end{align*}
		by our choice of $r$.
		To ease notation in further calculations, we fix $y$ to be an arbitrary element in $\bigcup_{\F} \M_2$, and we can continue our calculations as follows:
		\begin{equation}\label{tree_terms}
			\begin{split}
				\| X_{r}R_{r}^{-1}&\Phi_{r}(y)\! -\! X R^{-1}\Phi(y)\|_{\MG{2}} 
				\leq \|X_{r}R^{-1}_{r} - XR^{-1} \|_{op} \|\Phi_{r}(y) \|_{\MG{2}}  + \frac{1}{7}2^{-j} \\
				&= \|X_{r}R^{-1}_{r} - X_{r}R^{-1} + X_{r}R^{-1} - XR^{-1} \|_{op} \|\Phi_{r}(y) \|_{\MG{2}}  + \frac{1}{7}2^{-j} \\
				&\leq (\|R^{-1}_{r} - R^{-1} \|_{op} \| X_{r}\|_{op} + \|R^{-1}\|_{op}\|X_{r} - X \|_{op}) \|\Phi_{r}(y) \|_{\MG{2}} + \frac{1}{7}2^{-j} \\
				&= \|R^{-1}_{r} - R^{-1} \|_{op} \| X_{r}\|_{op}\|\Phi_{r}(y)\|_{\MG{2}} + \|R^{-1}\|_{op}\|X_{r} - X \|_{op}
				\|\Phi_{r}(y) \|_{\MG{2}} + \frac{1}{7}2^{-j}.
			\end{split}
		\end{equation}
		The formula above consists of three terms, where we need to calculate bounds for 
		the first two. 
		For the sake of structure and overview we consider each of the terms separately. We start with the term 
		$\|R^{-1}\|_{op}\|X_{r} - X \|_{op} \|\Phi_{r}(y) \|_{\MG{2}}$ and we notice that 
		$\| \Phi\MG{_r}(y) \|^2_{\MG{2}} = \sum_{1 \leq i \leq \ell} 1/(1 + \| y -\MG{ \tilde{y}^i(r)} \|_{\MG{2}}^2)^2 
		\leq \sum_{1 \leq i \leq \ell} 1= \ell$, thus 
		$\|\Phi_\MG{r}(y) \|_{\MG{2}} \leq \sqrt{\ell}$. Therefore
		\begin{align} \label{second_b}
			\|R^{-1}\|_{op}\|X_{r} - X \|_{op} \|\Phi_{r}(y) \|_{\MG{2}}
			\leq \|R^{-1} \|_{op} \| X - X_{r} \|_{op}\sqrt{\ell}
			\leq \frac{1}{7}2^{-j},
		\end{align} 
		where the last inequality follows by our choice of $r$, more precisely that 
		$\| \tilde x^{i}(r) - x^i \|_{\MG{2}} \leq 
		2^{-r} \leq \frac{1}{7}\frac{1}{\ell}2^{-k_R}2^{-j}$ for $i \in \{1,2,\dotsc,\MG{\ell}\}$
		where $\tilde x^{i}(r)$ is the approximation of $x^i$ 
		to accuracy $2^{-r}$ as defined in \eqref{xy_aprox},
		and that 
		$\|X_{r} - X\|_{op} \leq \sqrt{\ell} \max_{i \in \{1,2,\dotsc,\MG{\ell}\}} \| \tilde x^{i}(r) - x^i \|_{\MG{2}} \leq \sqrt{\ell} \cdot 2^{-r}$.
		
		We now move on to the second term 
		$\| R_{r}^{-1} - R^{-1} \|_{op} \| X_{r} \|_{op} \| \Phi_{r}(y) \|_{\MG{2}}$ of \eqref{tree_terms}. By the second resolvent identity we get that
		\begin{align*}
			\| R_{r}^{-1} - R^{-1} \|_{op} \| X_{r} \|_{op} \| \Phi_{r}(y) \|_{\MG{2}} \leq \|R_{r}^{-1}\|_{op} \|R^{-1}\|_{op} \|R_{r} - R 
			\|_{op} \| X_{r} \|_{op} \|\Phi_{r}(y) \|_{\MG{2}}.
		\end{align*}
		In order to continue our calculations we need to obtain an upper bound for $\|R^{-1}_{r} \|_{op}$ and $\|X_{r} \|_{op}$. For $\|X_{r}\|_{op}$ 
		this is fairly straightforward and we see that 
		\begin{equation}\|X_{r}\|_{op} = \|X_{r} - X + X \|_{op}  \leq \| X_{r} - X \|_{op} + \| X \|_{op} \leq 2^{-j} + \sqrt{\ell} \leq   1 + \sqrt{\ell}.\label{eq:XrBound}\end{equation}
		
		Deriving the upper bound for $\|R^{-1}_{r}\|_{op}$ is slightly more complicated. We start by a similar calculation as for $\|X_{r}\|_{op}$ and get that 
		$\|R_{r}^{-1} \|_{op} = \|R_{r}^{-1} - R^{-1} + R^{-1} \|_{op} \leq \|R_{r}^{-1} - R^{-1} \|_{op} + \| R^{-1} \|_{op}$. We continue our calculations by using the second resolvent identity:
		This implies that
		$\|R^{-1}_{r} \|_{op} \leq \|R^{-1}\|_{op} \| R^{-1}_{r} \|_{op} \| R_{r} - R \|_{op} + \|R^{-1} \|_{op}$ and hence
		$\|R^{-1}_{r} \|_{op} (1 - \|R^{-1}\|_{op} \|R_{r} - R \|_{op})  \leq \|R^{-1}\|_{op}$.
	    \MG{From \eqref{eq:Rr-ROpBound} and $\|R^{-1}\|_{op} \leq 2^{k_R}$ }we obtain that $\|R^{-1} \|_{op} \| R_{r} - R \|_{op} \leq 2^{-j}$.
		Thus, we get that $\|R^{-1}_{r} \|_{op}(1 - 2^{-j}) \leq \|R^{-1}_{r} \|_{op} (1 - \|R^{-1}\|_{op} \|R_{r} - R \|_{op}) \leq \|R^{-1}\|_{op}$, which gives us that $\|R^{-1}_{r} \|_{op} \leq \|R^{-1}\|_{op}/(1-2^{-j}) \leq 2^{k_R}/(1-2^{-j})$. At last, since $j \geq 1$ we get that 
		\begin{equation} \label{eq:RInvrBound}\|R^{-1}_{r} \|_{op} \leq 2^{k_R}/(1-2^{-j}) \leq 2\cdot 2^{k_R}.\end{equation}
		\MG{Using \eqref{eq:Rr-ROpBound}, \eqref{eq:XrBound}, \eqref{eq:RInvrBound} and 
		$\|\Phi_r(y)\|_{\MG{2}} \leq \sqrt{\ell} \leq 1+\sqrt{\ell}$, we can finish our calculations as follows:}
		\begin{align} \label{third_b}
			\|&R_{r}^{-1}\|_{op} \|R^{-1}\|_{op} \|R_{r} - R \|_{op} \| X_{r} \|_{op} \|\Phi_{r}(y) \|_{\MG{2}} \nonumber\\
			&\leq 2 \cdot 2^{k_R} \cdot 2^{k_R} \|R_{r}- R \|_{op} (1 + \sqrt{\ell})^2 \leq 2 \cdot 2^{2k_R} (2\sqrt{\MG{\ell}})^2 \|R_{r}- R \|_{op} \nonumber \\
			&\leq 2 \cdot 2^{2k_R}\cdot 4\ell \|R_{r}- R \|_{op} 
			\leq 2 \cdot 2^{2k_R} \cdot 4\ell \cdot  2^{-2k_R}\frac{1}{16}\frac{1}{\ell} 2^{-j}
			\leq \frac{1}{2} 2^{-j}
			\leq \frac{4}{7}2^{-j}. 
		\end{align}
		
		Finally, by combining the bounds of each of the terms obtained in \eqref{second_b} and 
		\eqref{third_b} with the estimate in \eqref{tree_terms}, 
		we can conclude that
		\begin{align*}
			\|  X_{r}R_{r}^{-1} \Phi_{r}(y) - XR^{-1}\Phi(y) \|_{\MG{2}} \leq \MG{6} \cdot \frac{1}{7}2^{-j}  
			\MG{\leq} 2^{-j}
		\end{align*}
		for arbitrary $y \in \bigcup_{\F} \M_2$, where $r$ satisfies the bound specified in 
		\eqref{strong_bound}. Since 
		$y \in \bigcup_{\F} \M_2$ was chosen arbitrarily we can conclude that the inequality above also holds true when we take the supremum over all 
		$y \in \bigcup_{\F} \M_2$, \MG{thus proving the claim.}

		At last, we observe that \MG{(using Lemma \ref{lem:radial_network} and the argument given in Remark \ref{rem:ConstructedNNInNNmN})} $X_{r} \circ  R^{-1}_{r} \circ \Phi_r \in \NN_{m,N}^{\F}$. We can conclude that part (ii) of the result holds. 		
		
		{\bf Proof of (iii).} The proof follows \MG{Step I} of the proof of (ii) closely.
		Let
		$
		\T^2 = \TT \MG{\catop} ((0,0),(\x,0)) \in \Omega. 
		$
		Note that since $S = \{ \T^2_{k} \}_{k \in \N} \subset \Omega$ consists of infinitely many distinct training sets, it follows, by the choice of $\Omega$, that 
		$S$ contains elements of the form $\T_b \MG{\catop} ((0,0),(\x \!+ \! \frac{\theta}{4^n}\MG{w},A(\x \!+ \! \frac{\theta}{4^n }\MG{w})))$ for infinitely many different values of $n$. Hence, for any $j \in \mathbb{N}$ there is an $i \geq j$ such that $\T_b \MG{\catop} ((0,0),(\x \!+ \! \frac{\theta}{4^i}\MG{w},A(\x \!+ \! \frac{\theta}{4^i }\MG{w}))) \in S$. Thus, by arguing as in \eqref{eq:argue}, we obtain \eqref{eq:bound1}. Finally, to see \eqref{eq:bremnes}, let $\Omega^{\flat} = S \cup \{\T^2\}$ and observe that $\Omega^{\flat}$ now contains $\TT \MG{\catop} ((0,0),(\x,0))$ and some subsequence of $\{\iota_n^1\}_{n \in \mathbb{N}}$ from \eqref{iota}. The rest of the argument is identical to the proof of \eqref{eq:failure2} in (ii) after replacing the computational problem $\{\Xi, \Omega^{\sharp}, \M', \Lambda\}$ in the proof of (ii) with $\{\Xi, \Omega^{\flat}, \M', \Lambda\}$.

		{\bf Proof of (iv).} Our aim is to show that given the computational problem 
		$\MG{\{\Xi,  \Omega, \M, \Lambda\}^{\Delta_1}}$ and the existence of a Turing machine (or Markov algorithm), that can compute approximations to $\Xi$ to accuracy 0.1, one can recursively transform this machine/algorithm to another Turing machine that can decide the halting problem. Let
		\begin{align*}
			\hat{\Omega} =  \{(T,\xi) \: | \: T \: \text{is a Turing machine}, \: \xi \: \text{is an input to } 
			T \},
		\end{align*}
		and let $\hat{\Xi}(T,\xi) = 1$ if $T$ halts on $\xi$ and $\hat{\Xi}( T,\xi) = 0$ if $T$ does not halt on $\xi$. We want to show that we can build a recursive mapping $\hat{\Gamma}$ such that $\hat{\Gamma}(T,\xi) = \hat{\Xi}(T,\xi)$.  
		We recall that $\{ \iota_r^1\}_{r \in \N} \cup \{\iota^0\}_{r \in \N} \: \subset \Omega$, with 
		$\iota^0 = \TT \MG{\catop}((0,0),(\x,0))$, and  
		$\iota_{r}^1 = \TT\MG{\catop}((0,0),(\frac{\theta}{4^{r}}w + \x, A(\frac{\theta}{4^{r}}w + \x)))$
		for $r \in \N$. Now, since $A \in \R^{m \times N}$ \MG{is fixed and} by assumption is computable, then (by Proposition \ref{prop:driving}) $\iota^0$ and 
		$\{\iota_{r}^1\}_{r \in \mathbb{N}}$ can be chosen to be computable -- as a sequence (that is, there is a Turing machine that can compute arbitrary approximations to all elements in the sequence).
		Hence, for any 
		\[
		f \in \Lambda =  \{f^k_{y,j},f^k_{x,i}: \Omega \to \R \: | \: i = 1, \dots, N \;, \; j = 1, \dots, m, \; k = 1, \dots, \ell\},
		\]	
		which we recall from \eqref{measurements}, and any $n \in \mathbb{N}$, there is a function $f_n: \Omega \to \mathbb{D}_n$ satisfying \eqref{inexact_input} so that $f_n$ provides a $2^{-n}$-approximation to $f$, and such that the mappings 
		\begin{equation}\label{eq:recursivef}
			\mathbb{N}^2 \ni (n,n') \mapsto f_n(\iota_{n'}^1), \quad \mathbb{N} \ni n \mapsto f_n(\iota^0) \, \text{ are recursive.}
		\end{equation}
		We will denote the approximation functions $f_n$, as above, corresponding to $f^k_{y,j},f^k_{x,i} \in \Lambda$ by \MG{$f^k_{y,j,n}$, $f^k_{x,i,n}$}. Now define the mapping 
		\MG{$g: \hat{\Omega}\times \N \times \{1,2,\dotsc,\ell\} \to \mathbb{D}^{m + N}$} \MG{(where we recall that $\mathbb{D} = \cup_{n=1}^{\infty} \mathbb{D}_n$ is the collection of dyadic numbers)} as follows. Let
		\MG{$g(T,\xi,n,k) = (\tilde x^k(n),\tilde y^k(n))$}, where 
		
		\MG{\begin{equation}\label{eq:def_iota}
				(\tilde x^k(n),\tilde y^k(n)) =
				\begin{cases}
					(\{f^k_{x,i,n+1}(\iota^{0})\}_{i=1}^N,
					\{f^k_{y,j,n+1}(\iota^{0})\}_{j=1}^m)
					&  \! \! \text{if $T$ does not halt on $\xi$ after $n$ steps,} \\
					(\{f^k_{x,i,n}(\iota^{1}_{n'}) \}_{i = 1}^N,
					\{f^k_{y,j,n}(\iota^{1}_{n'}) \}_{j = 1}^m) &  \! \! \text{if $T$ halts on $\xi$ after $n'$ steps, with $n' \: \leq \: n$.}
				\end{cases}
		\end{equation}}
	We now claim that $\tilde \iota = \tilde \iota(T,\xi) \coloneqq
		\MG{(\{(\tilde x^k(n),\tilde y^k(n))\}_{n=1}^{\infty})_{k=1}^{\ell}}$ is an element in $\tilde \Omega$. Indeed, we verify this claim 
		in two parts:
		
		\textbf{Case 1:} There exists an $n'\in \N$ such that $T$ halts on $\xi$ after $n'$ steps.
		We then claim that  $\tilde \iota = \tilde \iota(T,\xi)$ \MG{corresponds to} $\iota_{n'}^1 = \MG{(x^k,y^k)_{k=1}^{\ell}:=} \MG{\TT\catop((0,0),(\frac{\theta}{4^{n'}}\MG{w} + \x, 
		A(\frac{\theta}{4^{n'}}\MG{w} + \x)))}\in \Omega$ in the sense of \MG{\eqref{eq:linfestimateDelta1}} \MG{(where, of course, $x^k$ and $y^k$ depend on $n'$).} Indeed, for $n < n'$ we have that \MG{$(\tilde x^k(n),\tilde y^k(n)) =
		(\{f^k_{x,i,n+1}(\iota^{0})\}_{i=1}^N,
		\{f^k_{y,j,n+1}(\iota^{0})\}_{j=1}^m)$.} \MG{Thus, according to 
		\eqref{eq:linfestimateDelta1}, we need to check that $\|x^k - \tilde x^k(n)\|_{\infty}, \|y^k - y^k(n)\|_{\infty} \leq 2^{-n}$ for $k=1,2,\dotsc,\ell$ and $n \in \mathbb{N}$, which can be equivalently stated as the condition} $|f_{x,i,n+1}^k(\iota^0) - f_{x,i}^k(\iota_{n'}^1)| 
		\leq 2^{-n}$ and that $|f_{y,j,n+1}^k(\iota^0) - f_{y,j}^k(\iota_{n'}^1)| \leq 2^{-n}$ for each 
		$i = 1, \dots, N, j = 1, \dots, m$, $k = 1, \dots, \ell$, and \MG{$n < n'$}. 
		We start by verifying the claim for
		the $f_{x,i}^k$'s, and observe that 
		\begin{align*}
			|f_{x,i,n+1}^k(\iota^0) - f_{x,i}^k(\iota_{n'}^1)| 
			&\leq |f_{x,i,n+1}^k(\iota^0)- f_{x,i}^k(\iota^0)|+|f_{x,i}^k(\iota^0)- f_{x,i}^k(\iota_{n'}^1)| \\
			&\leq \frac{1}{2}\cdot 2^{-n} + 4^{-n} \leq 2^{-n},
		\end{align*}
		where the second to last inequality follows from the fact (by Proposition \ref{prop:driving}) that $\Omega$ satisfies conditions 
		(a)-(b) in Proposition \ref{computability_prop} . Verifying the claim for the $f_{y,j}^k$'s can be done in exactly the same way. For $n \geq n'$ \MG{and $k \in \{1,2,\dotsc,\ell\}$ we have that} \MG{$(\tilde x^k(n),\tilde y^k(n)) = (\{f^k_{x,i,n}(\iota^{1}_{n'})\}_{i = 1}^N,
			\{f^k_{y,j,n}(\iota^{1}_{n'}) \}_{j = 1}^m)$}, from which it immediately follows 
		that \MG{\eqref{eq:linfestimateDelta1}} is satisfied, by the definition of the $f^k_{x,i,n}$'s and 
		$f^k_{y,j,n}$'s.

		\textbf{Case 2:} T does not halt on $\xi$. We then claim that  $\tilde \iota 
		= \tilde \iota(T,\xi) = \MG{(\{(\tilde x^k(n),\tilde y^k(n))\}_{n=1}^{\infty})_{k=1}^{\ell}}$ \MG{corresponds to} 
		$\iota^0 = \TT \MG{\catop}((0,0),(\x,0)) \in \Omega$ in the sense of  \MG{\eqref{eq:linfestimateDelta1}}. Indeed, in this case 
	\MG{\[(\tilde x^k(n),\tilde y^k(n)) =
		(\{f^k_{x,i,n+1}(\iota^{0})\}_{i=1}^N,
		\{f^k_{y,j,n+1}(\iota^{0})\}_{j=1}^m),
		\] for each $n \in \N$ and $k \in \{1,2,\dotsc,\ell\}$}, and it
		immediately follows 
		that  \MG{\eqref{eq:linfestimateDelta1}} is satisfied, by the definition of the $f^k_{x,i,n}$'s and 
		$f^k_{y,j,n}$'s.

		Moreover, it is clear from \eqref{eq:recursivef} and the definition of $(\tilde x^k(n),\tilde y^k(n))$ that $g$ is recursive. Hence, we may now view $\tilde \iota = \tilde \iota(T,\xi)$ as an element in $\tilde \Omega$, which is computable (this is in the standard Turing model where one can access $\tilde \iota$ from an oracle tape), or we may view $\tilde \iota$ as being represented by an algorithm 
		-- as in the Markov model. In both cases, the rest of the argument is identical. Note that the mapping $(T,\xi) \mapsto \tilde \iota(T,\xi)$ is recursive since $g$ is recursive. 		
		
		By assumption, there exist a Turing machine 
		$\Gamma: \tilde \Omega \to \NN_{m,N}^{\F}$ such that 
		$
		\Gamma( \tilde \iota(T,\xi)) \in \mathcal{N}_{10^{-1}}(\Xi(\tilde \iota (T,\xi))),
		$ 	
		where, for a given set $S \subseteq \NN_{m,N}^{\F}$
		\[
		\mathcal{N}_{10^{-1}}(S) = \{ \mathbf{N}: \R^m \to \R^N \: | \: \inf_{\mathbf{N'} \in S} 
		\sup_{y \in \bigcup_{\mathcal{F}}\M_2}  
		\| \mathbf{N}(y) - \mathbf{N'}(y) \|_{\MG{2}} \leq 10^{-1}\}
		\]
		is the $10^{-1}$ neighbourhood of the set $S$.  
		Note that $\Gamma$ is actually defined only on the computable elements of $\tilde \Omega$, but we allow this slight abuse of notation for simplicity. In order to construct $\hat{\Gamma}: \hat{\Omega} \to \{0,1\}$ we proceed as follows. 		 
		Let $\mathcal{S}_{T,\xi}$ denote the string produced by $\Gamma( \tilde \iota(T,\xi))$ that uniquely determines the NN, as explained in Remark \ref{rem:computeNN}. By the same remark, there is a Turing machine $\Gamma^1$ that produces a Turing machine $\Gamma^2$ as follows: $\Gamma^2_{T,\xi} = \Gamma^1(\mathcal{S}_{T,\xi})$ such that $\Gamma^2_{T,\xi}$ computes $\mathbf{N}_{T,\xi} = \Gamma( \tilde \iota(T,\xi))$, that is
		$
		\|\Gamma^2_{T,\xi}(\tilde y, 2^{-n}) - \mathbf{N}_{T,\xi}(y)\|_{\MG{2}} \leq 2^{-n}, 
		$
		for all $n \in \mathbb{N}$, where $\tilde y$ is any $\Delta_1$ approximation to $y \in \mathbb{R}^m$. 
		We can now define
		\begin{align*}
			\hat{\Gamma}(T,\xi) = 
			\begin{cases}
				1 &\quad \text{if} \; \|\Gamma^2_{T,\xi}(\tilde 0, 2^{-10})\|_{\MG{2}} \leq \frac{1}{10} + \frac{1}{100}  ,\\
				0 &\quad \text{otherwise},
			\end{cases}
		\end{align*}		 
		where $\tilde 0 = \{0,0, \hdots \}.$ Note that $\hat \Gamma$ is recursive. Indeed, the mapping $(T,\xi) \mapsto \tilde \iota(T,\xi)$ is recursive (established above), and thus $(T,\xi) \mapsto \mathcal{S}_{T,\xi}$ (by the assumption that $\Gamma$ is a Turing machine) is recursive. Hence, the mapping $(T,\xi, \tilde y) \mapsto \Gamma^2_{T,\xi}(\tilde y, 2^{-10})$ is recursive. This implies finally that the mapping $(T,\xi) \mapsto \hat \Gamma(T,\xi)$ is recursive as \MG{exactly} evaluating \MG{the rational number} $\MG{\|\Gamma^2_{T,\xi}(\tilde 0, 2^{-10})\|^2_{\MG{2}}}$ can be done in finitely many arithmetic operations. 		 
		
		To finalise the argument we only need to show that		 	 
		$\hat{\Gamma}(T,\xi) = \hat{\Xi}(T,\xi)$ for all 
		$(T,\xi) \in \hat{\Omega}$: This follows from the following two arguments.

		\underline{$\hat{\Xi}(T,\xi) = 1 \implies \hat{\Gamma}(T,\xi) = 1$:} Since $\hat{\Xi}(T,\xi) = 1$, $T$ halts on input $\xi$. Thus, there exists an 
		$n' \in \N$ such that $T$ halts on \MG{$\xi$} after $n'$ steps. Hence, by \eqref{eq:def_iota} \MG{and the argument that follows in Case 1} we have that $\tilde \iota(T,\xi)$ is an inexact 
		representation of the element $\iota^1_{n'} \in \Omega$, and thus 
		$\mathbf{N}_{T,\xi} = \Gamma( \tilde \iota(T,\xi)) \in \mathcal{N}_{10^{-1}}(\Xi(\iota_{n'}^1))$, by assumption.
		Note that, as shown in the proof of Proposition \ref{prop:driving}, every optimal network $\mathbf{N}$ for the training set $\iota^1_{n'}$ satisfies $\mathbf{N}(0) = 0$. Hence, since $\mathbf{N}_{T,\xi} \in \mathcal{N}_{10^{-1}}(\Xi(\iota_{n'}^1))$ and $\Gamma^2_{T,\xi}$ computes $\mathbf{N}_{T,\xi}$, we must have $\|\Gamma^2_{T,\xi}(\tilde 0, 2^{-10})\|_{\MG{2}} \leq \frac{1}{10} + \frac{1}{100}$, thus $\hat{\Gamma}(T,\xi) = 1$.  
		
		\underline{$\hat{\Xi}(T,\xi) = 0 \implies \hat{\Gamma}(T,\xi) = 0$:} Since $\hat{\Xi}(T,\xi) = 0$, $T$ does not halt on input 
		$\xi$. Thus, by \eqref{eq:def_iota} \MG{and the argument that follows in Case 2}, we have that $\tilde \iota(T,\xi)$ is an inexact representation of $\iota^0 \in \Omega$, and 
		therefore, by assumption, $\mathbf{N}_{T,\xi} = \Gamma( \tilde \iota(T,\xi)) \in \mathcal{N}_{10^{-1}}(\Xi(\iota^0))$. 
		As shown in the proof of Proposition \ref{prop:driving}, every optimal network $\mathbf{N}$ for the training set $\iota^0$ satisfies $\mathbf{N}(0) = \frac{1}{2}v$, where $\|v\|_{\MG{2}} = 2 \kappa$ with $\kappa = \frac{3}{8}$. Hence, since $\mathbf{N}_{T,\xi} \in \mathcal{N}_{10^{-1}}(\Xi(\iota^0))$ and $\Gamma^2_{T,\xi}$ computes $\mathbf{N}_{T,\xi}$, we must have $\|\Gamma^2_{T,\xi}(\tilde 0, 2^{-10}) - \frac{1}{2}v \|_{\MG{2}} \leq \frac{1}{10} + 2^{-10}$. Thus, $\|\Gamma^2_{T,\xi}(\tilde 0, 2^{-10})\|_{\MG{2}} \geq 3/8 - \frac{1}{10} - 2^{-10} >  \frac{1}{10} + \frac{1}{100}$, and therefore $\hat{\Gamma}(T,\xi) = 0$. 
	\end{proof}

	\subsection{Formal statement and proof of Theorem \ref{thm_intro:second}} \label{proof2}
	
	In the language of the SCI hierarchy, Theorem \ref{thm_intro:second} has the following formal form. 
	
	\newcommand{\orker}[1]{\ker(#1)^{\bot}}		
	\begin{thm} \label{thm:second_formal}
		For any integers $N > m$, and any $\beta > 0$, 
		consider any fixed non-zero linear map 
		$A: \R^N \to \R^m$ such that the spectrum 
		$\mathrm{Sp}(AA^*) \subset [\beta^2,\infty)$. Then, for any rational 
		$\epsilon_1 \in (0, 3/16]$ and any integer \MG{$\ell \: >\; m +2$}, there exists a domain 
		$\Omega \subset T_{\ell}(A)$ (as described in \eqref{omega}) of training sets and a set of corresponding initial domains 
		$\{\M_1(\T) : \T \in \Omega\}$ (as described in \eqref{M1_dep}) such that the following occurs.
		For the mapping
		$\Xi:\Omega \rightrightarrows \NN_{m,N}^{\F}$ (as described in \eqref{eq:Xi}), we have that
		$\Xi(\mathcal{T}) \neq \emptyset$, for each $\T \in \Omega$.
		Moreover, each of the following happen simultaneously:
		\begin{itemize}
			\item[(i)] 
			For the computational problem $\{ \Xi, \Omega, \M', \Lambda\}^{\Delta_1}$ where $\M'$ is defined in \eqref{breakdown_metric}, we have 
			$\strbdepsp \geq \epsilon_1$ for $\mathrm{p} \in [0,1/2)$. In particular, no algorithm, not even randomized, can approximate any optimal neural network 
			$\mathbf{N}_{opt}^{\M_1(\T)} \in \Xi(\mathcal{T})$ \MG{(for all choices of $\T \in \Omega$ and inputs 
			$\tilde \T \in \tilde \Omega$ corresponding to $\T$)}  to accuracy $\epsilon_1$ (with probability greater than $p > 1/2$ in the randomized case).
			
			\item[(ii)]  There exists a $\T_1 \in \Omega$ and a $D>0$ such that $\Xi_{D}(\T_1) \neq \emptyset$ (recall that $\Xi_D$ was defined in \ref{eq:Xi_D}), that is, there exists an optimal neural network $ \mathbf{N}_{opt}^{\M_1(\T_1)}$
			for $\T_1$ such that
			the Lipschitz constant $L(\mathbf{N}_{opt}^{\M_1(\T_1)}) \leq D$. However, there exists a $\tilde \T_1 \in \tilde \Omega$ \MG{corresponding to $\T_1$} such that for any $\delta \in (0, \epsilon_1)$ and any algorithm $\Gamma$ such that \MG{(for any $\T \in \Omega$ and inexact input $\tilde \T \in \tilde \Omega$ corresponding to $\T$)}
			$\mathbf{N}_{\tilde \T,\epsilon} = \Gamma(\tilde \T,\epsilon)$  is a NN that approximates an optimal NN 
			$\mathbf{N}_{opt}^{\M_1(\T)} \in \Xi(\mathcal{T})$ to accuracy 
			$\epsilon \in (\epsilon_1, 2\epsilon_1 - \delta]$ \MG{(in the sense of \eqref{eq:Xi_norm})} we have the following: 
			The vector $0 \in \pi_2(\T_1)$ but, for every $\eta > 0$ and $K >0$, there exists a $y \neq 0$ with $\|y\|_{\MG{2}} \leq \eta$ and
			\begin{align}\label{eq:blow_up}
				\frac{\|\mathbf{N}_{\mathcal{\tilde T}_1,\epsilon}(0) - 
					\mathbf{N}_{\mathcal{\tilde T}_1,\epsilon}(y)
					\|_{\MG{2}}}{\|y\|_{\MG{2}}} > K.
			\end{align}
			\item[(iii)] There exists an algorithm $\Gamma$ such that for all $\T\in \Omega$, \MG{all
			$\tilde \T \in \tilde \Omega$ which correspond to $\T$ and all dyadic}
			$\epsilon > 2 \epsilon_1$ we have that $\mathbf{N}_{\tilde \T,\epsilon} = 
			\Gamma(\tilde \T,\epsilon)$  is a NN that approximates an optimal neural network 
			$\mathbf{N}_{opt}^{\M_1(\T)} \in \Xi(\mathcal{T})$ to accuracy
			$\epsilon$ \MG{(in the sense of \eqref{eq:Xi_norm})} and is such that the Lipschitz constant 
			$L(\mathbf{N}_{\tilde \T,\epsilon})$
			is uniformly bounded by $2/\beta$.
		\end{itemize}
	\end{thm}
	\MG{\begin{rem}\label{rem:whereFail}
		The result in part (i) is valid on the metric space $\M'$, i.e. the failure to approximate $\mathbf{N}_{opt}^{\M_1(\T)}$ can be observed by considering some element $z \in \bigcap_{\mathcal{F}} \M_2$. By contrast, the $\epsilon$ approximations in parts (ii) and (iii) should be interpreted as an approximation to an optimal neural network on $\bigcup_{\mathcal{F}} \M_2$ as per \eqref{metric}. This is intentional -- the negative result in part (i) is strongest when we measure the error on the smallest possible subset of the training data.
	\end{rem}}
	\begin{rem}
		Part (ii) shows that there exists an optimal neural network $\mathbf{N}_{opt}^{\M_1(\T_1)}$ for $\T_1$ that has a bounded Lipschitz constant, but the Lipschitz constant of any computed $\epsilon$-approximation ($\epsilon \in (\epsilon_1, 2\epsilon_1 - \delta]$) blows up. Moreover, this blow up is on the training data in the sense of \eqref{eq:blow_up} \MG{since $0 \in \pi_2(\T_1)$}. 
	\end{rem}

	\begin{proof}[Proof of Theorem \ref{thm:second_formal}]
		Before proving (i)-(iii), we will provide a preliminary setup that will be done in two steps.
		
		{\bf Step I -- Constructing $\Omega$:}
		Our aim is to construct a domain $\Omega \subseteq T_{\ell}(A)$ of the form 
		$\Omega = \{\iota_n^1\}_{ n\in \N} \cup \{\iota^0\}$, 
		where $\{\iota_n^1\}_{ n\in \N}, \{\iota^0\} \subseteq T_{\ell}(A)$ and where 
		$\{\iota_n^1\}_{ n\in \N}$ and $\iota^0$ are training sets that satisfy assumptions (a)-(b) in Proposition \ref{computability_prop}.
		
		Indeed, we use a very similar construction as in the proof of Proposition \ref{prop:driving}, 
		but with slight modifications. Let $\x \in \ker(A) \neq \{0\}$ (as $N > m$) be a fixed element such that 
		$\| \x \|_{\MG{2}} = 4 \epsilon_1$, and let 
		$w \in \ker(A)^{\bot}$ be some fixed unit vector such that 
		$\| \alpha e_i - \frac{\theta}{4^n}A(w)\|_{\MG{2}} \geq \alpha$ for all 
		$i = 1, \dots, \MG{m}$ and for all $n \in \N$, where $e_i$ are the standard unit vectors in 
		$\R^m$ for \MG{each} $i = 1, \dots,m$, where $\alpha$ is some fixed positive rational number with  $\alpha \leq \frac{1}{4}\min\{1/\beta,\beta\}$, and where $\theta$ is a positive rational number with $\theta \leq \alpha\|A\|^{-1}_{op}/8$. We shall use throughout the obvious inequality $\alpha \leq 1/4$ \MG{and that $\theta \leq \alpha/(8\beta) \leq 1/32$.}
		
		We start by proving that 
		there exists a $w \in \ker(A)^{\bot}$ that satisfies the conditions stated above. 
		Indeed, we observe that any
		vector $z \in (-\infty,0)^m$ satisfies $\| \alpha e_i - z\|_{\MG{2}} \geq \alpha$. Thus, it suffices to prove that we can pick an $w \in \ker(A)^{\bot}$ with $\|w\|_{\MG{2}} = 1$ such that 
		$A(w) \in (-\infty,0)^m$. To do this, consider the map $A_0: \ker(A)^{\bot} \to \R^m$ , defined so that $A_0$ is the restriction of $A$ to $\ker(A)^{\bot}$. Since $\mathrm{Sp}(AA^*)$ does not contain $0$, this map
		is bijective with $A^{\dagger} = (A_0)^{-1}$, where $A^{\dagger}$ is the pseudo inverse of $A$. Thus, we may choose an (arbitrary)
		$z_0 \in (-\infty,0)^m$ and define $w = \frac{1}{\|A^{\dagger}(z_0)\|_{\MG{2}}} A^{\dagger}(z_0)$. Then 
		$w \in \ker(A)^{\bot}$ with $\|w\|_{\MG{2}} = 1$, moreover 
		$A(w) = \frac{1}{\|A^{\dagger}(z_0)\|_{\MG{2}}} AA^{\dagger}(z_0) 
		= \frac{1}{\|A^{\dagger}(z_0)\|_{\MG{2}}} z_0 \in  (-\infty,0)^m$, from which it follows that 
		$\| \alpha e_i - \frac{\theta}{4^n}A(w)\|_{\MG{2}} \geq \alpha$ for all 
		$i = 1, \dots, m$ \MG{and for all $n \in \mathbb{N}$}.

		We now construct the relevant training sets and $\Omega$. To do so, we prove that there \MG{exist} infinitely many choices of \MG{tuples} $\TT' \MG{\in T_{\ell-m-2}}(A)$ satisfying each of the following simultaneously:
		\begin{enumerate}[label=(P\arabic*):,ref = (P\arabic*),series=propertiesOfTraining]
			\item For all $(x,y) \in 
			\TT'$ we have $y = Ax$ with $x \in \ker(A)^{\bot}$,\label{Prop:y=Ax}
			\item The \MG{tuples} $\pi_1(\TT')$ and $\pi_2(\TT')$ are of \MG{length $\ell - m - 2$ with each entry distinct}, \label{Prop:EnoughPoints}
			\item $\pi_2(\TT') \cap B_{\alpha}(\alpha e_i) = \emptyset$ for all $i = 1, \dots, m$,  \label{Prop:AwayFromAlphaEi}
			\item $\pi_2(\TT') \cap B_{\alpha}(0) = \emptyset$. \label{Prop:AwayFrom0}
			\item$\pi_1(\TT') \cap  \{ A^{\dagger}(\alpha e_i)\}= \emptyset$ for all 
			$i = 1, \dots, m$, \label{Prop:AwayFromTheBasis}
			\item $\pi_2(\TT') \cap \{ A(\frac{\theta}{4^n} w + \x)\} = \emptyset$ for all $n \in \N$. \label{Prop:AwayFromThetaw}
		\end{enumerate}
		\noindent where $B_{\alpha}(z) = 
		\{ x \in \R^m \: : \: \|x - z \|_{\MG{2}} < \alpha$ \},
		and where we recall \eqref{eq:pi1def} and \eqref{eq:pi2def}.
		
		To see that there exists (infinitely many) \MG{tuples} $\TT'$ satisfying conditions (P1)-(P6) above, first note that there are (infinitely many) possible choices for $z \in \R^m$ with \MG{$\|z \|_{\MG{2}} = \gamma\min\{1,\beta\}$ where $\gamma = \gamma_z \in (2/3,1]$}.
		Let $Z$ be an arbitrary set of $\MG{\ell}-m-2$ such distinct $z$ such that $z \neq A(\theta w/4^n  + v)$ for all $n \in \mathbb{N}$. This is possible since the set $\{A(\theta w/4^n  + v) \, \vert \, n \in \mathbb{N}\}$ is countable, whereas the set $\MG{\{z \, \vert \, \|z \|_{\MG{2}} = \gamma\min\{1,\beta\}, \gamma \in (2/3,1]\}}$ is uncountable; in fact, we see that there are infinitely many choices of such $Z$. \MG{We arbitrarily order the set $Z$ so that $Z:= \{z^1,z^2,\dotsc,z^{\ell-m-2}\}$, and define \[
		\TT' = ((A^{\dagger}(z^1), z^1), (A^{\dagger}(z^2), z^2), \dotsc (A^{\dagger}(z^{\ell-m-2}), z^{\ell-m-2})). \]
		We claim that properties \ref{Prop:y=Ax} - \ref{Prop:AwayFromThetaw} are satisfied by $\TT'$.}

		We recall that 
		$A_0$ as defined above is bijective with 
		$A^{\dagger} = (A_0)^{-1}$. In particular, for $z \in Z$ we must have
		$A^{\dagger}(z) \in \ker(A)^{\bot}$ with 
		$\|A^{\dagger}z\|_{\MG{2}} \leq \frac{1}{\beta}\|z\|_{\MG{2}} \leq 1$, where the bound on the norm \MG{of $A^{\dagger}$}
		follows from the assumption that $\text{Sp}(AA^*)\subset [\beta^2,\infty)$ and that 
		$A^{\dagger} = A^*(AA^*)^{-1}$. 
		We set $x = A^{\dagger}z$ and observe that, 
		with $y = Ax \MG{= z}$, we get that $\|y\|_{\MG{2}} = \|AA^{\dagger}z \|_{\MG{2}} = \|z\|_{\MG{2}} 
		= \MG{\gamma_z}\min\{1,\beta\} \MG{\leq 1}$. These inequalities imply \MG{$\TT' \in T_{\ell - m - 2}(A)$}, and that \ref{Prop:y=Ax} holds.
		
		To see that \ref{Prop:EnoughPoints} holds, note that $\pi_1(\TT') = A^{\dagger}(Z)$. But $A^{\dagger}$ is a bijection from $\R^m \to \ker(A)^{\bot}$ and thus $|\pi_1(\TT')| = |Z| = \ell-m-2$. Similarly, $\pi_2(\TT') = A (\pi_1(\TT')) = A_0(\pi_1(\TT'))$ since $\pi_1(\TT') \subset \ker(A)^{\bot}$. But as $A_0$ is a bijection we must therefore have \MG{$|\pi_2(\TT')| = |\pi_1(\TT')| = \ell-m-2$}, thus proving that \ref{Prop:EnoughPoints} holds. 
		
		Next, we observe that for $y \in \pi_2(\TT')$, we must have $2\alpha \leq \frac{1}{2}\min\{1,\beta\} < \|y\|_{\MG{2}}$.
		Therefore, we may conclude that $y \notin \bigcup_{i = 1}^m B_{\alpha}(\alpha e_i)$, 
		since all the elements in 
		$\bigcup_{i = 1}^m B_{\alpha}(\alpha e_i)$ have norm less than $2\alpha$, thus proving \ref{Prop:AwayFromAlphaEi}. Also, $y \notin B_{\alpha}(0)$, since all the elements in $B_{\alpha}(0)$ have norm less than $\alpha$, which proves \ref{Prop:AwayFrom0}. For $(x,y) \in \TT'$, we must have
		$x \neq A^{\dagger}(\alpha e_i)$ for any $i = 1, \dots m$: to see this, we argue by contradiction and see that if this does not hold then
		$y =  AA^{\dagger}(\alpha e_i) = \alpha e_i$, 
		which implies that $\|y \|_{\MG{2}} = \alpha$, contradicting the fact that 
		$2\alpha < \|y\|_{\MG{2}}$, from which we conclude \ref{Prop:AwayFromTheBasis}. Finally, \ref{Prop:AwayFromThetaw} is immediate from the assumption that $z \in Z$ implies \MG{$z \notin \{A(\theta w/4^n  + v) \, \vert \, n \in \mathbb{N}\}$}. Next, let $\TT = (\TT', (\alpha A^{\dagger}(e_1),\alpha e_1),\dots,
		(\alpha A^{\dagger}(e_m),  \alpha e_m))$,
		where $e_1, \dots, e_m$ denotes the standard basis in $\R^m$,
		and where $\TT'   = ((x^k,y^k))_{k =1}^{\ell - m - 2}$
		satisfies points \ref{Prop:y=Ax}-\ref{Prop:AwayFromThetaw}.
		
		We next define $\iota^0$ and the sequence $\{\iota_n^1\}_{n\in \N}$ as follows: for each $n \in \N$\MG{,} let
		\begin{align*}
			\iota_n^1 = \TT \MG{\catop}((0,0),(\frac{\theta}{4^{n } }w + \x, A(\frac{\theta}{4^{n}}w + \x))) 
			\quad
			\text{and}
			\quad
			\iota^0 = \TT \MG{\catop}((0,0),(\x,0)).
		\end{align*}
		At last, we define the corresponding sets $\M_1(\iota_n^1)$ and $\M_1(\iota^0)$, 
		just as in the proof of Proposition \ref{prop:driving}, as follows: we let 
		$\mathcal{M}_1(\iota_n^1) = \pi_1(\iota_n^1)$ and 
		$\mathcal{M}_1(\iota^0) = \pi_1(\iota^0)$. We define $\Omega = \{\iota_n^1\}_{ n\in \N} \cup \{\iota^0\}$.
		
		%	In other words, for $\iota_n^1$ we have the following: the first $\ell-2$ elements 
		%		are equal to the the elements in 
		%		$\T_b$, $(x^{\ell-1},y^{\ell-1}) = (0,0)$, and $(x^{\ell},y^{\ell}) 
		%		= (\x \!+ \! \frac{\theta}{4^n}w,A(\x \!+ \! \frac{\theta}{4^n }w))$. 
		%		Similarly, for $\iota_n^2$ we have the following,
		%		the first $\ell-2$ elements are equal to the the elements in 
		%		$\T_b$, $(x^{\ell-1},y^{\ell-1}) = (0,0)$, and $(x^{\ell},y^{\ell})=(0,\x)$.
		%	Using \ref{Prop:AwayFromThetaw}, we see that 
		%	$\iota_n^1, \iota_n^2 \subseteq T_{\ell}(A)$ for all $n \in \N$.
		\newcommand{\SecOmegaAndXiNonEmpty}{Step II -- Showing that $\Omega$ is a subset of \MG{$T_\ell(A)$} and that $\Xi(\T) \neq \emptyset$ for all $\T \in \Omega$:}
		\textbf{\SecOmegaAndXiNonEmpty}
		We start by showing that \MG{$\pi_1 (\TT \catop ((0,0)))$ and $\pi_2(\TT \catop ((0,0)))$} each have \MG{$\ell-1$} distinct points. To see this, since $A$ is of full rank we know that $A^{\dagger}(e_i) \neq A^{\dagger}(e_j)$ for any 
		$i \neq j$, since the columns of $A^{\dagger}$ are linearly independent. Moreover (using \ref{Prop:AwayFromAlphaEi} and \ref{Prop:AwayFromTheBasis}) we have chosen 
		$\TT'$ such that $\alpha A^{\dagger}(e_i) \notin \pi_1(\TT')$ and 
		$\alpha e_i \notin \pi_2(\TT')$ for any $i = 1, \dots, m$, and such that (using \ref{Prop:EnoughPoints} and \ref{Prop:AwayFrom0}) the \MG{tuples}
		$\pi_1(\TT')$ and $\pi_2(\TT')$ each have $\ell - m -2$ non-zero distinct points. 
		Using this we conclude that the \MG{tuples $\pi_1 (\TT \catop ((0,0)))$ and $\pi_2(\TT \catop ((0,0)))$} must each have $\ell-1$ distinct points. 
		
		Next, we observe that
		$A(\frac{\theta}{4^{n}}w + \x) \notin \pi_2(\TT)$, due to \ref{Prop:AwayFromThetaw} above and \MG{the fact that the choice of $w$ implies that}
		$\| \alpha e_i - A(\frac{\theta}{4^n} w  + \x)\|_{\MG{2}}=\| \alpha e_i - \frac{\theta}{4^n}A(w)\|_{\MG{2}} \geq \alpha$ for all 
		$i = 1, \dots, \MG{m}$ and for all $n \in \N$. Moreover,
		$\frac{\theta}{4^{n } }w + \x \notin \pi_1(\TT')$, due to \ref{Prop:y=Ax} (specifically, the fact that $\MG{0 \neq }\x \in \ker(A)$ implies that $\frac{\theta}{4^{n } }w + \x  \notin \orker{A}$). Thus $\frac{\theta}{4^{n } }w + \x \notin \pi_1(\TT)$
		since if 
		$\frac{\theta}{4^{n } }w + \x = \alpha A^{\dagger}(e_i)$, 
		for some $n \in \N$ and $i \in  \{1, \dots, m\}$,
		then $\alpha AA^{\dagger}(e_i) = A(\frac{\theta}{4^{n } }w + \x)$, which implies that $\alpha e_i = \frac{\theta}{4^n}A(w)$ which contradicts $\| \alpha e_i - \frac{\theta}{4^n}A(w)\|_{\MG{2}} \geq \alpha$. 
		Hence, we may conclude that the \MG{tuples} $\pi_1(\iota_n^1)$ and $\pi_2(\iota_n^1)$ each have $\ell$ distinct points. 
		
		\MG{To conclude that $\iota^1_n \in T_{\ell}(A)$, we must also show that for every $(x,y) \in \iota^1_n$ we have $\|x\|_{\MG{2}}, \|y\|_{\MG{2}} \leq 1$. This has already been proven for $(x,y)$ in $\T'_b$: if instead $(x,y) \in \TT$ but $(x,y) \notin \T'_b$ then $(x,y) = (\alpha A^{\dag}e_i,\alpha e_i)$ for some $i$. But then $\|x\|_{\MG{2}} \leq \alpha \|A^{\dag}\|_{op} \leq \alpha/\beta\leq 1/4$ and $\|y\|_2 = \alpha \leq 1/4$. If $(x,y) = (0,0)$ the bounds $\|x\|_{\MG{2}}, \|y\|_{\MG{2}} \leq 1$ are trivial. Finally, if $(x,y) = (\theta w /4^n  + v, A(\theta w/4^n + v)) = (\theta w /4^n +v, \theta Aw/4^n)$ then \[
			\|x\|_{\MG{2}} \leq \frac{\theta \|w\|_{\MG{2}}}{4^n} + \|v\|_{\MG{2}} \leq \frac{1}{32 \cdot 4^n} + 4\epsilon_1 \leq \frac{1}{32 \cdot 4} +\frac{3}{4} \leq 1
			\]
			by the bounds $\epsilon_1 \leq 3/16$, $\theta \leq 1/32$. Also $\|y\|_{\MG{2}} \leq \theta \|A\|_{op}/4^n \leq \alpha/(8\cdot 4^n) \leq 1$. Thus $\iota^1_n \in T_{\ell}(A)$.}
		
		The argument that \MG{$\iota^0 \in T_\ell(A)$} is slightly different: we know that \MG{$\TT \catop ((0,0))$} has exactly $\MG{\ell}-1$ distinct points (otherwise both $\pi_1(\MG{\TT \catop ((0,0))})$ and $\pi_2(\MG{\TT \catop ((0,0))})$ \MG{would} have fewer than $\MG{\ell}-1$ points). We claim that $(v,0) \notin \MG{\TT \catop ((0,0))}$: if $v \in \pi_1(\TT')$ then by \ref{Prop:y=Ax} \MG{we would get} $Av =0  \in \pi_2(\TT')$, which contradicts \ref{Prop:AwayFrom0}. If $v = \alpha A^{\dagger}e_i$ for some $i$ then $0 = Av = \alpha AA^{\dagger}e_i = \alpha e_i$ which is a contradiction, so $(v,0) \notin \TT$. Finally, it is obvious that $(v,0) \neq (0,0)$. 
		
		\MG{Once again, all that remains is to prove that for $(x,y) \in \iota^0$, $\|x\|_{\MG{2}}, \|y\|_{\MG{2}} \leq 1$. This is mostly an identical argument to the one presented for $\iota^1_n$: here, the only difference is that we must argue that if $(x,y) = (v,0)$ then $\|x\|_{\MG{2}}, \|y\|_{\MG{2}} \leq 1$. But $\|v\|_{\MG{2}} = 4\epsilon_1 \leq 3/4 \leq 1$ and $\|y\|_{\MG{2}} = 0 \leq 1$.} Thus, $\iota^0 \in \MG{T_\ell(A)}$ and we conclude that $\Omega \subset \MG{T_\ell(A)}$.
		
		We now prove the fact that $\Xi(\iota_n^1) \neq \emptyset$ for all $n \in \N$ \MG{and that $\Xi(\iota^0) \neq \emptyset$} in a similar fashion as in the proof of Proposition \ref{prop:driving}. Starting with \MG{$\Xi(\iota_n^1)$}, since $\pi_1(\iota^1_n)$ and $\pi_2(\iota^1_n)$ both have \MG{$\ell$} points we can use Lemma \ref{lem:radial_network} to conclude that, for each $n \in \mathbb{N},$ there is a neural network $\mathbf{N}_n: \mathbb{R}^m \rightarrow \mathbb{R}^N$  such that for each pair 
		$(\xi_n, \eta_n) \in \iota_n^1$ we have $\mathbf{N}_n(\eta_n) = \xi_n$. It is clear from Definition 
		\ref{optmap:trainingset} that this is an optimal map for the inverse problem 
		$(A,\mathcal{M}_1(\iota_n^1))$. 
		
		By the same argument we can find a neural network 
		$\mathbf{N}: \mathbb{R}^m \rightarrow \mathbb{R}^N$ so that $\mathbf{N}(0) = v/2$ and such that for any $(\xi, \eta) \in \TT$ we have 
		$\mathbf{N}(\eta) = \xi$.  
		We claim that $\mathbf{N}$ is optimal for $(A,\mathcal{M}_1(\iota^0))$. Indeed, note that
		\MG{\begin{equation}\label{0_optimal}
				\inf_{\varphi \colon\! \mathcal{M}_2(\iota^0)  \rightrightarrows \R^N} 
				\sup_{x \in \mathcal{M}_1(\iota^0)} d_1^{H}(\varphi(Ax), x) \geq  
				\inf_{\varphi \colon\! \mathcal{M}_2(\iota^0)  \rightrightarrows \R^N} 
				\max_{x =0,x=v} d_1^{H}(\varphi(Ax), x) = \|v\|_{\MG{2}}/2. 
		\end{equation}}
		However, by the definition of $\mathbf{N}$ \MG{(and since $\mathbf{N}$ is single-valued)} we have that
		\begin{equation}
			\sup_{x \in \mathcal{M}_1(\iota^0)} d_1^{\MG{H}}(\mathbf{N}(Ax), x) 
			= \max_{x=0,x=v} \MG{\|\mathbf{N}(Ax) - x\|_{\MG{2}}}  = \|v\|_{\MG{2}}/2,
		\end{equation}
		proving our claim. This means that  \MG{$\Xi(\iota_n^1) \neq \emptyset$ for all $n \in \N$ and that $\Xi(\iota^0) \neq \emptyset$}. 
		Thus the computational problem $\{ \Xi, \Omega, \M', \Lambda\}$ is now well defined.

		{\bf Proof of (i):} By the same argument as in the proof of 
		Proposition \ref{prop:driving}, we can prove that 
		$\{\iota_n^1\}_{ n\in \N}$ and $\iota^0$ satisfy assumptions 
		(a)-(b) in Proposition \ref{computability_prop} with $\kappa = 2 \epsilon_1$.
		Part (i) follows from this.
		
		{\bf Proof of (ii):}  We start by defining the relevant objects $\T_1$, $D$, $\mathbf{N}_{opt}^{\M_1(\T)}$.  We set  $\iota^0 \coloneqq\MG{ \TT \catop((0,0),(\x,0))}$ and $\T_1 \coloneqq \iota^0$. The existence of $\mathbf{N}_{opt}^{\M_1(\T)} \in \Xi(T)$ was discussed earlier in the proof when $\Omega$ was constructed; in particular, this proof relied on Lemma \ref{lem:radial_network}. From that Lemma we see that $\mathbf{N}_{opt}^{\M_1(\T)}$ is Lipschitz continuous; choose $D$ to be any valid Lipschitz constant for this network. We prove the remainder of (ii) in the following three parts:

		\emph{Part (1)}: We now define $\tilde{\T}_1$ \MG{which we will claim corresponds to $\T_1$.} To do so, we begin by defining $\Delta_1$ information $\hat{\Lambda}$. \MG{We write} $\Lambda = \{f_k \, \vert \, k \in I\}$ where $I$ is a finite index set and let (for $m,n \in \N$ and $k \in I$) $d_{k}^{m,n}$ \MG{be} a dyadic number such that 
		$|d_{k}^{m,n} - f_k(\iota_n^1)| \leq 2^{-m}$, and let $c_k^m$ be a dyadic number such that $|c_k^m - f_k(\iota^0)| \leq 2^{-m}/\sqrt{2}$ for $k \in I$.
		We then define $\hat \Lambda$ to be the set of functions $f_{k,m}$, for $k \in I$ and $m \in \N$ where each $f_{k,m}$ is given by
		\begin{align}\label{fmk}
			f_{k,m}(\iota_n^1) =
			\begin{cases}
				d_{k}^{m,n} &\quad \text{if} \; 1 \leq n \leq m \\
				c_k^m  &\quad \text{if} \; n > m,
			\end{cases}
			\quad \text{and} \quad
			f_{k,m}(\iota^0) = c_k^m.
		\end{align}
		Then the set $\hat{\Lambda} = \{ f_{k,m} \: : m \in \N, \: f_k \in \Lambda \}$ provides $\Delta_1$-information for $\Omega = \{\iota_n^1\}_{n \in \N} \cup \MG{\{\iota^0\}}$ in the sense of 
		\eqref{inexact_input}. Indeed, for $\iota_n^1$ with $1 \leq n \leq m$ this is immediate by the definition of $d_{k}^{m,n}$. For $n > m$ we first note (as in the proof of Proposition \ref{prop:driving}) \MG{that $\|\frac{\theta}{4^n} w + v - v \|_{\infty} \leq \theta \|w\|_{\MG{2}}/4^n \leq 1/4^n$ and that $\|A(\theta w/4^n + v) - 0\|_{\infty} \leq \theta \|Aw\|_{\MG{2}}/4^n \leq \theta\|A\|_{op}/4^n \leq 1/4^n$. Hence $|f_k(\iota^1_n) - f_k(\iota^0)| \leq 4^{-n}$ for each $k \in I$ and thus}
		$|f_{k,m}(\iota_n^1) - f_k(\iota_n^1)| \leq \MG{|c_k^m - f_k(\iota^0)| + |f_k(\iota^0) - f_k(\iota_n^1)| }
		\leq 2^{-m}/\sqrt{2} + 4^{-n} \leq 2^{-m}/\sqrt{2} + 2^{-m}/4 < 2^{-m}$.
		
		Similarly, for $\iota^0$ we get that
		$|f_{k,m}(\iota^0) - f_k(\iota^0)| = |c_k^m - \MG{f_k(\iota^0)}| \leq 2^{-m}/\sqrt{2} < 2^{-m}$.
		Thus, $\hat \Lambda$ yields $\tilde \iota^0,\tilde \iota_n^1 \in \tilde \Omega$ (recall Definition \ref{def:delta1} and \eqref{def:inexact_rep}) that correspond to $\iota^0,\iota_n^1 \in \Omega$ respectively. We choose $\tilde{\T}_1 = \tilde \iota^0$.

		\emph{Part (2)}: \MG{Choose an arbitrary $\delta \in (0, \epsilon_1)$ and fix $\epsilon \in (\epsilon_1, 2\epsilon_1 - \delta]$.} We prove that for any algorithm $\Gamma$, such that $\Gamma(\tilde \iota,\MG{\epsilon}) \neq \text{NH}$ 
		for all \MG{$\tilde \iota \in \tilde \Omega$}, there is an $N_0 \in \N$ such that for all 
		$n > N_0$ we have that $\Gamma(\tilde \iota_n^1,\epsilon) = \Gamma(\tilde \iota^0,\epsilon)$ \MG{where $\tilde \iota_n^1$ and $\tilde \iota^0$ are defined via the $\Delta_1$ information provided by $\hat \Lambda$}. Let 
		\begin{align*}
			N_0 = D_{\Gamma}(\tilde \iota^0) \coloneqq \sup \{ m \in \N \: : \: \exists f_k \in \Lambda 
			\: \text{with} \: f_{k,m} \in \hat{\Lambda}_{\MG{\Gamma(\cdot,\epsilon)}}(\tilde \iota^0)\}.
		\end{align*}
		\MG{Note that $N_0$ depends on relevant fixed parameters like $\epsilon$.} By the assumption that $\Gamma(\tilde \iota^0,\epsilon) \neq \text{NH}$ we get that 
		$N_0$ must be finite. Recall from Definition \ref{def:delta1}, that we define $\tilde f_{k,m} \in \tilde \Lambda$
		to act on the elements $\tilde \iota \in \tilde \Omega$ by 
		$\tilde f_{k,m}(\tilde \iota) = f_{k,m}(\iota)$ where $\iota \in \Omega$ is such that $\tilde \iota$ corresponds to $\iota$.
		Then, by \eqref{fmk}, we have that $\tilde f(\tilde \iota^0) = \tilde f(\tilde \iota_n^1)$ for all 
		$\tilde f \in \hat{\Lambda}_{\Gamma(\cdot,\epsilon)}(\tilde \iota^0)$ whenever $n > N_0$. Thus, by part (iii) of 
		Definition \ref{def:General_alg} it follows that 
		$\hat{\Lambda}_{{\MG{\Gamma(\cdot,\epsilon)}}}(\tilde \iota_n^1) = \hat{\Lambda}_{{\MG{\Gamma(\cdot,\epsilon)}}}(\tilde \iota^0)$ whenever 
		$n > N_0$. Consequently, by parts (i) and (ii) of Definition \ref{def:General_alg} it follows that 
		$\Gamma(\tilde \iota_n^1,\epsilon) = \Gamma(\tilde \iota^0,\epsilon)$ whenever $n > N_0$.

		\emph{Part (3)}: Finally, recall that we chose $\mathcal{\tilde T}_1 = \tilde \iota_0$. \MG{Let $K>0$ and $\eta > 0$ and 
			suppose that there is an algorithm $\Gamma$ such that 
			$\mathbf{N}_{\tilde \T,\epsilon} = \Gamma(\tilde \T,\epsilon)$  is a NN that approximates an optimal neural network $\mathbf{N}_{opt}^{\M_1(\tilde \T)} \in 
			\Xi(\tilde \T)$ to accuracy $\epsilon \in (\epsilon_1, 2\epsilon_1 - \delta]$ for all 
			$\tilde \T \in \tilde \Omega$.} We will now show \eqref{eq:blow_up}. Let 
		$n_1$ be such that $n_1 > N_0$, where $N_0$ is chosen according to part (\MG{2}), and such that 
		$\left\|A\left(v+\frac{\theta w}{4^{n_1}}\right)\right\|_{\MG{2}} = \|\frac{\theta}{4^{n_1}}A(w)\|_{\MG{2}} < \eta$ (recall $w$ is a unit vector in $\orker{A}$) and $4^{n_1}(2\delta - 4^{-n_1}\theta) > K$ (which is always possible since $\delta$ is strictly positive). Then, by part (2), we must have that $\mathbf{N}_{\mathcal{\tilde T}_1,\epsilon} = \Gamma(\mathcal{\tilde T}_1,\epsilon) = \Gamma(\tilde \iota^0,\epsilon) = \Gamma(\tilde \iota_{n_1}^1,\epsilon)$. Hence, since $\Gamma$ always provides $\epsilon$ approximations to optimal NNs, it follows that $\mathbf{N}_{\mathcal{\tilde T}_1,\epsilon}$ is an $\epsilon$-approximation to some $\mathbf{N}_{opt}^{\iota_{n_1}^1} \in \Xi(\iota_{n_1}^1)$ in the sense of \eqref{eq:Xi_norm}. Thus (using the observation that we have made several times that the optimal neural network for $\iota^1_{n_1}$ must exactly interpolate the training points in $\iota^1_{n_1}$), we have that $\|\mathbf{N}_{\mathcal{\tilde T}_1,\epsilon}(0)\|_{\MG{2}} \leq \epsilon$, as well as 
		$\|\mathbf{N}_{\mathcal{\tilde T}_1,\epsilon}(A(\frac{\theta}{4^{n_1}}w +\x))\|_{\MG{2}} \geq 
		\|\frac{\theta}{4^{n_1}}w +\x\|_{\MG{2}} - \epsilon \geq \|\x\|_{\MG{2}} - 4^{-n_1}\theta - \epsilon$. We can therefore derive that
		\begin{align*}
			\frac{\|
				\mathbf{N}_{\mathcal{\tilde T}_1,\epsilon}(A(\frac{\theta}{4^{n_1}}w +\x)) - 
				\mathbf{N}_{\mathcal{\tilde T}_1,\epsilon}(0)\|_{\MG{2}}}{\|\frac{\theta}{4^{n_1}}A(w)\|_{\MG{2}}} &\geq \frac{\|\x\|_{\MG{2}} - 2\epsilon - 4^{-n_1}\theta}{4^{-n_1}} \\
			&\geq 
			\frac{4\epsilon_1 - 4\epsilon_1 + 2 \delta - 4^{-n_1}\theta}{4^{-n_1}} = 4^{n_1}(2\delta - 4^{-n_1}\theta) > K.
		\end{align*}
		Clearly, $y = A(\frac{\theta}{4^{n_1}}w +\x) = A(\frac{\theta}{4^{n_1}}w)$ has $\|y\|_{\MG{2}} < \eta$, and thus we have established \eqref{eq:blow_up}.

		{\bf Proof of (iii):} In order to prove part (iii) we will show that there exists an algorithm 
		$\Gamma: \tilde{\Omega} \times \R_+ \to \NN_{m,N}^{\F}$ 
		that works on inexact input, such that for each $\tilde{\T} \in \tilde{\Omega}$, $\Gamma$ constructs a neural network 
		$\Gamma(\tilde{\T}, \epsilon)$ such that
		\begin{align} \label{claim}
			\sup_{x \in \bigcup_\mathcal{F} \mathcal{M}_1, y = Ax} \| \Gamma(\tilde{\T},\epsilon)(y) - x \|_{\MG{2}} 
			\leq \epsilon, 
		\end{align}
		(where analogously to \eqref{eq:Cups}, we define  $\bigcup_{\mathcal{F}} \mathcal{M}_1 := \bigcup_{(A,\mathcal{M}_1) \in \mathcal{F}} \mathcal{M}_1$) for any \MG{dyadic $\epsilon$ with} $\epsilon > 2\epsilon_1$, and such that $\Gamma(\tilde{\T} , \epsilon)$ is Lipschitz \MG{with constant at most $2/\beta$}.
		
		Our strategy is to construct an algorithm that produces a neural network which, on any input, applies the pseudo inverse $A^{\dagger}$ of $A$ and adds the bias term $\frac{1}{2}\x$. In order to achieve this we first run the following algorithm to approximate the pseudo inverse $A^{\dagger}$ of $A$. 
		For $(\tilde x,\tilde y) \in \tilde \T$ we use the same convention as in \eqref{xy_aprox} and denote by 
		$\tilde x(n)$ and $\tilde y(n)$ the $n$th entries in the sequences $\tilde x$ and $\tilde y$ respectively so that
		$\|\tilde x(n) - x \|_{\MG{2}} \leq 2^{-n}$ and $\|\tilde y(n) - y\|_{\MG{2}} \leq 2^{-n}$.\\
		{\it Subroutine} {\bf Approximate $A^{\dagger}$} 
		
		\indent \underline{Inputs}: $k \in \N$ and a training set $\tilde \T = ((\tilde x^1,\tilde y^1),(\tilde x^2,\tilde y^2),\dotsc,\MG{(\tilde x^\ell,\tilde y^\ell)})$ corresponding to \\\indent $\T  = ((x^1,y^1),(x^2,y^2),\dotsc,\MG{(x^\ell,y^\ell)}\MG{)}$. \\
		\indent \underline{Oracles}: Oracles for all $(\tilde{x}, \tilde{y}) \in \tilde{\T}$,\\
		\indent \underline{Output}: An approximation $A^{\dagger}_k$ of $A^{\dagger}$ such that 
		$\| A^{\dagger}_k - A^{\dagger}\|_{op} \leq 2^{-k}$.
		\begin{enumerate}[leftmargin=10mm, label= \arabic*.,ref = \arabic*]
			\item Let $n \in \N$ be a fixed constant such that $2^{-n} \leq \frac{\alpha}{4}$.
			\item  For $i = 1,\dots,m:$\label{inst:ADaggerLoop}
			\begin{enumerate}[leftmargin=8mm, label = \alph*,ref = \ref{inst:ADaggerLoop}\alph*.]
				\item Search through $j \in \{1,2,\dotsc,\MG{\ell}\}$  until we find $j$ such that $(\tilde{x}^j,\tilde{y}^j) \in \tilde{\T}$ with $\| \alpha e_i - \tilde y^{j}(n) \|_{\MG{2}} 
				\leq 2^{-n+1}$ \MG{(note that this can be checked because $\alpha e_i - \tilde y^{j}(n)$ is a rational vector)}. With this $j$: \label{inst:ADaggerSearch}
				\begin{enumerate}[label = \roman*.]
					\item Set $r(k) := \lceil \log_2(m)\rceil + n + k$
					\item Set 
					$A^{\dagger}_{k,i} = \frac{1}{\alpha} \tilde x^{j}(r(k))$, 
				\end{enumerate}
			\end{enumerate} 
			\item Set $A^{\dagger}_k = (A^{\dagger}_{k,1}, \dots, A^{\dagger}_{k,m})$ and return $A^{\dagger}_k$ (note that the notation $A^{\dagger}_{k,i}$ should be interpreted as the $i$th column of the $k$th approximation to $A^{\dagger}$, rather than the entry at row $k$, column $i$).
		\end{enumerate}
		To see that this algorithm fulfils its desired functions, we must show that the algorithm terminates (and in particular, the search in instruction \ref{inst:ADaggerSearch} terminates) and that when the algorithm terminates it terminates with a correct value. We achieve this by making (and proving) the following claim:
		
		\textbf{Claim:} Suppose that $(\tilde x, \tilde y) \in \tilde \T$ corresponds to $(x,y) \in \T$ \MG{with $\T \in \Omega$}. Then for each $i=1,2,\dotsc,m$, we claim that $(x,y) = (\alpha A^{\dagger}(e_i),\alpha e_i)$ if and only if $\|\alpha e_i - \tilde y(n) \|_{\MG{2}} \leq 2^{-n+1}.$ 
		
		To see this, note that if $\|\alpha e_i - \tilde y(n) \|_{\MG{2}} \leq 2^{-n+1}$ then
		\begin{align*}	
			\|\alpha e_i - y\|_{\MG{2}} \leq \|\alpha e_i - \tilde{y}(n)\|_{\MG{2}} 
			+ \|\tilde{y}(n) - y\|_{\MG{2}} 
			\leq 2^{-n+1}+ 2^{-n} \leq \frac{3}{4}\alpha < \alpha\MG{.}
		\end{align*}
		In particular, $(x,y) \notin \TT'$ (since any $(x,y) \in \TT'$ cannot have $\|y-\alpha e_i\|_{\MG{2}} < \alpha$ by Property  \ref{Prop:AwayFromAlphaEi}) and $(x,y) \notin \{(0,0),(\frac{\theta}{4^{\MG{t}} }w + \x, A(\frac{\theta}{4^{\MG{t}}}w + \x)), (\x,0)\} $ \MG{for any $t \in \mathbb{N}$} (since $w$ is chosen so that, \MG{for every $t \in \mathbb{N}$}, $\|\alpha \MG{e}_i -A(\frac{\theta}{4^{t}}w + \x)\|_{\MG{2}} \geq \alpha$). 
		
		Therefore $(x,y) \in \{(\alpha A^{\dagger}(e_1),\alpha e_1),\dots,
		(\alpha A^{\dagger}(e_m),  \alpha e_m)\}$
		and thus $(x,y) = (\alpha A^{\dagger}(e_i),\alpha e_i)$ (since $\|\alpha \MG{e_i} - \alpha e_j\|_{\MG{2}}\geq \alpha$ whenever $i \neq j$) as claimed.
		The converse is straightforward: if $(x,y) = (\alpha A^{\dagger}(e_i),\alpha e_i)$, then
		$
		\|\alpha e_i - \tilde{y}(n)\|_{\MG{2}} = \|y - \tilde{y}(n)\|_{\MG{2}} \leq 
		2^{-n} \leq 2^{-n+1}.
		$
		This ends the proof of the claim.
		
		Note that the claim immediately implies that subroutine Approximate $A^{\dagger}$ terminates: since $(\alpha A^{\dagger}(e_i),\alpha e_i) \in \T$ \MG{for any $\T \in \Omega$}, we have $(\alpha A^{\dagger}(e_i),\alpha e_i) = (x^j,y^j) \in \T$ for some $j \in \{1,2,\dotsc,\MG{\ell}\}$ and thus by the claim $\|\alpha e_i - \tilde y^{j}(n)\|_{\MG{2}} \leq 2^{-n+1}$ for this $j$, so the search in instruction \ref{inst:ADaggerSearch} terminates. Moreover, the converse result proven in the claim also implies that when the search is complete for a given $i$, the value $j$ found in instruction \ref{inst:ADaggerSearch} must be such that $(\tilde{x}^\MG{j}, \tilde{y}^\MG{j})$ corresponds to $(\alpha A^{\dagger}(e_i),\alpha e_i)$. In particular, $\| \tilde x^{j}(r(k)) - A^{\dagger}(\alpha e_i) \|_{\MG{2}} 
		\leq 2^{-r(k)}$. Thus for $i=1,2,\dotsc,m$ we have 
		\MG{$\| A^{\dagger}_{k,i}  - A^{\dagger}(e_i) \|_{\MG{2}} = \| \tilde x^{j}(r(k)) - A^{\dagger}(\alpha e_i) \|_{\MG{2}}/\alpha
		\leq 2^{-r(k)}/\alpha
		$} and hence
		\begin{align*}
			\| A^{\dagger}_k - A^{\dagger} \|_{op} 
			&\leq \sum_{i = 1}^m \| A^{\dagger}_{k,i} - A^{\dagger}(e_i) \|_{\MG{2}} 
			\leq \sum_{i = 1}^m \frac{1}{\alpha} \cdot 2^{-r(k)}  
			\leq m \cdot \frac{1}{\alpha} \cdot \frac{1}{m} \cdot 2^{-n} \cdot 2^{-k} \\
			&\leq m \cdot \frac{1}{\alpha} \cdot \frac{1}{m} \cdot \frac{\alpha}{4} \cdot 2^{-k} 
			\leq \frac{1}{4} \cdot 2^{-k} < 2^{-k}.
		\end{align*}
		Next, we run the following algorithm to approximate the bias term $\x$.\\
		{\it Subroutine} {\bf Approximate $v$} 
		
		\indent \underline{Inputs}: A dimension $m \in \N$, and a $k \in \N$. \\
		\indent \underline{Oracles}: Oracles for all $(\tilde{x}, \tilde{y}) \in \tilde{\T}$,\\
		\indent \underline{Output}: An approximation $v^k$ of $v$ such that 
		$\| v^k - v\|_{\MG{2}} \leq 2^{-k}$.
		\begin{enumerate}[leftmargin=10mm, label= \arabic*.,ref = \arabic*]
			\item \label{inst:kPrimeDef}Let $k' \in \N$ be such that $2^{-k'} \leq \frac{\alpha}{2} \cdot 2^{-k}$ .
			\item Let $r\in \mathbb{N}$ be such that $2^ {-r} \leq \epsilon_1$ and $2^{-r} \leq \alpha/4$.
			\item \label{inst:vSearch}Search through $j \in \{1,2,\dotsc,\MG{\ell}\}$  until we find $j$ such that $(\tilde{x}^j,\tilde{y}^j) \in \tilde{\T}$ with $\|\tilde x^{j}(r)\|^2_{\MG{2}} \geq 9\epsilon_1^2$ and $\|\tilde y^{j}(r)\|_{\MG{2}}^2 \leq  (2^{-r} + \alpha/8)^2$ \MG{(note that these conditions can be checked computationally as $\alpha,\epsilon_1 \in \mathbb{Q}$ and $\tilde x^{j}(r), \tilde{y}^j(r)$ are rational vectors)}. With this $j$ set $v^k = \tilde x^{j}(k') - A^{\dagger}_{\MG{k'}} \tilde y^{j}(k')$ (where $A^{\dagger}_{\MG{k'}}$ is produced by the subroutine `Approximate $A^{\dagger}$').
		\end{enumerate}

		To see this algorithm achieves its desired purpose, we start by making a similar claim as in the proof of the correctness of `Approximate $A^{\dagger}$':
		
		\textbf{Claim:} Suppose that $(\tilde x, \tilde y) \in \tilde \T$ corresponds to $(x,y) \in \T$. Then we claim that $(x,y) = (v+\theta w/4^n,A(v+\theta w/4^n))$ (for some $n \in \mathbb{N}$) or $(x,y) = (v,0)$ if and only if $\|\tilde x(r)\|^2_{\MG{2}} \geq 9\epsilon_1^2$ and $\|\tilde y(r)\|_{\MG{2}}^2 \leq  (2^{-r} + \alpha/8)^2$. 
		
		To see this, first assume that $(\tilde x,\tilde y) \in \tilde \T$ are such that $\|\tilde x(r)\|_{\MG{2}}^2 \geq 9\epsilon_1^2$ and $\|\tilde y(r)\|_{\MG{2}}^2 \leq  (2^{-r} + \alpha/8)^2$. Then 
		\[
		\|x\|_{\MG{2}} \geq  \|\tilde x(r)\|_{\MG{2}} - \|x - \tilde x(r) \|_{\MG{2}} \geq 3\epsilon_1 - \epsilon_1 \geq 2\epsilon_1 
		\]
		and
		\begin{equation*}
			\|y \|_{\MG{2}} \leq \|\tilde y(r)\|_{\MG{2}} + \|y - \tilde y(r)\|_{\MG{2}} \leq 2^{-r} + \frac{\alpha}{8}  + 2^{-r} 
			\leq \left(\frac{1}{2} + \frac{1}{8}\right)\alpha < \alpha.
		\end{equation*}
		In particular, $(x,y) \notin \TT'$ (since any $(x,y) \in \TT'$ cannot have $\|y\|_{\MG{2}} < \alpha$ by Property  \ref{Prop:AwayFrom0}), $(x,y) \notin \{(\alpha A^{\dagger}(e_1),\alpha e_1),\dots,
		(\alpha A^{\dagger}(e_m),  \alpha e_m))\}$ for the same reason, and $(x,y) \notin \{(0,0)\}$ since $\|x\|_2 \geq 2\epsilon_1$. Thus $(x,y) \in \{(\frac{\theta}{4^{n } }w + \x, A(\frac{\theta}{4^{n}}w + \x))\, \vert \, n \in \mathbb{N}\} \cup  \{(\x,0)\}$ as claimed.
		
		For the converse, if $(x,y) \in \{(\frac{\theta}{4^{n } }w + \x, A(\frac{\theta}{4^{n}}w + \x))\, \vert \, n \in \mathbb{N}\} \cup  \{(\x,0)\}$ then either $\|x\|^2_{\MG{2}} = \|\theta w/4^n\|^2_{\MG{2}} + \|\x\|^2_{\MG{2}} $ or $\|x\|^2_{\MG{2}} = \|\x\|^2_{\MG{2}}$ (where we have used the fact that $\x$ and $w$ are orthogonal as $w \in \ker(A)$ and $\x \in \ker(A)^{\bot}$). In either case, $\|x\|_{\MG{2}} \geq \|v\|_{\MG{2}} = 4\epsilon_1$ and hence $\|\tilde {x}(r)\|_{\MG{2}} \geq \|x\|_{\MG{2}} - 2^{-r} \geq 4\epsilon_1 - 2^{-r} \geq 3\epsilon_1 > 2\epsilon_1$. Similarly, since $\x \in \ker(A)$ we must have $\|y\|_{\MG{2}} = \theta \|A(4^{-n}w)\|_{\MG{2}}$ or $y = 0$. In the first case, $\|y\|_2 \leq \alpha\|A\|^{-1}_{op}\|A\|_{op} 4^{-n} \|w\|_{\MG{2}}/8 \leq \alpha\|w\|_{\MG{2}}/8\leq \alpha/8$, whereas in the latter $y = 0$ and so $\|y\|_{\MG{2}} \leq \alpha/8$ holds trivially. Hence in any case
		$\|\tilde{y}(r)\|_{\MG{2}} \leq \|\tilde{y}(r) - y\|_{\MG{2}} + \|y\|_{\MG{2}} \leq 2^{-r}  + \alpha/8$, completing the proof of the claim.
		This ends the proof of the claim.
		
		In particular, this result immediately implies that \MG{instruction \ref{inst:vSearch} from `Approximate $v$'} terminates correctly; indeed, every $\T \in \Omega$ has at least one of $(v+\theta w/4^n,A(v+\theta w/4^n)) \in \T$ (for some $n \in \mathbb{N}$) or $(v,0) \in \T$. 
		To see the algorithm outputs a correct approximation to $v$, we use the converse result proven in the claim. More precisely, the claim implies that when the algorithm terminates we must have that $(\tilde x^j, \tilde y^j)$ corresponds to $(x^j,y^j)$ with $(x^j,y^j) = (v+\theta w/4^n,A(v+\theta w/4^n))$ \MG{for some $n \in \mathbb{N}$} or $(x^j,y^j) = (v,0)$. In either case, $\|y^j\|_2 \leq \alpha/8 $ (by the same argument used in the claim). Using the fact that $\theta w/4^n \in \orker{A}$ and $v \in \ker(A)$, we also have $x^j - A^{\dagger}(y^j) = v  + \theta w/4^n - A^{\dagger}(A(v  + \theta w/4^n)) = v + \theta w/4^n - \theta \MG{w}/4^n = v$ (in the case where $(x^j,y^j) = (v+\theta w/4^n,A(v+\theta w/4^n))$) and $x^j - A^{\dagger}(y) = v  - A^{\dagger}(A(v)) = v$ (in the case where $(x^j,y^j) = (v,0)$). \MG{Thus we have established that $x^j - A^{\dagger}(y^j)  = v$ and that $\|y^j\|_2 \leq \alpha/8$.}
		Therefore
		\MG{\begin{align*} \|v^k - v\|_{\MG{2}} &= \|\tilde x^j(k') - A^{\dagger}_{k'}(\tilde y^j(k')) - (x^j- A^{\dagger}(y^j))\|_{\MG{2}}\\ &\leq \|\tilde x^j(k')- x^j \|_{\MG{2}} + \|A^{\dagger}_{k'}(\tilde y^j(k')) -A^{\dagger}(\tilde y^j(k'))\|_{\MG{2}} + \|A^{\dagger}(\tilde y^j(k')) -A^{\dagger}(y^j)\|_{\MG{2}} \\&\leq 2^{-k'}+ \|A^{\dagger}_{k'} - A^{\dagger}\|_{op} \|\tilde y^j(k') \|_{\MG{2}}  + \|A^{\dagger}\|_{op}\|\tilde y^j(k') -y^j\|_{\MG{2}}\\
		&\leq 2^{-k'} + 2^{-{k'}} (2^{-{k'}} + \|y^j\|_{\MG{2}}) + \|A^{\dagger}\|_{op}2^{-{k'}}.
		\end{align*}}
		\MG{Now, since $\alpha \leq 1/4$ and by the definition of $k'$ from instruction \ref{inst:kPrimeDef} from `Approximate $v$', we get $2^{-k'} \leq \alpha 2^{-k}/2\leq 1/8$. Hence 
			$
				\|v^k - v\|_{\MG{2}} \leq \frac{\alpha2^{-k}}{2} \left(1+\frac{1}{8} + \|y^j\|_2 +  \|A^{\dagger}\|_{op}\right).
			$
		The definition of $\alpha$ implies that $\|A^{\dagger}\|_{op} \leq \beta^{-1} \leq 1/\alpha$ and as above we have $\|y^j\|_{\MG{2}} \leq \alpha/8 \leq 1/32$. Thus
		\begin{equation*}
			\|v^k - v\|_{\MG{2}} \leq \frac{\alpha2^{-k}}{2} \left(1+\frac{1}{8} + \frac{1}{32} +  \alpha^{-1}\right)= 2^{-k}\left(\frac{37\alpha}{64} + \frac{1}{2}\right) \leq 2^{-k}
		\end{equation*}
		where the final inequality again uses $\alpha \leq 1/4$.}

		We are now ready to present our algorithm $\Gamma: \tilde{\Omega} \times \R_+ 
		\to \NN_{m,N}^{\F}$. 
		For a given dyadic $\epsilon \in (2\epsilon_1,\infty)$, let $\epsilon_2 > 0$ be dyadic and such that 
		$2\epsilon_1 + 2\epsilon_2 \leq \epsilon$, and note that such an $\epsilon_2$ always exists because $\epsilon$ is strictly bigger then $2\epsilon_1$. Then, for any 
		$\tilde{\T} \in \tilde{\Omega}$, $\Gamma(\tilde{\T}, \epsilon)$ 
		produces the single layer neural network $A_{k}^{\dagger}(\cdot) + \frac{1}{2} \x^k$, where $k$ is chosen so that $2^{-k} \leq \min\{\epsilon_2,\alpha\}$ (this can be done recursively by iterating through natural numbers until this is satisfied) and where $A_{k}^{\dagger}$ and $\x^k$ are found using the pseudo codes for `Approximate $A^{\dagger}$' and `Approximate $\x$' respectively. In particular, such $A^{\dagger}_k$ and $\x^k$ must satisfy
		$\|A^{\dagger}_k - A^{\dagger} \|_{op} \leq \min\{\epsilon_2,\alpha\}$ and
		$\| \x^k - \x \|_{\MG{2}} \leq \epsilon_2$.

		It only remains to show that for \MG{all $\T \in \Omega$ and $\tilde \T \in \tilde \Omega$ corresponding to $\T$}there exists a neural network $\mathbf{N}_{opt}^{\M_1(\T)} \MG{\in \Xi(\T)}$ with
		\begin{align*}
			\sup_{y \in \bigcup_{\F} \M_2} \| \Gamma(\tilde{\T}, \epsilon)(y) - 
			\mathbf{N}_{opt}^{\M_1(\T)}(y) \|_{\MG{2}} \leq \epsilon,
		\end{align*}
		where 
		$
		\Gamma(\tilde{\T}, \epsilon)(y) = A_k^{\dagger}(y) + \frac{1}{2}  \x^k
		$ 	
		, and that the Lipschitz constant of  $\Gamma(\tilde{\T}, \epsilon)$ is uniformly bounded so that
		$L(\Gamma(\tilde{\T}, \epsilon))  \leq \frac{2}{\beta}$. 
		We start by proving the former. For a given $\T$, we let $\mathbf{N}_{opt}^{\M_1(\T)}$ be the neural network constructed in the argument in `\SecOmegaAndXiNonEmpty'.
		Notice that
		\begin{align*}
			\sup_{y \in \bigcup_{\F} \M_2} \| \Gamma(\tilde{\T}, \epsilon)(y) - 
			\mathbf{N}_{opt}^{\M_1(\T)}(y) \|_{\MG{2}} 
			=& \sup_{y \in \bigcup_{\F} \M_2} \| A_k^{\dagger}(y) + \frac{1}{2} \x^k - 
			\mathbf{N}_{opt}^{\M_1(\T)}(y) \|_{\MG{2}}  \\
			\leq& \sup_{y \in \bigcup_{\F} \M_2} \|A^{\dagger}(y) \!+\! \frac{1}{2} \x \! - \! \mathbf{N}_{opt}^{\M_1(\T)}(y) \|_{\MG{2}}
			+ \| A_k^{\dagger} - A^{\dagger} \|_{op} \\&+ \frac{1}{2} \| \x^k \! -\! \x \|_{\MG{2}},
		\end{align*}
		where the last inequality follows since 
		$A^{\dagger}_k(y) = A^{\dagger}(y) + A^{\dagger}_k(y) - A^{\dagger}(y)$ and $\|y \|_{\MG{2}} \leq 1$ when $y \in \bigcup_{\F} \M_2$. Now, we consider the term 
		$\|A^{\dagger}(y) + \frac{1}{2} \x -  \mathbf{N}^{\M_1(\T)}_{opt}(y) \|_{\MG{2}}$ and show that it must always be less than or equal to $2\epsilon_1$. 
		Indeed, we notice that if $(x,y) \in \iota_n^1$ for some $n \in \N$, then either
		$A^{\dagger}(y) = x = \mathbf{N}_{opt}^{\M_1(\T)}(y)$ in the case where 
		$(x,y) \in \TT \cup \{(0,0)\}$ (where we have used Property \ref{Prop:y=Ax}),
		or $A^{\dagger}(y) = \frac{\theta}{4^n}w$ and 
		$\mathbf{N}_{opt}^{\M_1(\T)}(y) = \frac{\theta}{4^n}w + \x \MG{= x}$, in the case where $(x,y) = (\frac{\theta}{4^n}w + \x,A(\frac{\theta}{4^n}w + \x))$. 
		In the first case we get that
		\begin{align*}
			\|A^{\dagger}(y) + \frac{1}{2} \x -  \mathbf{N}_{opt}^{\M_1(\T)}(y) \|_{\MG{2}}  
			= \| x +  \frac{1}{2} \x - x \|_{\MG{2}} = \frac{1}{2}\| \x \|_{\MG{2}}
			= \frac{1}{2}4\epsilon_1 = 2\epsilon_1, 
		\end{align*}
		and in the second case we get that 
		\begin{align*}
			\|A^{\dagger}(y) + \frac{1}{2} \x - \mathbf{N}_{opt}^{\M_1(\T)}(y) \|_{\MG{2}}  
			= \| \frac{\theta}{4^n}w +  \frac{1}{2} \x - \frac{\theta}{4^n}w  - \x \|_{\MG{2}} = \frac{1}{2}\| \x \|_{\MG{2}}
			= \frac{1}{2}4\epsilon_1 = 2\epsilon_1.
		\end{align*}
		On the other hand, if $(x,y) \in \iota^0$ we either have that
		$A^{\dagger}(y) = x = \mathbf{N}_{opt}^{\M_1(\T)}(y)$, in the case where $y \in \TT$, or 
		$A^{\dagger}(y) = 0$ and $\mathbf{N}_{opt}^{\M_1(\T)}(y) = \frac{1}{2} \x$ in the case where $(x,y) \in \{(0,0), (\x,0)\}$. In the first case we get the exact same bound as above, while in the second case we get that
		\begin{align*}
			\|A^{\dagger}(y) + \frac{1}{2} \x -  \mathbf{N}_{opt}^{\M_1(\T)}(y) \|_{\MG{2}}  = \| \frac{1}{2} \x - \frac{1}{2}\x \|_{\MG{2}} = 0 \leq 2\epsilon_1.
		\end{align*} 
		%	On the other hand, if $y \in A(\M_2)$ we either have that
		%	$A^{\dagger}(y) = x = \mathbf{N}_{opt}^{\M_1(\T)}(y)$, in the case where $y \in \pi_2(\TT)$, or 
		%	$A^{\dagger}(y) = 0$ and $\mathbf{N}_{opt}^{\M_1(\T)}(y) = \frac{1}{2} \x$, in the case where $y = 0$ or $y = \x$. In the first case we get the exact same bound as above, while in the second case we get that
		%	\begin{align*}
			%		\|A^{\dagger}(y) + \frac{1}{2} \x -  \mathbf{N}_{opt}^{\M_1(\T)}(y) \|_{\MG{2}}  = \| \frac{1}{2} \x - \frac{1}{2}\x \|_{\MG{2}} = 0.
			%	\end{align*} 
		Using this, and the fact that $\|A_k^{\dagger} - A^{\dagger} \|_{op} \leq \epsilon_2$, we arrive at the conclusion that
		\begin{align*}
			\sup_{y \in \bigcup_{\F} \M_2} \| \Gamma(\tilde{\T}, \epsilon)(y) 
			- \mathbf{N}_{opt}^{\M_1(\T)}(y) \|_{\MG{2}} 
			&\leq \sup_{y \in \bigcup_{\F} \M_2} \! \|A^{\dagger}(y) \! +\! \frac{1}{2} \x \!- 
			\! \mathbf{N}_{opt}^{\M_1(\T)}(y) \|_{\MG{2}} 
			+ \| A_k^{\dagger} \! -\!  A^{\dagger} \|_{op}  +\frac{1}{2} \| \x^k - \x \|_{\MG{2}}\\
			&\leq 2\epsilon_1 \! + \epsilon_2 + \epsilon_2 = 2\epsilon_1 + 2\epsilon_2 
			\leq \epsilon .
		\end{align*}
		Moreover, since $\Gamma(\tilde{\T}, \epsilon)$ is a neural network,
		we can conclude that $\Gamma(\tilde{\T}, \epsilon)$ produces a neural network which approximates $\mathbf{N}_{opt}^{\M_1(\T)}$ 
		to accuracy $\epsilon$ for all 
		$\tilde{\T} \in \tilde{\Omega}$ \MG{which correspond to $\T$}.
		Finally, we address the Lipschitz constant of $\mathbf{N}_{\mathcal{\tilde \T},\epsilon} = 
		\Gamma(\mathcal{\tilde{T}},\epsilon)$. We have
		\begin{align*}
			L(\mathbf{N}_{\mathcal{\tilde \T},\epsilon}) \phantom{lll}
			&=\sup\limits_{y \neq \hat y}
			\frac{\| \mathbf{N}_{\mathcal{\tilde \T},\epsilon}(y)-\mathbf{N}_{\mathcal{\tilde \T},\epsilon}(\hat y)\|_{\MG{2}}}{\|y - \hat y\|_{\MG{2}}} \\
			&\leq \sup\limits_{y \neq \hat y}
			\frac{\|A^{\dagger}_k \|_{op}\|y - \hat y \|_{\MG{2}}}{\|y -\hat y\|_{\MG{2}}}
			\leq \|\MG{A^{\dagger}} \|_{op} + 2^{-k} \leq \frac{1}{\beta} + \frac{1}{\beta} \leq 
			\frac{2}{\beta},
		\end{align*}
		since $2^{-k} \leq \alpha \leq \frac{1}{4}\min\{\frac{1}{\beta},\beta\} \leq \frac{1}{\beta}$. This completes the proof of part (iii).
	\end{proof}
	\bibliographystyle{abbrv}
	\bibliography{NNBarriersFOCM} 

@article{Felix,
	author = {Grohs, Philipp and Voigtlaender, Felix},
        journal = {Foundations of Computational Mathematics},
	number = {4},
	pages = {1085--1143},
	title = {Proof of the Theory-to-Practice Gap in Deep Learning via Sampling Complexity bounds for Neural Network Approximation Spaces},
	volume = {24},
	year = {2024},
}

@article{Hertling,
title = {{A Banach-Mazur computable but not Markov computable function on the computable real numbers}},
journal = {Annals of Pure and Applied Logic},
volume = {132},
number = {2},
pages = {227-246},
year = {2005},
issn = {0168-0072},
author = {Peter Hertling},
}

@misc{CRP,
	title={{On the consistent reasoning paradox of intelligence and optimal trust in {AI}: The power of 'I don't know'}}, 
	author={Alexander Bastounis and Paolo Campodonico and Mihaela van der Schaar and Ben Adcock and Anders C. Hansen},
	year={2024},
	eprint={2408.02357},
	archivePrefix={arXiv},
	primaryClass={cs.AI},
	url={https://arxiv.org/abs/2408.02357}, 
}

@article{Ziegler,
title = {Computability in linear algebra},
journal = {Theoretical Computer Science},
volume = {326},
number = {1},
pages = {187-211},
year = {2004},
author = {Martin Ziegler and Vasco Brattka}
}

@book{pour1989computability,
  title={Computability in Analysis and Physics},
  author={Pour-El, M.B. and Richards, J.I.},
  isbn={9783540500353},
  lccn={lc88024880},
  series={Perspectives in Mathematical Logic},
  url={https://books.google.co.uk/books?id=vNruAAAAMAAJ},
  year={1989},
  publisher={Springer Berlin Heidelberg}
}

@article{Haagerup,
 author = {Uffe Haagerup and Steen Thorbj{\o}rnsen},
 journal = {Annals of Mathematics},
 number = {2},
 pages = {711--775},
 publisher = {Annals of Mathematics},
 title = {A New Application of Random Matrices:
          $\text{Ext}(C_{\text{red}}^{\ast}(F_{2}))$
          Is Not a Group},
 volume = {162},
 year = {2005}
}

@article{Alex1,
	Author = {Bastounis, A and Cucker, F. and Hansen, A. C.},
	Title = {When can you trust feature selection? -- {I}: {A} condition-based analysis of {LASSO} and generalised hardness of approximation},
	Year = {2023},
	journal = {arXiv:2312.11425}
}

@article{Alex2,
	Author = {Bastounis, A and Cucker, F. and Hansen, A. C.},
	Title = {When can you trust feature selection? -- {II}: {O}n the effects of random data on condition in statistics and optimisation},
	Year = {2023},
	journal = {arXiv:2312.11429}
}

@article{Vasseur2019BlowUpST,
  title={Blow-Up Solutions to 3{D} {E}uler are Hydrodynamically Unstable},
  author={Alexis F. Vasseur and Misha M. Vishik},
  journal={Communications in Mathematical Physics},
  year={2019},
  volume={378},
  number = {1},
  pages={557--568},
  url={https://api.semanticscholar.org/CorpusID:201058517}
}

@article{Hou2,
	Author = {Chen, Jiajie and Hou, Thomas Y.},
	Title = {Stable nearly self-similar blowup of the 2{D} {B}oussinesq and 3{D} {E}uler equations with smooth data},
	Year = {2022},
	journal = {arXiv:2210.07191v2}
}

@article{Hou3,
	Author = {Chen, Jiajie and Hou, Thomas Y.},
	Title = {Stable nearly self-similar blowup of the 2{D} {B}oussinesq and 3{D} {E}uler equations with smooth data {II}: {R}igorous {N}umerics},
	Year = {2023},
	journal = {arXiv:2305.05660}
}

@article{hales2005proof,
  title={A proof of the {K}epler conjecture},
  author={Hales, Thomas C},
  journal={Annals of mathematics},
  pages={1065--1185},
   volume = {162},
    number = {3},
  year={2005},
}

@inproceedings{hales2017formal,
  title={A formal proof of the {K}epler conjecture},
  author={Hales, Thomas and Adams, Mark and Bauer, Gertrud and Dang, Tat Dat and Harrison, John and Le Truong, Hoang and Kaliszyk, Cezary and Magron, Victor and McLaughlin, Sean and Nguyen, Tat Thang and others},
  booktitle={Forum of Mathematics, Pi},
  volume={5},
  year={2017},
  organization={Cambridge University Press}
}

@article{fefferman1994_2,
	author = {C. Fefferman and L. Seco},
	journal = {Advances in Mathematics},
	number = {1},
	pages = {1--185},
	title = {On the {D}irac and {S}chwinger corrections to the ground-state energy of an atom},
	volume = {107},
	year = {1994}}

@incollection{fefferman1996interval,
  title={Interval arithmetic in quantum mechanics},
  author={Fefferman, Charles and Seco, Luis },
  booktitle={Applications of interval computations},
  pages={145--167},
  year={1996},
  publisher={Springer}
}

@article {Markov,
    AUTHOR = {Ce\u{\i}tin, G. S.},
     TITLE = {Algorithmic operators in constructive metric spaces},
   JOURNAL = {Trudy Mat. Inst. Steklov.},
  FJOURNAL = {Akademiya Nauk SSSR. Trudy Matematicheskogo Instituta imeni V.
              A. Steklova},
    VOLUME = {67},
      YEAR = {1962},
     PAGES = {295--361},
      ISSN = {0371-9685},
   MRCLASS = {02.72},
  MRNUMBER = {152426},
MRREVIEWER = {J.\ C.\ Shepherdson},
}

@inproceedings{Shoenfield,
    AUTHOR = {Kreisel, G. and Lacombe, D. and Shoenfield, J. R.},
     TITLE = {Partial recursive functionals and effective operations},
 BOOKTITLE = {Constructivity in mathematics: {P}roceedings of the colloquium
              held at {A}msterdam, 1957 (edited by {A}. {H}eyting)},
    SERIES = {Stud. Logic Found. Math.},
     PAGES = {290--297},
 PUBLISHER = {North-Holland, Amsterdam},
      YEAR = {1959},
   MRCLASS = {02.00},
  MRNUMBER = {108443},
MRREVIEWER = {M.\ O.\ Rabin},
}

@article{Cubitt2, 
title={Undecidability of the Spectral Gap}, 
volume={10}, DOI={10.1017/fmp.2021.15}, 
journal={Forum of Mathematics, Pi}, 
author={Cubitt, Toby and Perez-Garcia, David and Wolf, Michael M.}, 
year={2022}, 
pages={e14}
}

@article{Cubitt1,
		Author = {Cubitt, Toby S. and Perez-Garcia, David and Wolf, Michael M.},
	Journal = {Nature},
	Number = {7581},
	Pages = {207--211},
	Title = {Undecidability of the spectral gap},
	Ty = {JOUR},
	Url = {https://doi.org/10.1038/nature16059},
	Volume = {528},
	Year = {2015}
	}

@proceedings{AIM,
  editor    = "Fefferman, C. and  Hansen, A. C. and Jitomirskaya, S.",
  title     = "Computational mathematics in computer assisted proofs",
  series    = "American Institute of Mathematics Workshops",
  publisher = "American Institute of Mathematics",
 YEAR = {2022},
     PAGES = {1--7},
     note = {Available online at \url{https://aimath.org/pastworkshops/compproofsvrep.pdf}}
}

@article{fannjiang2020numerics,
  title={The numerics of phase retrieval},
  author={Fannjiang, Albert and Strohmer, Thomas},
  journal={Acta Numerica},
  volume={29},
  pages={125--228},
  year={2020}
}

@article{candes2013phaselift,
  title={Phaselift: Exact and stable signal recovery from magnitude measurements via convex programming},
  author={Cand{\`{e}}s, Emmanuel J and Strohmer, Thomas and Voroninski, Vladislav},
  journal={Communications on Pure and Applied Mathematics},
  volume={66},
  number={8},
  pages={1241--1274},
  year={2013},
  publisher={Wiley Online Library}
}

@inproceedings{yarotsky2018optimal,
	title={Optimal approximation of continuous functions by very deep {R}e{LU} networks},
	author={Yarotsky, Dmitry},
	booktitle={Conference on learning theory},
	pages={639--649},
	year={2018},
	organization={PMLR}
}

@article{petersen2018optimal,
	title={Optimal approximation of piecewise smooth functions using deep {R}e{LU} neural networks},
	author={Petersen, Philipp and Voigtlaender, Felix},
	journal={Neural Networks},
	volume={108},
	pages={296--330},
	year={2018},
	publisher={Elsevier}
}

@article{boelcskei2019optimal,
	title={Optimal approximation with sparsely connected deep neural networks},
	author={B\"olcskei, Helmut and Grohs, Philipp and Kutyniok, Gitta and Petersen, Philipp},
	journal={SIAM Journal on Mathematics of Data Science},
	volume={1},
	number={1},
	pages={8--45},
	year={2019},
	publisher={SIAM}
}

@article{Arora1_Goedel,
title = "Probabilistic checking of proofs: A new characterization of {NP}",
author = "Sanjeev Arora and Shmuel Safra",
year = "1998",
volume = "45",
pages = "70--122",
journal = "Journal of the ACM",
publisher = "Association for Computing Machinery (ACM)",
number = "1",
}

@article{Bast_SIAM_NEWS,
  title={From Global to Local: Getting More from Compressed Sensing},
  author={Bastounis, Alexander and Adcock, Ben and Hansen, Anders C.},
  journal={SIAM News},
  volume={50},
  number={8},
  pages={1--4},
  year={2017}
}

@article{Bast_2017,
author = {Bastounis, Alexander and Hansen, Anders C.},
title = {On the Absence of Uniform Recovery in Many Real-World Applications of Compressed Sensing and the Restricted Isometry Property and Nullspace Property in Levels},
journal = {SIAM Journal on Imaging Sciences},
volume = {10},
number = {1},
pages = {335-371},
year = {2017}
}

@article{colbrook2021computingSIREV,
  title={Computing spectral measures of self-adjoint operators},
  author={Colbrook, Matthew J. and Horning, Andrew and Townsend, Alex},
  journal={SIAM Review},
  volume={63},
  number={3},
  pages={489--524},
  year={2021},
  publisher={SIAM}
}

@book{MarkovModel,
	title={Computation and Automata},
	author={Salomaa, A. and Cambridge University Press and Rota, G.C. and Doran, B. and Lam, T.Y. and Flajolet, P. and Ismail, M. and Lutwak, E.},
	isbn={9780521302456},
	lccn={84017571},
	series={EBL-Schweitzer},
	url={https://books.google.co.uk/books?id=IblDi626fBAC},
	year={1985},
	publisher={Cambridge University Press}
}

@article{Ben_Artzi2022,
	title={Computing scattering resonances},
	author={Ben-Artzi, Jonathan and Marletta, Marco and R{\"o}sler, Frank},
	journal={Journal of the European Mathematical Society},
	volume={25},
	number={9},
	pages={3633--3663},
	year={2022}
}

@article{godel1931formal,
  title={{{\"U}ber formal unentscheidbare S{\"a}tze der Principia Mathematica und verwandter Systeme I}},
  author={G{\"o}del, Kurt},
  journal={Monatshefte f{\"u}r mathematik und physik},
  volume={38},
  number={1},
  pages={173--198},
  year={1931},
  publisher={Springer}
}

@article{Smale_Weinberger, 
author = {Niyogi, P. and Smale, S. and Weinberger, S.}, 
title = {A Topological View of Unsupervised Learning from Noisy Data}, 
year = {2011}, 
issue_date = {May 2011}, 
publisher = {Society for Industrial and Applied Mathematics}, 
volume = {40}, 
number = {3}, 
issn = {0097-5397}, 
journal = {SIAM Journal on Computing}, 
pages = {646--663}
  }

@ARTICLE{Boyer_2016,
  author={J. {Bigot} and C. {Boyer} and P. {Weiss}},
  journal={IEEE Transactions on Information Theory}, 
  title={An Analysis of Block Sampling Strategies in Compressed Sensing}, 
  year={2016},
  volume={62},
  number={4},
  pages={2125-2139},
  doi={10.1109/TIT.2016.2524628}}

@article{BOYER_ACHA_2019,
title = "Compressed sensing with structured sparsity and structured acquisition",
journal = "Applied and Computational Harmonic Analysis",
volume = "46",
number = "2",
pages = "312 - 350",
year = "2019",
issn = "1063-5203",
doi = "https://doi.org/10.1016/j.acha.2017.05.005",
url = "http://www.sciencedirect.com/science/article/pii/S1063520317300489",
author = "C. Boyer and J. Bigot and P. Weiss",
}

@inproceedings{DezFaFr-17,
	Author = {S. Moosavi-Dezfooli and A. Fawzi and O. Fawzi and P. Frossard},
	Booktitle = {IEEE Conference on computer vision and pattern recognition},
	Month = {July},
	Pages = {86-94},
	Title = {Universal Adversarial Perturbations},
	Year = {2017}}

@inproceedings{SzZ-14,
	Author = {Christian Szegedy and Wojciech Zaremba and Ilya Sutskever and Joan Bruna and Dumitru Erhan and Ian J. Goodfellow and Rob Fergus},
	Booktitle = {International Conference on Learning Representations},
	Title = {Intriguing properties of neural networks},
	Year = {2014}}

@article{hammernik2018learning,
  title={Learning a variational network for reconstruction of accelerated {MRI} data},
  author={Hammernik, Kerstin and Klatzer, Teresa and Kobler, Erich and Recht, Michael P and Sodickson, Daniel K and Pock, Thomas and Knoll, Florian},
  journal={Magnetic Resonance in Medicine},
  volume={79},
  number={6},
  pages={3055--3071},
  year={2018},
  publisher={Wiley Online Library}
}

@inproceedings{tyukin2020adversarial,
  title={On adversarial examples and stealth attacks in artificial intelligence systems},
  author={Tyukin, I.Y. and Higham, D.J. and Gorban, A.N.},
  booktitle={2020 International Joint Conference on Neural Networks (IJCNN)},
  pages={1--6},
  year={2020},
  organization={IEEE}
}

@article{Choi2,
 author={Choi, Charles}, 
  journal={IEEE Spectrum}, 
  title={Some {AI} Systems May Be Impossible to Compute}, 
  year={2022},
  volume={March}
}

@article{THGW21,
	author = {Tyukin, Ivan Y and Higham, Desmond J and Bastounis, Alexander and Woldegeorgis, Eliyas and Gorban, Alexander N},
	title = {The feasibility and inevitability of stealth attacks},
	journal = {IMA Journal of Applied Mathematics},
	volume = {89},
	number = {1},
	pages = {44-84},
	year = {2023},
	month = {10},
	issn = {0272-4960},
	doi = {10.1093/imamat/hxad027},
	url = {https://doi.org/10.1093/imamat/hxad027},
	eprint = {https://academic.oup.com/imamat/article-pdf/89/1/44/58326072/hxad027.pdf},
}

@article{webb2021spectra,
  title={Spectra of {J}acobi operators via connection coefficient matrices},
  author={Webb, Marcus and Olver, Sheehan},
  journal={Communications in Mathematical Physics},
  volume={382},
  number={2},
  pages={657--707},
  year={2021},
  publisher={Springer}
}

@book{Weinberger, 
author = {Weinberger, Shmuel}, 
title = {Computers, Rigidity, and Moduli: The Large-Scale Fractal Geometry of Riemannian Moduli Space}, 
year = {2004}, 
isbn = {0691118892}, 
publisher = {Princeton University Press}, 
address = {USA} 
}

@article{adcock2016generalized,
  title={Generalized sampling and infinite-dimensional compressed sensing},
  author={Adcock, Ben and Hansen, Anders C},
  journal={Foundations of Computational Mathematics},
  volume={16},
  number={5},
  pages={1263--1323},
  year={2016},
  publisher={Springer}
}

@article{heaven2019deep,
  title={Why deep-learning {AI}s are so easy to fool},
  author={Heaven, Douglas},
  journal={Nature},
  volume={574},
  number={7777},
  pages={163--166},
  year={2019},
  publisher={Nature Publishing Group}
}

@article{Choi,
 author={Choi, Charles}, 
  journal={IEEE Spectrum}, 
  title={7 Revealing Ways {AI}s Fail}, 
  year={2021},
  volume={September}
}

@article{devore2020neural,
  title={Neural network approximation},
  author={DeVore, Ronald and Hanin, Boris and Petrova, Guergana},
  journal={Acta Numerica},
  volume={30},
  pages={327--444},
  year={2021},
  publisher={Cambridge University Press}
}

@ARTICLE{Baraniuk2020,
  author={G. {Ongie} and A. {Jalal} and C. A. {Metzler} and R. G. {Baraniuk} and A. G. {Dimakis} and R. {Willett}},
  journal={IEEE Journal on Selected Areas in Information Theory}, 
  title={Deep Learning Techniques for Inverse Problems in Imaging}, 
  year={2020},
  volume={1},
  number={1},
  pages={39-56},
  doi={10.1109/JSAIT.2020.2991563}}

@ARTICLE{Genzel,
  author={Genzel, Martin and Macdonald, Jan and M{\"a}rz, Maximilian},
  journal={IEEE Transactions on Pattern Analysis and Machine Intelligence}, 
  title={Solving Inverse Problems With Deep Neural Networks - Robustness Included}, 
  year={2022},
  volume={},
  number={},
  pages={1-1},
  doi={10.1109/TPAMI.2022.3148324}}

@article{Hallucination_NatureM2,
	Author = {Hoffman, David P. and Slavitt, Isaac and Fitzpatrick, Casey A.},
	Da = {2021/02/01},
	Id = {Hoffman2021},
	Isbn = {1548-7105},
	Journal = {Nature Methods},
	Number = {2},
	Pages = {131--132},
	Title = {The promise and peril of deep learning in microscopy},
	Volume = {18},
	Year = {2021},
	}

@article{finlayson2019adversarial,
  title={Adversarial attacks on medical machine learning},
  author={Finlayson, Samuel G and Bowers, John D and Ito, Joichi and Zittrain, Jonathan L and Beam, Andrew L and Kohane, Isaac S},
  journal={Science},
  volume={363},
  number={6433},
  pages={1287--1289},
  year={2019},
  publisher={American Association for the Advancement of Science}
}

@article{arridge2019solving,
  title={Solving inverse problems using data-driven models},
  author={Arridge, Simon and Maass, Peter and {\"O}ktem, Ozan and Sch{\"o}nlieb, Carola-Bibiane},
  journal={Acta Numerica},
  volume={28},
  pages={1--174},
  year={2019},
  publisher={Cambridge University Press}
}

@article{mccann2017convolutional,
  title={Convolutional neural networks for inverse problems in imaging: {A} review},
  author={McCann, Michael T and Jin, Kyong Hwan and Unser, Michael},
  journal={IEEE Signal Process Magazine},
  volume={34},
  number={6},
  pages={85--95},
  year={2017},
  publisher={IEEE}
}

@article{jin17,
  title={Deep convolutional neural network for inverse problems in imaging},
  author={Jin, Kyong Hwan and McCann, Michael T and Froustey, Emmanuel and Unser, Michael},
  journal={IEEE Transactions on Image Processing},
  volume={26},
  number={9},
  pages={4509--4522},
  year={2017},
  publisher={IEEE}
}

@article{donoho2006compressed,
  title={Compressed sensing},
  author={Donoho, David L},
  journal={IEEE Transactions on Information Theory},
  volume={52},
  number={4},
  pages={1289--1306},
  year={2006},
  publisher={IEEE}
}

@article{candes2006robust,
  title={Robust uncertainty principles: Exact signal reconstruction from highly incomplete frequency information},
  author={Cand{\`e}s, Emmanuel J and Romberg, Justin and Tao, Terence},
  journal={IEEE Transactions on Information Theory},
  volume={52},
  number={2},
  pages={489--509},
  year={2006},
  publisher={IEEE}
}

@article{adcock2020gap,
  title={The gap between theory and practice in function approximation with deep neural networks},
  author={Adcock, Ben and Dexter, Nick},
  journal={SIAM Journal on Mathematics of Data Science},
  volume={3},
  number={2},
  pages={624--655},
  year={2021},
  publisher={SIAM}
}

@article{gottschling2023troublesome,
author = {Gottschling, Nina M. and Antun, Vegard and Hansen, Anders C. and Adcock, Ben},
title = {The Troublesome Kernel: On Hallucinations, No Free Lunches, and the Accuracy-Stability Tradeoff in Inverse Problems},
journal = {SIAM Review},
volume = {67},
number = {1},
pages = {73-104},
year = {2025}
}

@article {antun2020instabilities,
	author = {Antun, Vegard and Renna, Francesco and Poon, Clarice and Adcock, Ben and Hansen, Anders C.},
	title = {On instabilities of deep learning in image reconstruction and the potential costs of {AI}},
	volume = {117},
	number = {48},
	pages = {30088--30095},
	year = {2020},
	publisher = {National Academy of Sciences},
	issn = {0027--8424},
	journal = {Proceedings of the National Academy of Sciences}
}

@book{AdcockHansenBook,
	author = {Adcock, B. and Hansen, A. C.},
	isbn = {9781108421614},
	publisher = {Cambridge University Press},
	title = {Compressive Imaging: Structure, Sampling, Learning},
	url = {https://books.google.co.uk/books?id=cAMOzgEACAAJ},
	year = {2021}}

@article{colbrook_spec_meas,
title={Computing Spectral Measures and Spectral Types},
  author={Colbrook, Matthew J.},
  journal={Communications in Mathematical Physics},
  volume={384},
  number={1},
  pages={433--501},
  year={2021},
  DOI={10.1007/s00220-021-04072-4}
}

@article{Hansen2016ComplexityII,
  title={Complexity Issues in Computing Spectra, Pseudospectra and Resolvents},
  author={Anders C. Hansen and Olavi Nevanlinna},
  journal={Banach Center Publications},
  year={2016},
  volume={112},
  pages={171-194}
}

@article{SCI,
      title={Computing Spectra -- {O}n the Solvability Complexity Index Hierarchy and Towers of Algorithms}, 
      author={Jonathan Ben-Artzi and Matthew J. Colbrook and Anders C. Hansen and Olavi Nevanlinna and Markus Seidel},
      year={2020},
      eprint={1508.03280},
      journal={arXiv:1508.03280v5},
      primaryClass={cs.CC}
}

@article {smale_question,
    AUTHOR = {Smale, Steve},
     TITLE = {On the efficiency of algorithms of analysis},
   JOURNAL = {Bulletin of the American Mathematical Society},
    VOLUME = {13},
      YEAR = {1985},
    NUMBER = {2},
     PAGES = {87--121}
}

@article {Smale2,
    AUTHOR = {Smale, S.},
     TITLE = {The fundamental theorem of algebra and complexity theory},
   JOURNAL = {Bulletin of the American Mathematical Society},
    VOLUME = {4},
      YEAR = {1981},
    NUMBER = {1},
     PAGES = {1--36},
      ISSN = {0273-0979},
     CODEN = {BAMOAD},
   MRCLASS = {65H05 (12D10 58-02 58C30 68C25 90C30)},
  MRNUMBER = {590817 (83i:65044)},
}

@article{CRAS,
	author = {Jonathan Ben-Artzi and Anders C. Hansen and Olavi Nevanlinna and Markus Seidel},
	doi = {http://dx.doi.org/10.1016/j.crma.2015.08.002},
	issn = {1631-073X},
	journal = {Comptes Rendus Mathematique},
	number = {10},
	pages = {931 - 936},
	title = {New barriers in complexity theory: On the solvability complexity index and the towers of algorithms},
	url = {http://www.sciencedirect.com/science/article/pii/S1631073X15002137},
	volume = {353},
	year = {2015}}

@article{Matt1,
author = {Colbrook, Matthew and Hansen, Anders C.},
year = {2023},
  pages={4639--4718},
  volume = {25}, 
 number = {12},
title = {The foundations of spectral computations via the {S}olvability {C}omplexity {I}ndex hierarchy},
journal = {Journal of the European Mathematical Society}
}

@article{ben2021computing,
  title={Computing the sound of the sea in a seashell},
  author={Ben-Artzi, Jonathan and Marletta, Marco and R{\"o}sler, Frank},
  journal={Foundations of Computational Mathematics},
  pages={1--35},
  volume = {22}, 
 number = {3},
 year={2021},
  publisher={Springer}
}

@article{Sudan_Overview_2009,
author = {Sudan, Madhu},
title = {Probabilistically Checkable Proofs},
year = {2009},
issue_date = {March 2009},
publisher = {Association for Computing Machinery},
address = {New York, NY, USA},
volume = {52},
number = {3},
journal = {Communications of the ACM},
pages = {76-84},
numpages = {9}
}

@article{Sudan, 
author = {Bellare, Mihir and Goldreich, Oded and Sudan, Madhu}, 
title = {Free Bits, {PCP}s, and Nonapproximability -- Towards Tight Results}, 
year = {1998}, 
issue_date = {June 1998}, 
publisher = {Society for Industrial and Applied Mathematics}, 
volume = {27}, 
number = {3}, 
journal = {SIAM Journal on Computing}, 
pages = {804-915}, 
numpages = {112}
}

@article{Hastad_Acta,
	Author = {H{\aa}stad, Johan},
	Da = {1999/03/01},
	Doi = {10.1007/BF02392825},
	Id = {H{\aa}stad1999},
	Isbn = {1871-2509},
	Journal = {Acta Mathematica},
	Number = {1},
	Pages = {105--142},
	Title = {Clique is hard to approximate within $n^{1-\epsilon}$},
	Ty = {JOUR},
	Volume = {182},
	Year = {1999}
}

@article{Juditsky_2011,
  author    = {Anatoli B. Juditsky and
               Fatma Kilin{\c{c}}{-}Karzan and
               Arkadi Nemirovski},
  title     = {Verifiable conditions of $\ell_{1}$-recovery
               for sparse signals with sign restrictions},
  journal   = {Mathematical Programming},
  volume    = {127},
  number    = {1},
  pages     = {89--122},
  year      = {2011}
  }

@article{Juditsky_2012,
author = {Anatoli Juditsky and Fatma Kilin{\c{c}}{-}Karzan and Arkadi Nemirovski and Boris Polyak},
title = {{Accuracy guaranties for $\ell_{1}$ recovery of block-sparse signals}},
volume = {40},
journal = {The Annals of Statistics},
number = {6},
publisher = {Institute of Mathematical Statistics},
pages = {3077--3107},
year = {2012}
}

@article{Fefferman_Klartag2,
author = {Charles Fefferman and Bo'az Klartag},
title = {{Fitting a $C^m$-Smooth Function to Data II}},
volume = {25},
journal = {Revista Matematica Iberoamericana},
number = {1},
publisher = {Real Sociedad Matem{\'a}tica Espa{\~n}ola},
pages = {49--273},
year = {2009},
doi = {rmi/1236864106},
URL = {https://doi.org/}
}

@Article{Fefferman_Klartag,
 Author = {Charles L. {Fefferman} and Bo'az {Klartag}},
 Title = {{Fitting a \(C^m\)-smooth function to data. I}},
 FJournal = {{Annals of Mathematics. Second Series}},
 Journal = {{Annals of Mathematics}},
 ISSN = {0003-486X; 1939-8980/e},
 Volume = {169},
 Number = {1},
 Pages = {315--346},
 Year = {2009},
 Publisher = {Princeton University, Mathematics Department, Princeton, NJ},
 Language = {English},
 MSC2010 = {41A05 65D15},
 Zbl = {1175.41001}
}

@book{realRAM,
  title={Computational Geometry: An Introduction},
  author={Preparata, F.P. and Shamos, M.I.},
  isbn={9781461210986},
  series={Monographs in Computer Science},
  url={https://books.google.co.uk/books?id=\_p3eBwAAQBAJ},
  year={2012},
  publisher={Springer New York}
}

@article{von_Neumann, 
author = {von Neumann, John}, 
title = {First Draft of a Report on the {EDVAC}}, 
year = {1993}, 
issue_date = {October 1993}, 
publisher = {IEEE Educational Activities Department}, 
address = {USA}, 
volume = {15}, 
number = {4}, 
issn = {1058-6180}, 
journal = {IEEE Annals of the History of Computing}, 
pages = {27-75}, 
numpages = {49} }

@incollection {Smale_McMullen,
    AUTHOR = {Smale, S.},
     TITLE = {{The work of {C}urtis {T} {M}c{M}ullen}},
 BOOKTITLE = { Proceedings of the International Congress of Mathematicians I, Berlin},
    SERIES = {Doc. Math. J. DMV},
     PAGES = {127--132},
      YEAR = {1998},
}

@book{bishop1967foundations,
  title={Foundations of Constructive Analysis},
  author={Bishop, E.},
  lccn={67022952},
  series={McGraw-Hill Series in higher mathematics},
  url={https://books.google.co.uk/books?id=o2mmAAAAIAAJ},
  year={1967},
  publisher={McGraw-Hill}
}

@article{Cucker_Smale97, 
author = {Cucker, Felipe and Smale, Steve}, 
title = {Complexity Estimates Depending on Condition and Round-off Error}, 
year = {1999}, 
issue_date = {Jan. 1999}, 
publisher = {Association for Computing Machinery}, 
address = {New York, NY, USA}, 
volume = {46}, 
number = {1}, 
doi = {10.1145/300515.300519}, 
journal = {Journal of the ACM}, 
pages = {113-184}, 
 }

@article{Lovasz_JACM_96, 
author = {Feige, Uriel and Goldwasser, Shafi and Lov\'{a}sz, Laszlo and Safra, Shmuel and Szegedy, Mario}, 
title = {Interactive Proofs and the Hardness of Approximating Cliques}, 
year = {1996}, 
issue_date = {March 1996}, 
publisher = {Association for Computing Machinery}, 
address = {New York, NY, USA}, 
volume = {43}, number = {2}, 
issn = {0004-5411}, 
url = {https://doi.org/10.1145/226643.226652}, 
doi = {10.1145/226643.226652}, 
journal = {Journal of the ACM}, 
pages = {268-292}, numpages = {25}, 
 }

@article{Arora_JACM_98_2, 
author = {Arora, Sanjeev and Lund, Carsten and Motwani, Rajeev and Sudan, Madhu and Szegedy, Mario}, 
title = {Proof Verification and the Hardness of Approximation Problems}, 
year = {1998}, 
issue_date = {May 1998}, 
publisher = {Association for Computing Machinery}, 
address = {New York, NY, USA}, 
volume = {45}, 
number = {3}, 
 journal = {Journal of the ACM}, 
month = may, 
pages = {501-555}, 
numpages = {55}
}

@book{Arora2007,
	author = {Sanjeev Arora and Boaz Barak},
	publisher = {Princeton University Press},
	title = {Computational Complexity - A Modern Approach},
	year = {2009}}

@book{Ko1991ComplexityTO,
	author = {K. Ko},
	booktitle = {Progress in Theoretical Computer Science},
	title = {Complexity Theory of Real Functions},
  publisher = {Birkhauser},
	year = {1991}}

@book{lovasz1987algorithmic,
	author = {Lovasz, L.},
	isbn = {9780898712032},
	lccn = {86061532},
	publisher = {Society for Industrial and Applied Mathematics},
	series = {CBMS-NSF Regional Conference Series in Applied Mathematics},
	title = {An Algorithmic Theory of Numbers, Graphs and Convexity},
	url = {https://books.google.co.uk/books?id=sJ3mBHTU55QC},
	year = {1987}}

@article{Chambolle_Alg,
	author = {Chambolle, Antonin},
	doi = {10.1023/B:JMIV.0000011325.36760.1e},
	issn = {1573-7683},
	journal = {Journal of Mathematical Imaging and Vision},
	number = {1},
	pages = {89--97},
	title = {An Algorithm for Total Variation Minimization and Applications},
	url = {http://dx.doi.org/10.1023/B:JMIV.0000011325.36760.1e},
	volume = {20},
	year = {2004}}

@article{Chambolle_Alg2,
	acmid = {1969036},
	address = {Norwell, MA, USA},
	author = {Chambolle, Antonin and Pock, Thomas},
	doi = {10.1007/s10851-010-0251-1},
	issn = {0924-9907},
	issue_date = {May 2011},
	journal = {Journal of Mathematical Imaging and Vision},
	month = may,
	number = {1},
	numpages = {26},
	pages = {120--145},
	publisher = {Kluwer Academic Publishers},
	title = {A First-Order Primal-Dual Algorithm for Convex Problems with Applications to Imaging},
	url = {http://dx.doi.org/10.1007/s10851-010-0251-1},
	volume = {40},
	year = {2011}}

@article{BSS_Machine,
	author = {Blum, Lenore and Shub, Mike and Smale, Steve},
	coden = {BAMOAD},
	doi = {10.1090/S0273-0979-1989-15750-9},
	fjournal = {American Mathematical Society. Bulletin. New Series},
	issn = {0273-0979},
	journal = {Bulletin of the American Mathematical Society},
	mrclass = {68Q10 (03D15 03D75 03D80 03F60 68Q15 68Q25)},
	mrnumber = {974426 (90a:68022)},
	mrreviewer = {John Michael Robson},
	number = {1},
	pages = {1--46},
	title = {On a theory of computation and complexity over the real numbers: {NP}-completeness, recursive functions and universal machines},
	url = {http://dx.doi.org/10.1090/S0273-0979-1989-15750-9},
	volume = {21},
	year = {1989}}

@article{DeVore,
	author = {Cohen, Albert and Dahmen, Wolfgang and DeVore, Ronald},
	doi = {10.1090/S0894-0347-08-00610-3},
	fjournal = {Journal of the American Mathematical Society},
	issn = {0894-0347},
	journal = {Journal of the American Mathematical Society},
	mrclass = {94A12 (15B52)},
	mrnumber = {2449058},
	mrreviewer = {Luoqing Li},
	number = {1},
	pages = {211--231},
	title = {Compressed sensing and best {$k$}-term approximation},
	url = {http://dx.doi.org/10.1090/S0894-0347-08-00610-3},
	volume = {22},
	year = {2009}}

@incollection{Smale_Acta_Numerica,
	address = {Cambridge},
	author = {Smale, Steve},
	booktitle = {Acta numerica, 1997},
	doi = {10.1017/S0962492900002774},
	mrclass = {65Y20 (68Q15 68Q25)},
	mrnumber = {1489262 (99d:65385)},
	mrreviewer = {Klaus Meer},
	pages = {523--551},
	publisher = {Cambridge Univ. Press},
	series = {Acta Numer.},
	title = {Complexity theory and numerical analysis},
	url = {http://dx.doi.org/10.1017/S0962492900002774},
	volume = {6},
	year = {1997}}

@article{Turing_Machine,
	author = {Turing, A. M.},
	doi = {10.1112/plms/s2-42.1.230},
	fjournal = {Proceedings of the London Mathematical Society},
	issn = {0024-6115},
	journal = {Proceedings of the London Mathematical Society},
	mrclass = {Contributed Item},
	mrnumber = {1577030},
	number = {1},
	pages = {230},
	title = {On {C}omputable {N}umbers, with an {A}pplication to the {E}ntscheidungsproblem},
	url = {http://dx.doi.org/10.1112/plms/s2-42.1.230},
	volume = {S2-42},
	year = {1936}}

@article{McMullen1,
	author = {McMullen, Curt},
	fjournal = {Annals of Mathematics. Second Series},
	journal = {Annals of Mathematics},
	number = {3},
	pages = {467-493},
	title = {Families of rational maps and iterative root-finding algorithms},
	volume = {125},
	year = {1987}}

@article{McMullen2,
	author = {McMullen, Curt},
	coden = {INVMBH},
	doi = {10.1007/BF01389368},
	fjournal = {Inventiones Mathematicae},
	issn = {0020-9910},
	journal = {Inventiones Mathematicae},
	mrclass = {58F12 (30D05)},
	mrnumber = {922801 (89d:58075)},
	mrreviewer = {J. S. Birman},
	number = {2},
	pages = {259--272},
	title = {Braiding of the attractor and the failure of iterative algorithms},
	url = {http://dx.doi.org/10.1007/BF01389368},
	volume = {91},
	year = {1988}}

@article{Doyle_McMullen,
	author = {Doyle, Peter and McMullen, Curt},
	coden = {ACMAA8},
	doi = {10.1007/BF02392735},
	fjournal = {Acta Mathematica},
	issn = {0001-5962},
	journal = {Acta Mathematica},
	mrclass = {12E12 (20D06 20F19)},
	mrnumber = {1032073 (91d:12004)},
	mrreviewer = {Doru {\c{S}}tef{\u{a}}nescu},
	number = {3-4},
	pages = {151--180},
	title = {Solving the quintic by iteration},
	url = {http://dx.doi.org/10.1007/BF02392735},
	volume = {163},
	year = {1989}}

@article{Breaking,
	address = {Cambridge, UK},
	author = {B. Adcock and A. C. Hansen and C. Poon and B. Roman},
	day = {001},
	doi = {10.1017/fms.2016.32},
	journal = {Forum of Mathematics, Sigma},
	month = {001},
	pages = {1-84},
	publisher = {Cambridge University Press},
	title = {BREAKING THE COHERENCE BARRIER: A NEW THEORY FOR COMPRESSED SENSING},
	url = {https://www.cambridge.org/core/article/div-class-title-breaking-the-coherence-barrier-a-new-theory-for-compressed-sensing-div/455E5F506912B78E9647CD7A7488530B},
	volume = {5},
	year = {2017}}

@article{Hansen_JAMS,
	author = {Hansen, Anders C.},
	doi = {10.1090/S0894-0347-2010-00676-5},
	fjournal = {Journal of the American Mathematical Society},
	issn = {0894-0347},
	journal = {Journal of the American Mathematical Society},
	mrclass = {47A10 (46N40 47A75 65J10)},
	mrnumber = {2726600 (2012a:47009)},
	mrreviewer = {A. B{\"o}ttcher},
	number = {1},
	pages = {81--124},
	title = {On the solvability complexity index, the {$n$}-pseudospectrum and approximations of spectra of operators},
	url = {http://dx.doi.org/10.1090/S0894-0347-2010-00676-5},
	volume = {24},
	year = {2011}}

@article{pinkus,
author = {Pinkus, Allan},
title = {Approximation theory of the {MLP} model in neural networks},
journal = {Acta Numerica},
volume = {8},
publisher = {Cambridge University Press},
year = {1999},
}

@article{comp,
    title={ The extended {S}male's 9th problem -- {O}n computational barriers and paradoxes in estimation, regularisation, computer-assisted proofs and learning. },
    author = {Bastounis, A. and Hansen, A. C. and {Vla\v{c}i\'{c}}, V.},
    eprint={2110.15734},
      journal={arXiv:2110.15734},
    year={2021}
}

@book{meshfree,
author = {Fasshauer, Gregory F.},
title = {Meshfree Approximation Methods with MATLAB},
year = {2007},
isbn = {9789812706348},
publisher = {World Scientific Publishing Co., Inc.},
address = {USA},
}

@article{Matt2,
         Author = {Matthew J. Colbrook and Vegard Antun and Anders C. Hansen},
	Doi = {10.1073/pnas.2107151119},
	Eprint = {https://www.pnas.org/doi/pdf/10.1073/pnas.2107151119},
	Journal = {Proceedings of the National Academy of Sciences},
	Number = {12},
	Pages = {e2107151119},
	Title = {The difficulty of computing stable and accurate neural networks: {O}n the barriers of deep learning and {S}male's 18th problem},
	Url = {https://www.pnas.org/doi/abs/10.1073/pnas.2107151119},
	Volume = {119},
	Year = {2022}
	}

@book{Complexity_and_real_comp,
author = {Blum, Lenore and Cucker, Felipe and Shub, Michael and Smale, Steve},
title = {Complexity and Real Computation},
year = {1997},
isbn = {0387982817},
publisher = {Springer-Verlag},
address = {Berlin, Heidelberg}
}

@article{Turing,
    author = {Turing, A. M.},
    title = "{I.-Computing machinery and intelligence}",
    journal = {Mind},
    volume = {LIX},
    number = {236},
    pages = {433-460},
    year = {1950},
    month = {10},
    issn = {0026-4423},
    doi = {10.1093/mind/LIX.236.433},
    url = {https://doi.org/10.1093/mind/LIX.236.433},
    eprint = {https://academic.oup.com/mind/article-pdf/LIX/236/433/30123314/lix-236-433.pdf},
}

@INPROCEEDINGS{res1,
  author={He, Kaiming and Zhang, Xiangyu and Ren, Shaoqing and Sun, Jian},
  booktitle={2016 IEEE Conference on Computer Vision and Pattern Recognition (CVPR)}, 
  title={Deep Residual Learning for Image Recognition}, 
  year={2016},
  volume={},
  number={},
  pages={770-778},
  doi={10.1109/CVPR.2016.90}
  }

@inproceedings{res2,
author = {Sinha, Ayan and Lee, Justin and Li, Shuai and Barbastathis, George},
ournal = {Digital Holography and Three-Dimensional Imaging},
year = {2017},
publisher = {Optica Publishing Group},
title = {Solving inverse problems using residual neural networks},
doi = {10.1364/DH.2017.W1A.3},
url = {https://opg.optica.org/abstract.cfm?URI=DH-2017-W1A.3}
}
	
\end{document}